\documentclass[11pt]{amsart}
\usepackage[utf8]{inputenc}
\usepackage{fancyhdr}
\usepackage{multicol}

\usepackage{mathtools}
\usepackage{enumitem}
\usepackage{hyperref}
\usepackage[american]{babel}
\usepackage{amssymb}
\usepackage{geometry}
\newtheorem{theorem}{Theorem}[section]
\newtheorem{lemma}[theorem]{Lemma}
\newtheorem{sublemma}[theorem]{Sublemma}
\newtheorem{corollary}[theorem]{Corollary}

\newtheorem{definition}[theorem]{Definition}
\newtheorem{proposition}[theorem]{Proposition}
\newtheorem{remark}[theorem]{Remark}
\newtheorem{condition}[theorem]{Condition}

\newlist{enumroman}{enumerate}{1}
\setlist[enumroman,1]{label=(\roman*)}
\newcommand{\htop}{h_{\tiny{\mbox{top}}}}
\newcommand{\AAA}{{\mathcal{A}}}
\newcommand{\BB}{{\mathcal{B}}}
\newcommand{\CC}{{\mathcal{C}}}
\newcommand{\DD}{{\mathcal{D}}}
\newcommand{\EE}{{\mathcal{E}}}
\newcommand{\HH}{{\mathcal{H}}}

\newcommand{\OO}{{\mathcal{O}}}
\newcommand{\KK}{{\mathcal{K}}}
\newcommand{\LL}{{\mathcal{L}}}
\newcommand{\MM}{{\mathcal{M}}}
\newcommand{\MMM}{{\mathbb{M}}}
\newcommand{\NN}{{\mathcal{N}}}
\newcommand{\QQ}{{\mathcal{Q}}}
\newcommand{\RR}{{\mathcal{R}}}
\newcommand{\VV}{{\mathcal{V}}}
\newcommand{\WW}{{\mathcal{W}}}
\newcommand{\complex}{{\mathbb{C}}}
\newcommand{\SSS}{{\mathbb{S}}}
\newcommand{\real}{{\mathbb{R}}}
\newcommand{\integer}{{\mathbb{Z}}}
\newcommand{\naturall}{{\mathbb{N}}}
\newcommand{\FFF}{{\mathbb{F}}}
\newcommand{\supp}{{\mathrm{supp}\,}}
\newcommand{\tr}{{\mathrm{tr}\,}}
\newcommand{\dd}{{\mathrm{d}}}
\newcommand{\dimm}{{\mathrm{dim}\,}}
\newcommand{\diamm}{{\mathrm{diam}\,}}
\newcommand{\divv}{{\mathrm{div}\,}}
\newcommand{\voll}{{\mathrm{vol}\,}}

\newcommand{\intt}{{\mathrm{int}\,}}
\newcommand{\D}{{\mathrm{D}}}

\newcommand{\id}{\mathrm{id}}
\newcommand{\iunit}{i}
\newcommand{\Op}{\mathrm{Op}}

\begin{document}
\title[Horocycle Averages  and Transfer Operators]{Horocycle Averages on Closed Manifolds
		and Transfer Operators}
	\author{Alexander Adam \and Viviane Baladi}
\address{\phantom{noaddress}}
	\email{a.adam.sci@protonmail.com}
\address{  
Sorbonne Universit\'e and Universit\'e Paris Cit\'e, CNRS,Laboratoire de Probabilit\'es, Statistique et Mod\'elisation,
F-75005 Paris, France}
\email{baladi@lpsm.paris}
	\date{\today}
	\thanks{AA thanks  Centre Henri Lebesgue and CIRM  Luminy for their warm hospitality
and the K. and A. Wallenberg foundation
		for  invitations to Lund.
We thank   Y. Guedes Bonthonneau, L. Flaminio, G.~Forni, S. Gou\"ezel, C.~Guillarmou,   M. J\'ez\'equel, C. Liverani,
and  S.  Shen for helpful  discussions and the referee for sharp remarks encouraging us to improve
the exposition.    Most of this work was done while
AA was working at IMJ-PRG, Paris and the Max Planck Institute for Mathematics, Bonn. This research is supported
by the European Research Council (ERC) under the European Union's Horizon 2020 research and innovation programme (grant agreement No 787304).}
	
\begin{abstract}
We study the ergodic integrals  of the horocycle flows $h_\rho$ of $C^r$ codimension one mixing Anosov flows.
In  dimension three, for any suitably bunched $C^3$ contact Anosov flow  with orientable strong-stable distribution $E_-$, we show $|\frac 1T\int_0^T\varphi\circ h_{\rho}(x)\dd\rho
	-\mu(\varphi)|\le \frac {C} {T^{\epsilon}} \|\varphi\|_{C^{3}}$ for some $\epsilon >0$, with $\mu$ the invariant
measure of $h_\rho$. We thereby implement  the toy model program of Giulietti--Liverani \cite{giul_liv_2017} in the natural setting of geodesic flows in variable negative curvature, where nontrivial resonances exist.
\end{abstract}
	
\maketitle


\section{Introduction}

 Anosov introduced a class of
$C^2$ flows
$g_{\alpha}\colon M\to M$, now bearing his name \cite{Anosov_1969}, on  closed (i.e. compact and boundaryless) orientable  manifolds\footnote{All manifolds are implicitly endowed with a Riemannian metric.} $M$ of  dimension $d\ge 3$.  We  focus on topologically mixing Anosov flows.
A special class of  mixing Anosov flows are those preserving a contact structure. Geodesic flows on the unit tangent bundle of a closed  manifold with (possibly variable) negative sectional curvature are well-studied classes of contact Anosov flows.

Every Anosov flow $g_\alpha$  admits a strong stable foliation, tangent 
to a vector bundle denoted $E_-$. 
If this foliation is orientable and has dimension $d_-$ equal to one, and if $g_\alpha$
is mixing, one associates
with  $g_{\alpha}$ another flow, the \textit{horocycle} flow $h_{\rho}\colon M\rightarrow M$, such that for every $x\in M$ the  trajectory $h_\real(x)$ is  a strong stable leaf (defined up to speed reparametrisation). 
 Horocycle flows were first introduced in the case of  geodesic Anosov flows, \cite[p.84]{Marcus_1977}, \cite{Hirsch_1975}.
In a setting more general than ours (with $d_-\ge 1$),  Bowen and Marcus  \cite{Bowen_1977} then proved that the
 horocycle flow is uniquely ergodic and minimal. Its invariant probability measure $\mu$ 
(related to, but distinct from, the
measure of maximal entropy of $g_\alpha$, see Remark~\ref{MME}) plays an important role below. 

Since the horocycle flow is induced by the Anosov flow, there exists  $\tau(\rho,\alpha,x)>0$ such that
\[g_{\alpha}\circ h_{\rho}(x)=h_{\tau(\rho,\alpha,x)}\circ g_{\alpha}(x)\, ,
\,\, \forall x \in M\, ,\,\,\forall \alpha\, ,\rho \in \real_+  .\]
We call  $\tau(\rho,\alpha,x)$ the \textit{renormalisation time}.
 In the setting of unit speed geodesic flows on compact (or more generally, finite volume) surfaces of constant negative 
curvature, renormalisation has been used effectively in the work of Flaminio and Forni \cite{Flaminio_2003} to study the \emph{horocycle integrals}
\[
\gamma_{x}(\varphi,T)\coloneqq \int_0^T\varphi\circ h_{\rho}(x)\dd\rho\, ,
\quad x\in M\, ,\quad T>0\, ,
\]
for $\varphi:M\to \real$ in Sobolev spaces of positive order.
Flaminio and Forni found that the speed of convergence 
of $\gamma_{x}(\varphi,T)/T$ to $\mu(\varphi)$ as $T\to\infty$ is controlled by 
\textit{invariant distributions} under the push-forward of the horocyclic vector field. These distributions are also eigendistributions  under the push-forward of the geodesic vector field, and
the  eigenvalues give the powers of $T$ appearing in the expansion of $\gamma_{x}(\varphi,T)/T-\mu(\varphi)$.

Their approach inspired Giulietti and Liverani \cite{giul_liv_2017} 
to study a toy model, replacing the Anosov flow with a hyperbolic diffeomorphism $F$, using the renormalisation dynamics as a key to study $\gamma_{x}(\varphi,T)$.
Letting $\htop$ be the  topological entropy of $F$, they show analogously (for the corresponding invariant measure $\mu$) that the speed of convergence to zero of $\gamma_{x}(\varphi,T)/T-\mu(\varphi)$
is controlled by eigenvalues in the annulus
$1<|z|<e^{\htop}$  (and the corresponding eigendistributions) of a weighted transfer operator of $F$.
 Unfortunately, in  the setting of \cite{giul_liv_2017},
 there are in fact no eigenvalues in the annulus $1<|z|<e^{\htop}$, see \cite{balkill}.
The approach of Giulietti and Liverani has been applied
successfully in the meantime by Faure--Gou\"ezel--Lanneau \cite{FGL} 
to the linear flow in the stable direction of  a two-dimensional linear pseudo-Anosov map, and by Butterley--Simonelli \cite{BuSi} to parabolic
flows on (3-dimensional)  Heisenberg nilmanifolds  which are renormalized by  partially hyperbolic automorphisms (circle extensions
of Anosov diffeomorphisms). In both these algebraic applications, nontrivial eigenvalues are present.

Giulietti and Liverani conjectured that a similar expansion exists for more general (non algebraic) Anosov flows than in \cite{Flaminio_2003}, e.g. for the geodesic flow 
of a  surface with variable negative  curvature \cite[Conjecture 2.12]{giul_liv_2017}. More precisely,
letting $\htop$ be the topological entropy of the time-one map $g_1$,
 we expect that there exists $\delta >0$ such that, for smooth enough observables $\varphi$, the following expansion holds\footnote{For Anosov flows, $\htop>0$, see  \cite{Anosov_1969}. For geodesic flows on finite volume negatively curved surfaces, $\htop=1$.} (analogously to  \cite{Flaminio_2003}, \cite{giul_liv_2017})
\begin{align}\label{eq:intro:I}
	\gamma_{x}(\varphi,T)=T \int \varphi\dd\mu+ \sum_{\delta < \Re \lambda < \htop}T^{\frac{\Re \lambda}{\htop}} \tilde c_\lambda(T,x)\OO_\lambda(\varphi)+
	{\EE}_{T,x}(\varphi).
\end{align}
In the above formula, $\EE_{T,x} = O (T^{\frac{\delta}{\htop}})$, uniformly in $x$,
the $\OO_\lambda$ are generalised eigendistributions associated to the eigenvalue $\lambda$ for the adjoint (or dual) of the generator  of a weighted transfer operator (see \eqref{1tr}) acting on an anisotropic Banach space, the real parameter
$\delta$ is an upper bound on the essential spectral bound of the generator, and   $\tilde c_\lambda(T,x)\in \complex$ 
satisfies $\sup_{x,T} |\tilde c_\lambda(T,x)| |\log T|^{-J_\lambda}<\infty$, where $J_\lambda\ge 0$ is the size of the largest   Jordan block of $\lambda$.

The main result of this work, Theorem~\ref{t:decomposition},  provides an asymptotic expansion \eqref{eq:intro:I} for $C^r$ time reparametrisations of the
unit speed horocycle flow of codimension one topologically mixing $C^r$ Anosov flows, if $r>2$ and  the distribution
$E_-$ is $C^{r-1}$, under an essential spectral gap  
condition ($\lambda^{s,t,p}_{min}<\htop$), and a weak Dolgopyat condition on the resolvent (Condition~\ref{cnd:wA}).
As a consequence, we get power-law convergence of the ergodic averages (Corollary~\ref{c:cor0}).
In Proposition~\ref{p:Anosov},
we show that the conditions of  Theorem~\ref{t:decomposition} hold for
 $C^3$ contact Anosov flows in dimension three,
with  orientable strong stable bundle $E_-$,  under 
the following bunching assumption:
Recalling that $d_-=1$, define 
$$
\lambda_+=\lim_{\alpha \to \infty} \frac  {\sup_x \log \|\D  g_{-\alpha}(x)|_{E_- }\|}{\alpha}\, \, , \, \,
\lambda_-=-\lim_{\alpha \to \infty} \frac {  \sup_x \log \|\D  g_{\alpha}(x)|_{E_-}\|}{\alpha}
\, , \, \,  \hat \varpi:=2\frac{ \lambda_-}{\lambda_+}\in (0,2]\,  .
$$
The bunching condition  is 
\begin{equation}\label{smalleta}
\hat\varpi >\frac{8}{5}\, .
\end{equation}
For constant negative curvature geodesic flows, we have $\hat \varpi =2$.
Assumption \eqref{smalleta} 
thus holds  for  geodesic flows with variable strictly negative curvature close enough to
a constant, but the reader is warned that it does not apply to generic  three-dimensional
contact Anosov flows.

For  compact  surfaces of constant negative curvature, Randol \cite{Randol_1974} proved  that there exist eigenvalues of the Laplacian arbitrarily close to $1$. This provides examples for which the expansion of Flaminio--Forni \cite{Flaminio_2003}, and thus the expansion in Theorem~\ref{t:decomposition}, 
is not reduced to $T \int \varphi\dd\mu+
	{\EE}_{T,x}(\varphi)$.
 
As in the work of Giulietti and Liverani \cite{giul_liv_2017}, the key idea to study $\gamma_{x}(\varphi,T)$ is to introduce the weighted semigroup of transfer operators, with generator $X+V$, defined by
\begin{equation}
\label{1tr}
\LL_{\alpha,V}: W_p^{s,t,q}(M)\rightarrow W_p^{s,t,q}(M), \quad \LL_{\alpha,V}\varphi= \phi_\alpha \cdot(\varphi\circ g_{-\alpha})\, , \,\, \phi_\alpha(x)=e^{\int_0^\alpha V\circ g_\beta(x)\dd \beta}\, ,\quad \alpha\ge 0 \, ,
\end{equation}
where the potential is $V=-\partial_\alpha \partial_{\rho}\tau(0,0,\cdot)$ (so that $\phi_\alpha=\partial_{\rho}\tau(0,-\alpha,\cdot)$),
and where $W_p^{s,t,q}(M)$ is an anisotropic Banach space with  regularity parameters $s<0<q\le  t<r-1+s$ and $p\in (1,\infty)$.  
In the case of the unit speed parametrisation of  $h_{\rho}$, we shall see that 
$\phi_\alpha=\det\D g_{-\alpha}{}|_{E^-}$ is just the Jacobian along the strong stable distribution at a negative time $-\alpha$, and
$V= \divv ( X|_{E_-})$.

\smallskip

The paper is organised as follows:  The transfer operator $\LL_{\alpha,V}$  is defined in Section~\ref{s:to} (for more general potentials).
  The new anisotropic Banach spaces $W_p^{s,t,q}(M)$ are constructed in Section~\ref{s:BS}
after introducing admissible cones for $C^r$ Anosov flows in Section~\ref{s:geo} (if $p=2$ we get
Hilbert spaces). 
These spaces are a flow analogue to the spaces constructed in  \cite{Baladi2007} to study hyperbolic diffeomorphisms. Anisotropic Banach spaces are now a standard tool  for  hyperbolic dynamics (see e.g. \cite{blank_ruelleperronfrobenius_2002, Liverani_2004, Baladi2012,   Giulietti_2013, Baladi_2018,  Gouezel_2008,   Tsujii_2016}). Although we do not study here the dynamical determinant or zeta function associated to the transfer operator $\LL_{\alpha,V}$, we believe that the  spaces introduced in the present work are well suited for 
this purpose (see \cite{Bbook}). Guedes Bonthonneau and Lefeuvre very recently \cite{GL} applied a (microlocal) flow
 implementation  of the spaces from  \cite{Baladi2007} to study some dynamical and geometric problems.

In Section~\ref{s:prop}, we establish  properties of the transfer operator semigroup, its generator $X+V$ and the resolvent $\RR_z$ (see \eqref{eq:res}). 
Most of these results do not require the contact assumption. Among those are norm bounds  yielding a Lasota--York inequality for the resolvent (Theorem~\ref{t:res:lasota}). Then, in Corollary~\ref{c:res:spectrum}, we obtain a strip in the spectrum of the generator containing at most countable eigenvalues of finite multiplicity. 
Proposition~\ref{cnd:A} puts the weak Dolgopyat Condition~\ref{cnd:wA} in
more standard form. These tools are used in Section~\ref{s:horo}
to show the above-mentioned main results,  
Theorem~\ref{t:decomposition} and Proposition~\ref{p:Anosov}. 
(Our proofs highlight
sufficient conditions for intrinsicness of
resonances and portability of Dolgopyat bounds on the resolvent when navigating
between different Banach spaces.)
Finally, Appendix~\ref{intparts} contains  (elementary) integration by parts lemmas, adapted
 from \cite{Baladi2007},  Appendix~\ref{fragrec}
recalls  the fragmentation/reconstitution lemmas 
from \cite{Baladi2007}, and Appendix~\ref{mollapp} is devoted to interpolation and mollifiers.

\smallskip 

We end this introduction with some remarks:

\smallskip
(a)
The conjecture that the distributions $\OO_v$  in the
expansion \eqref{eq:intro:I} are fixed by the (adjoint) of the horocycle flow, which was the starting point in \cite{Flaminio_2003}, remains  open for general codimension one mixing
Anosov flows (see  \cite[Remark~2.10]{giul_liv_2017}). For smooth contact Anosov flows with $d=3$, invariance  was proved by Faure--Guillarmou \cite{Faure_2017_2}.

\smallskip
(b) The anisotropic Banach spaces $W^{s,t,q}_p(M)$ in this paper are based 
on those in \cite{Baladi2007}. We could also define spaces $\BB^{s,t}(M)$
based on those in
\cite{Baladi2008} (or \cite[Chapter 5]{Bbook}). We expect that the following variational upper bound may be obtained\footnote{Here, $h_{\tilde \mu}$ is the entropy of an ergodic $g_1$ invariant probability measure, and $\chi_{\tilde \mu}(A)$  is the largest Lyapunov exponent of 
$A$.}  for the essential spectrum of the semigroup $\LL_{\alpha, V}$ on $\BB^{s,t}(M)$:
\begin{equation}\label{conj}\lambda_{\min}^{s,t}(X,V):= \sup_{{\tilde \mu}}
\Bigl \{h_{\tilde \mu}(g_1) + \chi_{\tilde \mu}\left (\frac{\phi_1}{\det (Dg_1|_{E_+})}\right ) +
\max\bigl \{t \chi_{\tilde \mu}(Dg_1|_{E_-}), |s| \chi_{\tilde \mu}(Dg_{-1}|_{E_+} )\bigr \}
\Bigr \}\, .
\end{equation} 
The above is in general better (even in the volume preserving case)  than the bound $\lambda_{\min}^{s,t,p}(X,V)$ we obtain
in Corollary~\ref{c:res:spectrum} (see \eqref{eq:res:lasota:A1}). Since $\lambda_{\min}^{0,0}(X,V)=\htop$, the essential spectral gap condition
$\lambda_{\min}^{s,t}<\htop$ would thus hold 
for $\BB^{s,t}(M)$ for arbitrarily small $s<0$ and $t>0$ (so that the assumptions of Proposition~\ref{p:Anosov} could be weakened accordingly, and $s'$ could be taken arbitrarily close to $0$ in
\eqref{smalls'}).  However, the scale $\BB^{s,t}(M)$ is  more messy\footnote{See also the caveat in \cite[Remark 5.18]{Bbook} regarding the lack of validity of \eqref{trick}.} to define, it is not
an interpolation scale, it does not include a Hilbert space,
and showing \eqref{conj} would
require a  thermodynamic analysis of the
sums over subcovers in the proof of 
Lemma~\ref{l:to:bound}. To keep the paper short, we restrict to the scale $W^{s,t,q}_p(M)$.

\smallskip
(c) The renormalisation time $\tau(\rho,\alpha,\cdot)$ inherits the smoothness of the invariant bundle $E_-$, which is only H\"older
in general. We add the extra assumption  that $E_-$ is smooth enough and that an essential spectral gap holds in Theorem~\ref{t:decomposition} (and Lemma~\ref{l:invmeas}), and we give settings where this is satisfied in Proposition~\ref{p:Anosov}. To work with anisotropic spaces with higher regularity (depending only on $r$), one could lift the dynamics to the Grassmannian \cite{Gouezel_2008, giul_liv_2017}. 
We have chosen to avoid the cumbersome corresponding technicalities for the sake of readability.

\smallskip
(d) Finally, we mention two directions of future research: First, the expansions of Flaminio and Forni \cite{Flaminio_2003} (or Faure--Tsujii  \cite{Faure_2013, FT3}, see
also \cite{CG}) are not limited to finite sets  of eigenvalues. Our methods do not currently allow to go beyond the smallest $\delta$
such that $\Sigma_\delta=\sigma (X+V)|_{W_p^{s,t,q}(M)} \cap \{ \Re z > \delta\}$ is finite ($\delta=1/2$ for \cite{Flaminio_2003}).
Second, even if  the Dolgopyat Condition~\ref{cnd:wA} holds for some $\delta<0$, we cannot improve the remainder  due to the term  with $\|\varphi\|_0$  in  
Lemma~\ref{l:horo:decomp:local}. Although it is hoped that this term is spurious, 
an analogous error term is present  in \cite[Thm 1.5]{Flaminio_2003} or \cite[Thm 2.8]{giul_liv_2017}.

\section{The Transfer Operators and the Banach Spaces}
\label{sec2}
\subsection{Transfer Operators Associated to a Flow $g_\alpha$ and Weight $\phi_\alpha$. The Generator $X+V$}\label{s:to}
In the entire paper, $M$ is a compact, boundaryless, connected, orientable, smooth  manifold of dimension $d\ge 3$, and  $r>1$ is fixed, while $g_\alpha\colon M\rightarrow M$, $\alpha\in\real$, is
a $C^r$ {\it Anosov flow} on $M$. By definition, 
  there is a $Dg_\alpha$-invariant splitting of the tangent space  
$TM=E_-\oplus E_+\oplus E_0$
of the tangent space such that for some  $C_*\ge 1$ and $0<\theta<1$, we have
\begin{equation}\label{eq:split:ineq'}	
		\|\D  g_\alpha{v}\|\le C_* \theta^\alpha\|v\| \, , \forall v\in E_-\, , \,\,
		\| \D  g_{-\alpha}{v}\| \le C_* \theta^{\alpha}\|v\| \, ,\forall v\in E_+ \, ,\,\forall \alpha \ge 0\, ,
\end{equation}		
while $E_0=\langle X \rangle$, where the $C^{r-1}$ vector field $X$ is the generator of the flow
defined by
\begin{align}\label{eq:flow:gen}
	 X \coloneqq\partial_{\alpha}g_{-\alpha}{}|_{\alpha=0}\, .
\end{align}
The {\it (strong) stable and unstable distributions} $E_-$ and $E_+$ 
are H\"older.  For $x\in M$, we split
$T_xM=E_{-,x}\oplus E_{+,x}\oplus E_{0,x}$.
The cotangent space $T^*M$ (the dual space of $TM$) is split analogously
\begin{align*}
	T^*M=E_-^*\oplus E_+^*\oplus E_0^*\, ,\,\, T^*_xM=E_{-,x}^*\oplus E_{+,x}^*\oplus E_{0,x}^*,\quad x\in M\, .
\end{align*} 
The splitting above is $(\D g_{\alpha})^{\tr}$-invariant and, up to taking
larger $C_*$, we have 
\begin{align}\label{eq:split:ineq}
\begin{cases}
 C_*^{-1} \|\xi\|\le\|(\D g_{-\alpha})^{\tr}\xi\|
\le C_*\|\xi\| \, , & \forall \xi\in E_0^* \, ,\,\,\forall \alpha \ge 0\, ,\\
\|(\D g_{-\sigma \alpha})^{\tr}\xi\|
\le C_* \theta^\alpha \|\xi\| \, , & \forall \xi \in E_\sigma^* \, ,\,
\sigma=\pm
 \, ,\,\, \forall \alpha \ge 0\, .
\end{cases}
\end{align}
The dimensions of the spaces $E_{\sigma,x}$ do not depend on $x$, and we set
$d_-:= \dimm E_{-}=\dimm E^*_{-}$ and $d_+:= \dimm E_{+}=\dimm E^*_{+}=d-1-d_-$.
Fixing a   {\it potential} $V\in C^{r-1}(M,\real)$,
we introduce the $\phi_\alpha$-weighted transfer operators
\begin{align}\label{eq:to}
	\LL_{\alpha,V}(\varphi):= \phi_\alpha\cdot(\varphi\circ g_{-\alpha})\, ,\quad\alpha\ge 0\, , \,\,  \varphi\in C^{r-1}(M)\, ,
\end{align}
where
\begin{align*}
\phi_\alpha(x)\coloneqq \exp(\int_0^{\alpha} V\circ g_{-\beta}(x)\dd\beta)\, ,\quad
\mbox{i.e. }
	V=\partial_{\alpha}\phi_{\alpha}{}|_{\alpha=0^+}\, .
\end{align*}
For an ``integrability''  parameter $p\in (1,\infty)$, and 
suitable anisotropic ``regularity'' parameters $s$, $t$, and $q$ (see \eqref{regul-sqt}), we will construct Banach spaces $W_p^{s,t,q}(M)$, containing $C^{r-1}(M)$ as a dense subspace,  on which the  operators $\LL_{\alpha,V}$  extend continuously to form a strongly continuous semigroup (Lemma~\ref{l:to:C0}). 
In particular, for all $\varphi\in {C^{r-1}(M)}$ 
\[
\partial_\alpha\LL_{\alpha,V}\varphi{}|_{\alpha=0^+}=X\varphi+V\varphi\in W_p^{s,t,q}(M)\, . 
\]
 The generator  of the semigroup is $X+V$, we denote by $\RR_z$  its resolvent
\begin{align}
\label{eq:res}
	\RR_z\varphi= (z-V-X)^{-1}\varphi\, ,
\quad z\not\in\sigma (X+V)|_{W_p^{s,t,q}},\quad\varphi\in W_p^{s,t,q}(M) \, ,
\end{align}
 where
$\sigma (X+V)|_{\BB}$ denotes
the spectrum of the operator $X+V$  on $\BB$. 
Theorem~\ref{t:res:lasota} will provide a Lasota--Yorke inequality for  $\RR_z$ for large 
$\Re z$. This gives a vertical  strip in the complex plane in which $\sigma (X+V)|_{W_p^{s,t,q}(M)}$ contains only isolated eigenvalues of finite multiplicity (Corollary~\ref{c:res:spectrum}).

\subsection{Cone Ensembles.  The Atlas $\AAA$.
Cone Hyperbolicity. Admissible Cones for $g_\alpha$}\label{s:geo}

A {\it cone} is a nonempty convex set $\CC\subset \real^d$ such that $\lambda \xi\in \CC$ for all $\xi\in \CC$
and $\lambda \in \real$. We say that a cone $\CC$ is $d'$-dimensional if $d'\ge 1$ is
the maximal dimension of a linear subset of $\CC$.
A cone $\CC$ is \textit{compactly included} in another cone $\CC'$, denoted by $\CC\Subset \CC'$,
if  $\bar \CC\subseteq \intt \CC'\cup\{0\}$. Two cones $\CC$ and  $\CC'$ are \textit{transversal} if  $\CC\cap \CC'=\{0\}$.

We identify $T^\ast M$ with $\real^d$, and for any $d'\ge 1$ we denote the norm
of $\xi\in \real^{d'}$ by $|\xi|=(\sum_j \xi_j^2)^{1/2}$. For   $\xi\in T^*_xM$, write $\xi=\xi^-+\xi^++\xi^0$, where $\xi^\sigma\in E^*_{\sigma,x}$ for $\sigma\in\{\pm, 0\}$. For $\gamma >0$,  define two transversal closed  cone fields on\footnote{These cones  have non-empty interior while \cite[Proposition 17.4.4]{Katok1995} uses ``flat" cones included in $E_+^*\oplus E_-^*$.}
$T^*M$,
of respective dimensions $d_-$ and $d_+$,   by
\begin{align}
\label{eq:cone}
		\CC_-^\gamma(x)\coloneqq\{\xi \mid
\, \max\{|\xi^+|,|\xi^0|\}\le \gamma |\xi^-| \}
\, ,\,  \,\,
		\CC_+^\gamma(x)\coloneqq\{\xi\mid\, \max\{|\xi^-|,|\xi^0|\}\le \gamma |\xi^+|\}\, ,
\end{align}
and define a one-dimensional  closed  cone field on $T^* M$ by
\begin{align*}
		\CC_0^\gamma(x):=\{\xi\in T^\ast_xM \mid
\, \max\{|\xi^-|,|\xi^+|\}\le \gamma |\xi^0| \}
\, .
\end{align*}
 For $\sigma\in \{\pm, 0\}$ we have $E^*_{\sigma,x}\subset \CC_\sigma^{\gamma}(x)$
and,  if $\gamma'>\gamma$, then $
\CC_\sigma^\gamma(x)\Subset \CC_\sigma^{\gamma'}(x)$.
Moreover $T^*_x M\subset \cup_\sigma \CC_\sigma^\gamma(x)$, if $\gamma \ge 1$
(any line through the origin must cross one side of the unit cube in
$\real^3$), while
$\CC_\sigma^{\gamma}(x)$ and $\CC_\tau^{\gamma}(x)$ are transversal
if $\sigma\ne \tau$ and $\gamma <1$.
Last but not least, 
the lemma below 
is the key to construct admissible cones:\footnote{No such property holds for $\CC^\gamma_0(x)$. The cones in \eqref{eq:cone} are strictly expanding and contracting,
respectively, and this is not true for $\CC^\gamma_0(x)$.}

\begin{lemma}\label{l:cone:inclusion}
Let $C_*\in[1, \infty)$ and $\theta\in (0,1)$ be the constants from \eqref{eq:split:ineq}. Then for any $\gamma, \gamma'\in (0,1)$
 and all $\alpha>0$  such that $C^2_*\theta^\alpha\gamma<\gamma'$, we have,
recalling
$\CC^-_\gamma(x)$ and $\CC^+_\gamma(x)$ from \eqref{eq:cone},
	$$
(\D g_{-\alpha})^{\tr}\CC_\gamma^-(x)\Subset \CC_{\gamma'}^-(g_{\alpha}(x))\, 
,
\, \, \, (\D g_{\alpha})^{\tr}\CC_{\gamma}^+(x)\Subset \CC_{\gamma'}^+(g_{-\alpha}(x))\, ,
\,\forall x \in M\, .
$$
\end{lemma}

\begin{proof}
	We show the first claim. The proof of the second claim  is analogous. Let $\xi=\xi^-+\xi^++\xi^0\in \CC_\gamma^-(x)$. We estimate 
	\begin{align*}
	\max\{	|(\D g_{-\alpha})^{\tr}_x\xi^+|,
|(\D g_{-\alpha})^{\tr}_x\xi^0|\}&
\le C_*\max\{|\xi^+|,|\xi^0|\}\le C_*\gamma|\xi^-|\le C^2_*\theta^\alpha\gamma|(\D g_{-\alpha})^{\tr}\xi^-|\, .
	\end{align*}
	 Since $\CC_{\gamma'-\epsilon}^-(g_{\alpha}(x))\Subset \CC_{\gamma'}^-(g_{\alpha}(x))$ for all $\epsilon\in (0,\gamma')$, we conclude.
\end{proof}

Next, we adapt to flows the cone ensembles for hyperbolic diffeomorphisms from \cite{Baladi2007, Bbook}:
\begin{definition}[Cone ensembles $\Theta$ for flows. Coverings $\widetilde \Phi$]
A cone ensemble of $\real^d$, with $d=d_-+d_++1$, $d_\pm\ge 1$, is a pair $\Theta=(\CC,\Phi)$,
where $\CC=(\CC_-,\CC_+,\CC_0)$ is a triplet of pairwise
transversal closed cones with nonempty interiors, of respective dimensions $d_-$, $d_+$, and one,
while $\Phi=(\Phi_-,\Phi_+,\Phi_0)$, where each $\Phi_\sigma$
 is a $C^\infty$
map from the unit sphere
$\SSS^d$ to $[0,1]$, such that 
$$
\Phi_-+\Phi_++\Phi_0\equiv 1 \, ,\, \quad \Phi_\sigma|_{\CC_\sigma\cap \SSS^d}\equiv 1\, , \,\,\sigma \in\{\pm, 0\} \, .
$$
In addition, we require that $\CC_0=\{\xi \mid |\xi|\le \gamma_0 |\xi_d|\}$
for some finite $\gamma_0$.

For two cone ensembles $\Theta$ and $\Theta'$ of $\real^d$, we say that\footnote{For our purposes, the second
 condition could be replaced by the weaker condition $\real^d \setminus (\CC_-\cup \CC_+)\Subset \CC'_0\cup \CC'_-$.}
$\Theta'<\Theta$ if
$$
\real^d \setminus (\CC_+\cup \CC_0)\Subset \CC'_-\,\,
\mbox{ and } \,\,\real^d \setminus \CC_+\Subset \CC'_0\cup \CC'_-\, .
$$	

Finally, for a cone ensemble $\Theta$, we say that a triplet $\widetilde \Phi=(\widetilde \Phi_+,\widetilde \Phi_-,\widetilde \Phi_0)$ is a covering of $\Theta$
if each ${\widetilde\Phi}_{\sigma}:\SSS^d \to [0,1]$ is  $C^\infty$, with 
 ${\widetilde \Phi}_{\sigma}|_{\supp\Phi_\sigma}\equiv 1$.
\end{definition}

\begin{definition}[Cone hyperbolicity]
	Let $K\subset\real^d$ be compact with nonempty interior, and let $\Theta=(\CC, \Phi)$, $\Theta'=(\CC',\Phi')$ be  cone ensembles.
A diffeomorphism $F\colon K\rightarrow F(K)$  is called cone-hyperbolic 
from $\Theta'$ to $\Theta$ (on $K$) if we have\footnote{For our purposes, the second
 condition could be replaced by 
$(D _xF)^{\tr}(\real^d \setminus (\CC'_-\cup \CC'_+))\Subset \CC_0\cup \CC_-$.}
	\begin{align}\label{eq:cone:hyper}
		(\D _xF)^{\tr}\bigl (\real^d \setminus (\CC'_+\cup \CC'_0)\bigr )\Subset \CC_-
\, , \,\, \, 
(\D _xF)^{\tr}\bigl ( \real^d \setminus \CC'_+\bigr )
\Subset \CC_0\cup \CC_-\, ,\qquad
\forall x \in K \, .
	\end{align}
\end{definition}
The conditions \eqref{eq:cone:hyper} ensure that no parts of higher regularity in the anisotropic Banach spaces  of Section~\ref{s:BS} are mapped to parts of lower regularity (see \eqref{eq:l:multop:0} in the proof of \eqref{eq:l:multop:III} below).

We next introduce  a crucial ingredient to construct the anisotropic Banach spaces.

\begin{lemma}[Admissible Cone Ensembles for  $g_{-\alpha}$]\label{l:excone}
There exists an atlas $\AAA$, formed of a  finite open cover   $\{V_\omega\subseteq M\mid \omega\in\Omega\}$
of $M$ and  $C^r$ local diffeomorphisms $\kappa_\omega\colon V_\omega\rightarrow \real^d$,
such that
\begin{align}\label{eq:chart0}
\mbox{setting }K_\omega:=\kappa_\omega(V_\omega)\, , 
\mbox{ we have } \min_{\omega\ne \omega'} d(K_\omega,K_{\omega'})>1\, ,  \mbox{ and } K_M:=\cup_\omega \overline{ K_\omega} \mbox{ is compact,} 
\end{align}
and,   fixing coordinates $(x_1,\dots,x_d)\in\real^d$ and recalling  \eqref{eq:flow:gen}, the flow box condition 
\begin{align}\label{eq:chart:flowbox}
	(\D \kappa_{\omega}) (X|_{V_\omega})=\partial_{x_d}|_{\kappa_\omega(V_\omega)}\, ,
\end{align}
holds, and,
further, setting $V_{\alpha,\omega\omega'}\coloneqq V_{\omega}\cap g_{\alpha}(V_{\omega'})$ for each $\alpha\in\real$ and  $\omega$, $\omega'\in\Omega$
such that $V_{\omega}\cap g_{\alpha}(V_{\omega'})\ne \emptyset$, and also defining
$F_{-\alpha,\omega\omega'}\colon \kappa_{\omega}(V_{\alpha,\omega\omega'})\rightarrow \kappa_{\omega'}(V_{\omega'})$
as
$$
 F_{-\alpha,\omega\omega'}:=\kappa_{\omega'}\circ g_{-\alpha}\circ\kappa_{\omega}^{-1}\, ,
$$
there exists $\alpha_0>0$ and, for each $\omega$,  there exist    cone ensembles
 $\Theta_\omega$
	such that  the cone $D\kappa_\omega^{tr}(\CC_{\sigma, \omega})$ in the cotangent space contains the normal subspace $E^*_\sigma$ and is bounded away from
$E^*_\tau$ for $\tau\ne\sigma$, and, for all $\alpha\ge \alpha_0$, the map $F:=F_{-\alpha,\omega\omega'}$ is cone-hyperbolic from $\Theta_{\omega}$ to $\Theta_{\omega'}$
on $K:=\kappa_\omega(V_{\alpha,\omega\omega'})$.
\end{lemma}

\begin{remark}[{\cite[Remark 4.12]{Bbook}}]\label{forLeibniz}
For any $\Theta'<\Theta$  the identity 
is cone-hyperbolic from $\Theta$ to $\Theta'$.  If $F$ is  cone-hyperbolic 
from $\Theta'$ to $\Theta$, then there exists $\widetilde \Theta'< \Theta'$
such that  $F$ is  cone-hyperbolic 
from $\widetilde \Theta'$ to $\Theta$, and  there exists $\widetilde \Theta>\Theta$
such that  $F$ is  cone-hyperbolic 
from $\Theta'$ to $\widetilde \Theta$.
Thus, Lemma~\ref{l:excone} implies that there exist    cone ensembles
$\Theta'_\omega<\Theta_\omega$
such that   for all $\alpha\ge \alpha_0$ the map $F_{-\alpha,\omega\omega'}$ is cone-hyperbolic from $\Theta'_{\omega}$ to $\Theta_{\omega'}$
on $\kappa_\omega(V_{\alpha,\omega\omega'})$. 
Finally, the proof of  Lemma~\ref{l:excone}  provides an atlas $\AAA$, cone ensembles $\Theta'_\omega< \Theta_\omega$, and  $\alpha_0>0$,  such that  $\kappa_{\omega}\circ\kappa_{\omega'}^{-1}$ is cone-hyperbolic from $\Theta_{\omega'}$
to $\Theta'_{\omega}$ and
 $F_{-\alpha,\omega\omega'}$ is cone-hyperbolic from $\Theta'_{\omega}$ to $\Theta_{\omega'}$
for all $\alpha \ge \alpha_0$.
(Such   pairs $\{\Theta_\omega\}$, $\{\Theta'_\omega\}$ are called adapted to $\AAA$ and $g_\alpha$. They are used in Lemma~\ref{moll}.)
\end{remark}

\begin{proof} [Proof of  Lemma~\ref{l:excone}.]
Let $cc(A)$ denote the convex closure of a set $A$.\footnote{Taking the convex closure may be useful for tiny $\gamma>0$.}  By uniform continuity of the stable and unstable distributions,  
 setting
\begin{align*}
	\CC_{\sigma,\omega}^\gamma\coloneqq
cc \bigl ( \bigcup_{x\in V_\omega}(\D \kappa_{\omega}^{-1})^{\tr}\CC_\sigma^\gamma(x)\bigr ) \, ,\quad \sigma\in\{\pm\}\, , \, \, \omega\in\Omega\, ,\,\,
\gamma \in (0,1)\, , 
\end{align*}
we may choose  (small) $V_\omega$, $\kappa_\omega$ satisfying \eqref{eq:chart0}--\eqref{eq:chart:flowbox}
and $\gamma^*\in(0,1)$, $\tilde \gamma^*\in (0,\gamma^*)$ such that
\begin{align}\label{eq:chart:comp}
	&(\D _x\kappa_\omega)^{\tr}\CC^\gamma_{-,\omega}\Subset \CC_-^{\gamma^*}(x)\, ,
 \,\, \,\qquad\,\, 
(\D _x\kappa_\omega)^{\tr}\CC^\gamma_{+,\omega}\Subset \CC_+^{\gamma^*}(x)\, ,\,
&\forall \, \omega \in\Omega\, ,\,  \, x\in V_\omega\, ,
\\
	\label{eq:excone:I}
	&	(\D _{\kappa_\omega(x)}\kappa_\omega^{-1})^{\tr} 
		\CC_-^{\tilde\gamma^*}(x)\Subset \CC_-^{\gamma,\omega}\, ,\quad (\D _{\kappa_\omega(x)}\kappa_\omega^{-1})^{\tr}
\CC_+^{\tilde{\gamma}^*}(x)\Subset \CC_+^{\gamma,\omega}\, ,\,
&\forall \, \omega \in\Omega\, ,\,  \, x\in V_\omega\, .
	\end{align}
For $C_*\ge 1$, $\theta<1$ as in \eqref{eq:split:ineq}, and $\gamma$, $\tilde \gamma^*$, $\gamma^*$ as above,
 let $\alpha_0>0$ be such that  $C^2_*\theta^\alpha\gamma^*< \tilde \gamma^*$ for  all $\alpha\ge \alpha_0$.
By	  \eqref{eq:chart:comp} and Lemma~\ref{l:cone:inclusion}, we have, using the  transversal closed cones $\CC^\gamma_{\pm,\omega}$,
	\begin{equation*}
(\D g_{\alpha})^{\tr}(\D _x\kappa_\omega)^{\tr}\CC^\gamma_{+,\omega}
\Subset(\D g_{\alpha})^{\tr}\CC_+^{\gamma^*}(x)
\Subset  \CC_+^{\tilde\gamma^*}(g_{-\alpha}(x))\,,\,\,
\forall \alpha \ge \alpha_0\, , \, 
\forall x\in {V}_{\omega}\, .
\end{equation*}
We proceed similarly for $(\D g_{-\alpha})^{\tr}(\D _x\kappa_\omega)^{\tr}\CC^\gamma_{-,\omega}$,
using \eqref{eq:excone:I} and Lemma~\ref{l:cone:inclusion}.
Therefore, for $\alpha\ge \alpha_0$ and $\omega$, $\omega'\in\Omega$ such that $V_{\alpha,\omega\omega'}\ne\emptyset$, we have
	\begin{align}	
\label{eq:excone:II}
		{(\D F_{-\alpha,\omega'\omega})^{\tr}}\CC^\gamma_{-,\omega}
\Subset \CC^\gamma_{-,\omega'}\quad\mbox{ and } \quad
{(\D F_{\alpha,\omega'\omega}\underline{})^{\tr}}\CC^\gamma_{+,\omega}=
	{(\D F_{-\alpha,\omega'\omega}^{-1}\underline{})^{\tr}}\CC^\gamma_{+,\omega}
& \Subset \CC^\gamma_{+,\omega'}\, .
	\end{align}
Thus, there  exist  cone ensembles
$
\Theta_\omega$
such that for any $\alpha\ge \alpha_0$ and $\omega$, $\omega'\in\Omega$ with $V_{\alpha,\omega\omega'}\ne\emptyset$, using  the first claim of \eqref{eq:excone:II},  the first  inclusion in \eqref{eq:cone:hyper} holds for $F=F_{-\alpha,\omega\omega'}$,
$\Theta'=\Theta_{\omega'}$, $\Theta=\Theta_\omega$, while 
using  the second claim of \eqref{eq:excone:II},
the second  inclusion in \eqref{eq:cone:hyper} holds.
\end{proof}

\subsection{The Anisotropic Banach Spaces}\label{s:BS}

We  shall use a Paley--Littlewood decomposition:
Let $\chi_+\colon \real_+\rightarrow [0,1]$ be a $C^\infty$ map such that $\chi_+|_{{[0,1]}}\equiv 1$ and $\supp\chi_+\subseteq[0,2]$.  Set 
\begin{align}\label{eq:pldecomp}
	\Psi_0(\xi)\coloneqq\chi_+(|\xi|)\mbox{ and }\Psi_n(\xi)
\coloneqq\chi_+(|2^{-n}\xi|)-\chi_+(|2^{1-n}\xi|)\, , n\ge 1\, , \,
\xi \in \real^d \, .
\end{align}
This defines a $C^\infty$ partition of unity on $\real^d$ since 
$\sum_{n=0}^\infty\Psi_n(\xi)=
1$.
We have 
\begin{align*}
\Psi_n(\xi)=\Psi_1(2^{-n+1}\xi)\, , \mbox{ and thus }	\supp\Psi_n\subseteq\{\xi\in\real^d\mid 2^{n-1}\le|\xi|\le 2^{n+1}\}\, ,
\, \forall n \ge 1 \, .
\end{align*}

Given a cone ensemble $\Theta=(\CC,\Phi)$, we set 
\begin{align*}
	\Psi_{\sigma,0}(\xi):= \frac{\chi_+(\xi)} 3 \, ,\, \, \, \,	\Psi_{\sigma,n}(\xi):= \Psi_n(\xi)\Phi_{\sigma}\biggl (\frac{\xi}{|\xi|}\biggr ) \, ,\, \, \sigma\in\{\pm, 0\}\, ,\,n\ge 1\, .
\end{align*}
This also defines a $C^\infty$ partition of unity on $\real^d$ since 
$\sum_{n=0}^\infty\sum_{\sigma\in \{\pm, 0\}}\Psi_{n,\sigma}(\xi)=
1$.

Writing the inverse Fourier transform as
$\FFF^{-1}v(x)\coloneqq(2\pi)^{-d}\int_{\real^d}e^{\iunit \xi x}v(\xi)\dd\xi$,
where
$\xi x\coloneqq \langle\xi,x \rangle$
denotes the scalar product, we have
\begin{align*}
\|\FFF^{-1}\Psi_n\|_{L_1}=\|\FFF^{-1}\Psi_1\|_{L_1}< \infty,\quad	\|\FFF^{-1}\Psi_{\sigma,n}\|_{L_1}=\|\FFF^{-1}\Psi_{\sigma,1}\|_{L_1}<\infty\, , \, 
\sigma\in \{\pm ,0\}\, ,n\ge 1 \, .
\end{align*}
Analogous estimates hold for $\FFF^{-1}\Psi_{\sigma,0}$ and $\FFF^{-1}\Psi_0$.

Using the  convolution $v_1\ast v_2(x):=\int_{\real^d}v_1(x-y)\, v_2(y)\, \dd y$
of two distributions $v_1$, $v_2$, we associate
to any $\Psi$ with $\FFF^{-1}\Psi\in L_1(\real^d)$
a {\it pseudo-differential operator with symbol $\Psi$} via
\begin{equation}\label{pseu}
\Psi^{\Op}v(x):=((\FFF^{-1}\Psi) \ast v)(x)=\FFF^{-1}(\Psi \cdot \FFF v)(x)=
(2\pi)^{-d} \int_K\int_{\real^d}  e^{\iunit (x-y)\xi} \Psi(\xi) v(y) \dd \xi \dd y \, .
\end{equation}
Young's inequality, $\|v_1\ast v_2\|_{L_p}\le\|v_1\|_{L_1}\|v_2\|_{L_p}$
for all $p\in (1,\infty)$, gives
\begin{equation}\label{minimar}
\|\Psi^{\Op}v\|_{L_p}\le \| \FFF^{-1}\Psi\|_{L_1} \|v\|_{L_p} \,,\forall p \in (1,\infty) \, .
\end{equation}

 The (Bochner) space
$L_p(\real^d,\HH)$ associated to a  Hilbert space $\HH$ is defined by
$\|v\|_{L_p(\real^d,\HH)}\coloneqq
\|\|v\|_\HH\|_{L_p(\real^d)}$. 
The following  is a variant of the Marcinkiewicz theorem, generalising \eqref{minimar}:

\begin{theorem}[{See e.g. \cite[Thm 0.11.F]{Ta}}]
\label{l:multiplier}
Let $\HH_1$ and $\HH_2$ be Hilbert spaces, and let  $L(\HH_1,
\HH_2)$ be the space of bounded linear operators from $\HH_1$ 
to $\HH_2$ endowed with the operator norm. 
If $\QQ(\cdot)\in C^{\infty}(\real^d, L(\HH_1,\HH_2))$ satisfies
\begin{equation*}
\|\partial^\beta_\xi \QQ(\xi)\|_{L(\HH_1,
\HH_2)}\le
C_{\beta}(1+\|\xi\|^2)^{-|\beta|/2}\, , \,\mbox{ for each multi-index  } \beta\, ,
\end{equation*}
 then  the operator $\QQ^{\Op}$
defined for compactly supported continuous $a:\real ^d\to \HH_1$ by\footnote{With $\FFF^{-1}Q(x)= (2\pi)^{-d}\int_{\real^d}e^{\iunit x\xi} \QQ(\xi)\dd\xi$, for $x \in \real^d$, the notation \eqref{eq:BochnerOP} is compatible with \eqref{pseu}.}
	\begin{align}\label{eq:BochnerOP}
(\QQ^{\Op} a)(x) :=
 \frac{1}{(2\pi)^d}\int_{\real^d}\int_{\real^d}
e^{\iunit (x-y)\xi} \QQ(\xi)a(y)\, \dd y \, \dd\xi
\end{align}
extends for each  $1<p<\infty$ to a bounded operator from $L_p(\real^d,\HH_1)$ to $L_p(\real^d,\HH_2)$, 
and	\[\|{{Q}^\Op}\|_{L(L_p(\real^d,\BB_1), L_p(\real^d,\HH_2))}
\le\|\FFF^{-1}Q\|_{ L_1(\real^d,\LL(\HH_1,\HH_2))}\, .\]
\end{theorem}

\smallskip

We will mostly consider the  three cases
$$\HH_1=\HH_2=\complex\, ; \,\,  \HH_1=\HH_2=\ell_2(\complex)\, ; \,\,
\HH_2=\ell^c_2\,  \mbox{ and  } \HH_1=\ell^c_2
\mbox{ or } \HH_1=\ell^{c'}_2\, ; 
$$
 where $\ell_2^{c}$, $\ell_2^{c'}$ are the Hilbert spaces associated,
for fixed 
\begin{equation}\label{regul-sqt}
-(r-1)<s< 0< q\le t<r-1
\end{equation}
 and  $-(r-1)<s'<s$,
 $-(r-1)<q'\le q$,  $-(r-1)<t'< t$, to  
\begin{align*}
\|a\|_{\ell_2^c}\coloneqq
\bigl ( \sum_{\sigma,n}4^{c(\sigma)n}|a_{\sigma,n}|^2\bigr )^{\frac 12}
\, , 
\,\,\, \|a\|_{\ell_2^{c'}}\coloneqq
\bigl ( \sum_{\sigma,n}4^{c'(\sigma)n}|a_{\sigma,n}|^2\bigr )^{\frac 12}\, ,
\end{align*}
  where we set
\begin{align}\label{eq:reg:I}
	c(-)\coloneqq s\, ,\quad c(+)\coloneqq t\, ,\quad c(0)\coloneqq q	\, ,
\quad 
c'(-)\coloneqq s'\, ,\quad c'(+)\coloneqq t'\, ,\quad c'(0)\coloneqq q'\, .
\end{align}

Set
$C^{\tilde r}_0(K)
:=\{ f :\real^d \to \complex\mid
f \mbox{ is }  C^{\tilde r}\, ,\,\, \supp(f)\subset K\}$ for $K\subseteq \real^d$  compact with nonempty interior and $\tilde r\in [0,\infty]$.
We  introduce the basic building block for our anisotropic spaces: 

\begin{definition}[Local Anisotropic Norm and Banach Space]
Fix a cone ensemble $\Theta$ and\footnote{In the present work, we shall either have $s=q=t$ or
$
t-(r-1)< s<0< t <r-1
$ with $q\le t$.}
\begin{equation}\label{asabove}
p\in (1,\infty)\, , \quad -(r-1)<s\le q \le t<r-1\, . 
\end{equation}
For  a compactly supported $C^\infty$ function $v:\real^d\to \complex$, set
	\[
\|v\|_{W_{p,\Theta}^{s,t,q}}:=
\bigl\|\bigl (\sum_{n=0}^\infty 4^{ns}|\Psi{}^{\Op}_{-,n}v|^2
+4^{nt}|\Psi{}^{\Op}_{+,n}v|^2
+4^{nq}|\Psi{}^{\Op}_{0,n}v|^2\bigr )^{\frac 12} \bigr \|_{L_p(\real^d)}
\, .
\]
 For $K\subset\real^d$ compact with nonempty interior, the Banach space $W_{p,\Theta}^{s,t,q}(K)$ is defined to be the completion of $C^{\infty}_0(K)$ under $\|\cdot\|_{W_{p,\Theta}^{s,t,q}}$.

We shall also use the auxiliary semi-norm
$
\|v\|_{W_{p,\Theta}^{q}}:= 
\bigl\|\bigl (\sum_{n=0}^\infty 4^{nq}
|\Psi{}^{\Op}_{0,n}v|^2\bigr )^{\frac 12} \bigr \|_{L_p(\real^d)}$.
\end{definition}

The definition of $\|v\|_{W_{p,\Theta}^{s,t,q}}$ is just like\footnote{This is a ``Sobolev'' (Triebel--Lizorkin) type norm. A ``Besov-H\"older'' version \cite[(2.3)]{Baladi2007} should work too, in particular for $p=\infty$.}
 \cite[(2.4)]{Baladi2007}, except that we have three cones instead of two.  
We record the following result for convenience;
the continuous inclusion claim is obvious, while
	 the compact inclusion claim --- which is in fact not used in the present work --- is proved exactly like \cite[Prop. 5.1]{Baladi2007}, using Arzel\`a--Ascoli: 

\begin{lemma}[Continuous and compact embeddings, local spaces]\label{l:contembed}
 Let $K\subset \real^d$ be compact with nonempty interior.	Let  $\Theta$ be a cone ensemble,
	fix $p$ and $s,q,t$ as in \eqref{asabove}. For any  $s'\le s$, $q'\le q$, and $t'\le t$, the inclusion
	$W_{p,\Theta}^{s,t,q}(K)\subseteq W^{s',t',q'}_{p,\Theta}(K)$
	is continuous. If  $s'<s$, $t'<t$, and
 $q'<q$, the inclusion
	$W_{p,\Theta}^{s,t,q}(K)\subseteq W_{p,\Theta}^{s',t',q'}(K)$
	is compact.
\end{lemma}

Next, letting  $\|v\|_{W_{p}^{t}}=\|(1+\Delta)^{t/2}v\|_{L_p}$ denote  the classical isotropic Sobolev (Triebel--Lizorkin) norm, the
arguments\footnote{See also \cite[(4.17)--(4.18), footnote 15]{Bbook}.} of
 \cite[App A]{Baladi2007}  give,
for $p$, $s$, $t$, $q$ as in \eqref{asabove}, a constant $C\in (1,\infty)$ with\footnote{Generalise \cite[(A.2), (A.4),(A.5)]{Baladi2007} to three cones, noting that, if $s=t(=q)$, then $\Theta$ plays no role for the norm noted $W^{\Theta, s,t,(q,)p}_{\dagger\dagger}$ there, so that the
condition $\Theta'>\Theta$ for \cite[(A.5)]{Baladi2007} is immaterial.}
\begin{equation}\label{compareiso}C^{-2} \|v\|_{W_{p}^{-r+1}}
\le 
C^{-1} \|v\|_{W_{p}^{s}}
\le 
\|v\|_{W_{p,\Theta}^{s,t,q}}\le C \|v\|_{W_{p}^{t}}
\le C^2 \|v\|_{W_{p}^{r-1}}\, , \, \, \forall
v \in C^\infty_0(K) \, .
\end{equation}
It follows that $C^{r-1}_0(K)\subset W_{p,\Theta}^{s,t,q}(K)$  (as a dense subset).

  We are finally ready to define our anisotropic space of distributions on $M$:

\begin{definition}[Anisotropic Banach space]\label{d:bs:ani}
Let  $\AAA=\{\VV=\{V_\omega\}_{\omega \in \Omega},\{\kappa_\omega:V_\omega\to \real^d\}_{\omega\in \Omega}\}$, $\alpha_0>0$, and
cone ensembles $\{\Theta_\omega\}_{\omega \in \Omega}$ 
admissible for $\{g_{-\alpha}\}_{\alpha\ge \alpha_0}$ 
be as given by
 Lemma~\ref{l:excone}.
 Fix  a $C^{r}$ partition of unity $\{\vartheta_\omega\colon M \rightarrow [0,1]\}_{\omega \in\Omega}$, subordinate to $\VV$, that is, with
$\supp \vartheta_\omega \subset V_\omega$.
 For $p$, $s$, $t$, $q$ as in \eqref{asabove}, we put\footnote{Integration with respect to $\alpha$ allows one to handle the  times $\alpha \in [0,\alpha_0]$ where the flow is not sufficiently hyperbolic. (This is similar to \cite[Definition 8.1]{Faure_2017} and, replacing the integral by a supremum, to
\cite[(3.6)]{Baladi2012}. See \cite[Lemma 6]{Butterley_2013}
for a slightly different trick.) The identity \eqref{Htrick} shows why the $L_2$ norm is natural.} for
$\varphi\in{C^{\infty}(M)}$ (extending $\vartheta_\omega\circ\kappa_\omega^{-1}$   from $\kappa_\omega(V_\omega)$ to $\real^d$ by zero),
\begin{align}\label{eq:norm}
\|\varphi\|_{{W_p^{s,t,q}}} := \bigl (\sum_{\omega\in\Omega}\int_0^{\alpha_{0}}
\|(\vartheta_\omega\cdot \LL_{\alpha,V} \varphi)\circ \kappa_\omega^{-1}\|_{W_{p,\Theta_\omega
}^{s,t,q}}^2\dd\alpha\bigr )^{\frac 12}\, .
\end{align}
The Banach space $W_p^{s,t,q}(M)$ is defined to be the completion of $C^{\infty}(M)$ under $\|\cdot\|_{W_p^{s,t,q}}$.
\end{definition}

We show in Lemma~\ref{interpolation} that the scale $W_p^{s,t,q}(M)$ is an interpolation
scale (the proof also shows this for the  scale $W_{p,\Theta}^{s,t,q}(K)$).
In Lemma~\ref{moll} we use mollifiers to approximate the identity.

Note that $W_p^{s,t,q}(M)$ depends on  the atlas $\AAA$, the cone ensembles $\Theta_\omega$, and $\alpha_0$.  It follows from \eqref{compareiso} and the bound
$\sup_{\alpha\in [0, \alpha_0]}\|\LL_{\alpha, V}\varphi\|_{C^{r-1}}\le C(\alpha_0)
\|\varphi\|_{C^{r-1}}$ that $C^{r-1}(M)\subset
W_p^{s,t,q}(M)$, as a dense subset. It is not hard to show that ${W_p^{s,t,q}}(M)$
is contained in the set of distributions of order $r-1+d/p$ on $M$.
For $p=2$, the space $W^{s,t,q}_p$ is a Hilbert space since it satisfies the parallelogram law
\cite[Lemma 15.2]{Blanchard_2015} 
\begin{equation}\label{Htrick}
\|\varphi_1+\varphi_2\|_{{W_2^{s,t,q}}}^2+
\|\varphi_1-\varphi_2\|_{{W_2^{s,t,q}}}^2
=2\|\varphi_1\|_{{W_2^{s,t,q}}}^2+2\|\varphi_2\|_{{W_2^{s,t,q}}}^2
\, .\end{equation}

Clearly, if $s'\le s$, $q'\le q$, and $t'\le t$, we have the continuous  injection
$W_{p}^{s,t,q}(M)\subseteq W^{s',t',q'}_{p}(M)$.
Due to the integration over $\alpha$ in the definition
\eqref{eq:norm} of the norm, the compact inclusion from Lemma~\ref{l:contembed} does not carry over\footnote{See e.g. \cite{aubin} and \cite{amann} for relevant results in this context.}  automatically to the  space ${W_p^{s,t,q}}(M)$ (despite the fact that $M$ is compact).  We provide a direct proof of the following lemma instead:

\begin{lemma}[Compact embeddings]\label{l:compact}
Fix $p\in (1,\infty)$, and $s$, $q$ , $t$ as in \eqref{asabove}.  If $s'<s$, $q'<q$, and $t'<t$
satisfy $-r-1<s'\le q'\le t'$ then  the inclusion
$W_p^{s,t,q}(M)\subset W^{s',t',q'}_p(M)$
is compact.
\end{lemma}

\begin{proof}
It suffices to show that the inclusion
$W_p^{s,t,q}(M)\subset W^{s',s',s'}_p(M)$ is compact.
Indeed,  for every $\epsilon >0$ there exists $C(\epsilon) <\infty$
such that, for any $v \in W^{s,t,q}_{p,\Theta}$, 
\begin{equation}\label{trick}
\|v \|_{ W^{s',t',q'}_{p,\Theta}}\le \epsilon
\|v \|_{ W^{s,t,q}_{p,\Theta}}+C(\epsilon)\|v \|_{ W^{s',s',s'}_{p,\Theta}} 
\end{equation}
(This is easy to prove along the lines of \cite[Remarks 2.22, 4.27]{Bbook}).
Thus, 
$$
\|\varphi \|_{ W^{s',t',q'}_{p}}\le \epsilon
\|\varphi \|_{ W^{s,t,q}_{p}}+ C(\epsilon) \|\varphi \|_{ W^{s',s',s'}_{p}}\,,\quad
\forall\varphi \in W^{s,t,q}_p(M)\, .
$$

We now show that $W_p^{s,t,q}(M)\subset W^{s',s',s'}_p(M)$ is compact. For a $C^r$ diffeomorphism $F:K\to F(K)$  and $f\in C^{r-1}_0(K)$,
we introduce the  local transfer operator
\begin{align}\label{eq:to:loc}
	\MM_{F,f}\colon C^{r-1}(F(K))\rightarrow C^{r-1}_0(K)\, ,\,\,
\MM_{F,f}(v)= f\cdot (v\circ F) \, .
\end{align}
The key  fact is that the operator $\MM_{ F,  f}$ is  bounded  for the classical Sobolev norm $W^{s}_p$ on $\real ^d$ if $s\in (-r-1,r-1)$ (apply e.g. the results of \cite[Chapter 2]{Bbook}), with norm depending only on  $\|f\|_{C^{r-1}}$ and $\|F\|_{C^r}$.
In particular,  since $F_{\alpha,\omega\omega'}:=\kappa_{\omega'}\circ g_{\alpha}\circ \kappa_\omega^{-1}$ and $f_{\alpha,\omega}:=(\vartheta_\omega \phi_{-\alpha})\circ \kappa_\omega^{-1}$
 are $C^r$, respectively $C^{r-1}$  (on 
 $K_\omega$)
uniformly in $\alpha \in[0,\alpha_0]$,  with $\phi_{-\alpha}\phi_\alpha\equiv 1$,\footnote{It is useful here that $|\phi_\alpha|$ is bounded away from zero.}  decomposing 
$$
(\vartheta_\omega \varphi)\circ \kappa_\omega^{-1}=
\bigl (\vartheta_\omega  \phi_{-\alpha}\sum_{\omega'} (\vartheta_{\omega'}\phi_\alpha [\varphi \circ g_{-\alpha}]) \circ \kappa_{\omega'}^{-1}
\circ \kappa_{\omega'}\circ g_\alpha\bigr ) \circ \kappa_\omega^{-1}\, ,
$$
  we have, setting $\MM_{\alpha_0}=\sup_{\omega, \omega',\alpha \in [0,\alpha_0]} \|\MM_{F_{\alpha,\omega\omega'}, f_{\alpha,\omega}}\|_{W^s_p}<\infty$, 
\begin{equation}\label{keyy}
  \| (\vartheta_{\omega}\varphi)\circ \kappa_{\omega}^{-1}\|_{W_{p
}^{s}}\le \MM_{\alpha_0}
\sum_{\omega'} \| (\vartheta_{\omega'} [\LL_\alpha \varphi])\circ \kappa_{\omega'}^{-1}\|_{W_{p
}^{s}}\, ,\, \forall \varphi \, ,\forall \omega\, ,\,\forall \alpha \in [0,\alpha_0] \, .
\end{equation}

Let now $\varphi_m$ be a sequence in the unit ball of $W_p^{s,t,q}(M)$. 
By definition \eqref{eq:norm} of the norm,   for every $m$ there exists $\alpha(m)\in [0,\alpha_0]$ such that
$$
  \|(\vartheta_{\omega'}\cdot \LL_{\alpha(m),\phi_{\alpha(m)}} \varphi_{m})\circ \kappa_{\omega'}^{-1}\|_{W_{p,\Theta_{\omega'}
}^{s,s,s}}^2\le \|(\vartheta_{\omega'}\cdot \LL_{\alpha(m),\phi_{\alpha(m)}} \varphi_{m})\circ \kappa_{\omega'}^{-1}\|_{W_{p,\Theta_{\omega'}
}^{s,t,q}}^2
\le \frac1 {\alpha_0 }\, , \,\, \forall \omega'\, .
$$
Thus, using \eqref{keyy}, and recalling  $C$ from
 \eqref{compareiso},
\begin{equation}\label{upp}
\| (\vartheta_{\omega}\varphi_m)\circ \kappa_{\omega}^{-1}\|_{W_{p
}^{s}}\le \MM_{\alpha_0} \cdot \sum_{\omega'}  \|(\vartheta_{\omega'}\cdot \LL_{\alpha(m),\phi_{\alpha(m)}} \varphi_{m})\circ \kappa_{\omega'}^{-1}\|_{W_{p
}^{s}}^2\le  M_{\alpha_0}\frac {C \cdot \# \Omega} {\alpha_0 }\, , \,\, \forall \omega'\, ,\, \forall m\, .
\end{equation}

Assume  for a contradiction that there is $\epsilon >0$ and,
for any $k_0\ge 1$, there are $k,\ell \ge k_0$ with 
	\begin{align}\label{eq:l:compact}
\|\varphi_{k}-\varphi_{\ell}\|_{W_p^{s',s',s'}}>\epsilon\, .
	\end{align} 
By definition,   this implies that there exists $\omega \in \Omega$ with
$$
\int_0^{\alpha_{0}}
\|(\vartheta_\omega\cdot \LL_{\alpha,V} (\varphi_{k}-\varphi_{\ell}))\circ \kappa_\omega^{-1}\|_{W_{p,\Theta_\omega
}^{s',s',s'}}^2\dd\alpha
> \frac{\epsilon^2}{ \# \Omega}\, .
$$
Now, using again the key fact, and setting $M'_{\alpha_0}=\sup_{\omega, \omega',\alpha \in [0,\alpha_0]} \|\MM_{F_{-\alpha,\omega\omega'}, f_{-\alpha,\omega}}\|_{W^{s'}_p}$, we get
\begin{align}
\nonumber \int_0^{\alpha_{0}}&
\|(\vartheta_\omega\cdot \LL_{\alpha,V} (\varphi_{k}-\varphi_{\ell}))\circ \kappa_\omega^{-1}\|_{W_{p,\Theta_\omega
}^{s',s',s'}}^2\dd\alpha\le C \alpha_0 \sup_{\alpha\in [0,\alpha_0]}
\| (\vartheta_{\omega}\cdot\LL_{\alpha,\phi_{\alpha}}(\varphi_{k}-\varphi_{\ell}))\circ \kappa_{\omega}^{-1}\|_{W_{p}^{s'}}^2\\
 &\le C \alpha_0 \cdot (M'_{\alpha_0})^2
\sum_{\omega '}\| (\vartheta_{\omega'}(\varphi_{k}-\varphi_{\ell}))\circ \kappa_{\omega'}^{-1}\|_{W_{p}^{s'}}^2
\label{contrad}
\,  .
\end{align}
Since \eqref{upp} implies that $(\vartheta_{\omega'}\cdot (\varphi_k-\varphi_\ell))\circ \kappa_{\omega'}^{-1}$ is a sequence in a bounded subset of $W_{p}^{s}(K_M)$ for each $\omega'$, and 
since the embedding $W_{p}^{s}(K_M)\subset W_{p}^{s'}(K_M)$ is compact,
we find $k_0$ such that \eqref{contrad} is 
smaller than $\epsilon^2/ \#\Omega$ for all
$k\ge \ell \ge k_0$ and thus the desired contradiction with \eqref{eq:l:compact}.
\end{proof}

\section{Properties of the Transfer Operator, the Generator, the Resolvent}\label{s:prop}

\subsection{Basic Estimates on the Local Anisotropic Space}
The  natural ordering $-<0<+$ on $\{-, +, 0\}$ is compatible with our choice $s=c(-)\le q=c(0)\le t=c(+)$
from \eqref{eq:reg:I}. Inspired by \cite{Baladi2007},  we introduce the following definition:

\begin{definition}[Arrow relation]\label{d:arrow}
For $K\subset\real^d$  compact with nonempty interior, let $F:K\to F(K)$ be a $C^r$
cone hyperbolic diffeomorphism from $\Theta'$ to $\Theta$  on $K$. For a covering $\widetilde \Phi'$
of $\Theta'$, set
\begin{align}\label{covv}
	|F|_{\tau}&:=
\sup_{\substack{x\in K\\\eta\in \supp\widetilde{\Phi}'_{\tau}}} 
|\D F(x)^{\tr}\eta|\,, 
\,\,\,\,\,
	|F^{-1}|_{\sigma}:=
\sup_{\substack{x\in F(K)\\\xi\in\supp\Phi_{\sigma}}}
|\D F^{-1}(x)^{\tr}\xi | \, .
\end{align}
Fix $s<0<q<t$. For $n$, $\ell\ge 0$, and $\sigma$, $\tau\in\{\pm, 0\}$, we say that	$(\tau,\ell)\hookrightarrow_K (\sigma,n)$ if

\begin{align}\label{obvv}
(2^{sn}\le |F|^t_+\mbox{ or } 2^{-q\ell} \le |F^{-1}|^{|s|}_-)\mbox{ and }
	\sigma\le \tau \mbox{ and }|F^{-1}|_{\sigma}^{-1}2^{-4}\le 2^{n-\ell}\le 2^4|F|_{\tau}
\, ,
\end{align}
and we say that $(\tau,\ell)\not\hookrightarrow_K (\sigma,n)$ otherwise.
\end{definition}

Recalling \eqref{eq:pldecomp}, let $\widetilde{\Psi}_0$, $\widetilde{\Psi}_1\in C^\infty$ be such that $\widetilde{\Psi}_{0}|_{\supp\Psi_0}\equiv 1$ and $\widetilde{\Psi}_{1}|_{\supp\Psi_1}\equiv 1$. Set  $\widetilde{\Psi}_n(\xi)\coloneqq \widetilde{\Psi}_1(2^{-n+1}\xi)$ for $n\ge 2$.
 With\footnote{\label{sixteen}In our application below, $|F|_+<1$ while $|F^{-1}|_->1$.}  \eqref{eq:split:estimate:I}, \eqref{eq:split:estimate:I'}, and \eqref{hook00}, the following lemma shows the usefulness of the arrow relation:

\begin{lemma}\label{forcomp}If $F$ is a $C^r$ cone hyperbolic diffeomorphism from $\Theta'$ to  $\Theta$ on $K$, there exist a covering $\widetilde{\Phi'}$ of $\Theta'$ and
a constant  $C_1=C_1(F,K)>0$ such
 that, setting,
\begin{align}\label{eq:suppext}
	\widetilde{\Psi'}_{\sigma,0}(\xi) := \frac{\chi_+(\xi)}3 \, ,\,\, 	\widetilde{\Psi'}_{\sigma,n}(\xi) := \widetilde{\Psi}_n(\xi) \widetilde{\Phi'}_\sigma\biggl (\frac{\xi}{|\xi|}\biggr )\, ,
\, \,\,\xi \in \real^d\, , \sigma\in\{\pm, 0\}\, , \,\, n\ge 1\, ,
\end{align} 
we have
	\begin{align}\label{eq:l:multop:III}
		\inf_{x\in K} d
\bigl (\supp\Psi_{\sigma, n}-\D F(x)^{\tr}\supp{\widetilde{\Psi}}' _{\tau, \ell}\bigr )\ge  C_12^{\max\{n,\ell\}}\,,\qquad
\forall (\tau,\ell)\not\hookrightarrow_K(\sigma,n)
 \, .
	\end{align}
\end{lemma}

\begin{proof}If $\sigma\le \tau$, then \eqref{eq:l:multop:III} follows from \eqref{obvv} (without using cone-hyperbolicity): Indeed, if $n\ge \ell$, 
	\begin{align*}
		d(\supp\Psi_{\sigma, n}, \D F(x)^{\tr}\supp{\widetilde{\Psi}}' _{\tau, \ell})&\ge 2^{n-1}-2^{\ell+2}|F|_{\tau}
		=2^{n-1}(1-2^{\ell-n+3}|F|_{\tau})
>2^{n-2}\, , \,\forall x \in K\, ,
\end{align*}
while if $n<\ell$, we have 
\begin{align*}
	d(\supp\Psi_{\sigma, n}
-\D F(x)^{\tr}\supp{\widetilde{\Psi}}' _{\tau, \ell}\bigr  )&\ge	2^{\ell-1}(2^{n-\ell}-2^3|F|_{\tau})>2^{\ell+2}
|F|_{\tau}\, , \,\forall x \in K\, ,
	\end{align*}

If $\sigma> \tau$,
then either $\tau=0$ and $\sigma=+$, or  $\tau=-$ and $\sigma\in\{0,+\}$. In both cases, cone-hyperbolicity of $F$ implies
\begin{align}\label{eq:l:multop:0}
		\bigcup_{x\in K} (\supp\Psi_{\sigma})\cap  (\D F(x)^{\tr}\supp {\Psi}'_{\tau}) =\{0\} \, ,
\end{align}
which is a trivial intersection of closed cones.	Hence  
there exists a covering $\widetilde \Phi'$ such that 	$\bigcup_{x\in K} (\supp\Psi_{\sigma})\cap  (\D F(x)^{\tr}\supp {\widetilde \Psi}'_{\tau}) =\{0\}$, and  \eqref{eq:l:multop:III} holds for suitable $C_1$. 
\end{proof}

For a $C^r$  diffeomorphism $F:K\to F(K)$,  cone-hyperbolic from
$\Theta'$ to $\Theta$ and a covering $\widetilde \Phi'$ of $\Theta'$ satisfying  \eqref{eq:l:multop:III} and $f\in C^{r-1}_0(K)$,
recalling the weighted composition operator
$\MM_{F,f}(v)=f\cdot (v \circ F)$ from 
\eqref{eq:to:loc}, set, 
for $a=(a_{\tau,\ell})\in L_p(\real^d,\ell_2^{c})$
(Lemmas~\ref{l:bound:I} and~\ref{l:bound:II} providing the necessary summability), 
\begin{align*}
	(Q_{\hookrightarrow_K}^{\Op}a)_{\sigma,n}&:= \Psi^{\Op}_{\sigma,n}\sum_{(\tau,\ell)\hookrightarrow_K  (\sigma,n)}\MM_{F,f}(a_{\tau,\ell})\, , \,\,\,\,
		(Q_{\not\hookrightarrow_K}^{\Op}a)_{\sigma,n}
:=  \Psi^{\Op}_{\sigma,n}
\sum_{(\tau,\ell)\not\hookrightarrow_K 
(\sigma,n)}\MM_{F,f}( \widetilde {\Psi}'^{\,\, \Op}_{\tau,\ell} a_{\tau,\ell} )\, .
\end{align*}Then, taking
$a_{\tau,\ell}:= {\Psi'}^{\,\, \Op}_{\tau,\ell}v$ for the ensemble $\Theta$,
we have 
\begin{align}\label{eq:multop}
\|\MM_{F,f}v\|_{W_{p,\Theta}^{s,t,q}}=
\|Q_{\not\hookrightarrow_K}^{\Op}a
+Q_{\hookrightarrow_K}^{\Op}a\|_{L_{p}(\real^d,\ell_2^c)} \, .
\end{align}

Lemma~\ref{l:bound:I} describes the $\hookrightarrow$ term in the decomposition above. 
It   will give
the ``contracting'' factor $C_\pm$ in Lemma~\ref{l:to:loc:bound}
for  $\sigma=\pm$,  while the term with $C_0$ in Lemma~\ref{l:to:loc:bound}  for $\sigma=0$
will become compact  for the resolvent,  see Lemmas~\ref{l:norm:loc:derivative} and~\ref{l:to:bound}).

\begin{lemma}[The Bounded Term]\label{l:bound:I}
Fix  $p\in(1,\infty)$ and $s$, $q$, $t$ as in \eqref{regul-sqt}. 
There exists  $C<\infty$, such that for each compact $K\subset \real^d$  with nonempty interior, each
$C^r$ cone-hyperbolic\footnote{Cone hyperbolicity is not really needed.}diffeo\-morphism $F:K\to F(K)$
from $\Theta'$ to $\Theta$, and each covering
$\widetilde \Phi'$ of $\Theta'$, 
\begin{align*}
\|Q_{\hookrightarrow_\KK}^{\Op}a\|_{L_{p}(\real^d,\ell_2^c)}
&\le C
 \max\{|F|_{+}^t,|F^{-1}|_{-}^{|s|}\}
\sup_K \bigl |f  |\det DF|^{-1/p} \bigr |\cdot
\|a\|_{L_{p}(\real^d,\ell_2^c)}\\
& \quad\quad
+ C
|F|_{0}^q
 \sup_K \bigl |f  |\det DF|^{-1/p}\bigr |
\cdot \,
\bigl \| \bigl (\sum_{\ell}4^{q\ell}
\bigl | a_{0,\ell}\bigr |^2\bigr )^{\frac 12} \bigr \|_{L_p}\, ,\,\,\forall f\in C^{r-1}_0(K)
\, .
\end{align*}
\end{lemma}

\begin{proof} 
Recall \eqref{eq:reg:I}. 
There exists $C<\infty$, independent of $F$,  such that 
\begin{align}
	\sum_{n:(\tau,\ell)\hookrightarrow_K (\sigma,n)}
2^{ c(\sigma)n-c(\tau)\ell}&=
\sum_{n:(\tau,\ell)\hookrightarrow_K (\sigma,n)}
2^{ (c(\sigma)-c(\tau))n+c(\tau)(n-\ell)}\nonumber
	\label{eq:split:estimate:I}\le
 \sum_{n:(\tau,\ell)\hookrightarrow_K  (\sigma,n)}2^{ c(\tau)(n-\ell)}\\
&\le C \max\{|F|_{+}^t,|F^{-1}|_{-}^{|s|}\}\, ,\
\forall  (\tau,\ell)\, , \, \sigma\, ,  \mbox{ with } (\sigma,\tau)\ne (0,0)\, .
\end{align}
Similarly, up to taking a larger constant $C<\infty$,  we have 
\begin{align}	\label{eq:split:estimate:I'}
	\sum_{\ell: (\tau,\ell)\hookrightarrow_K (\sigma,n)}&
2^{ c(\sigma)n-c(\tau)\ell}\le
 C \max\{|F|_{+}^t,|F^{-1}|_{-}^{|s|}\}\, ,\
\forall  (\sigma,n)\, , \, \tau\, ,  \mbox{ with } (\sigma,\tau)\ne (0,0)\, .
\end{align}
Theorem~\ref{l:multiplier} applied to $\HH_1=\HH_2=\ell_2^c$ and  $(\QQ b)_{\sigma, n}(\xi)=\Psi_{\sigma,n} (\xi) b_{\sigma, n}(\xi)$ gives $D_1$ 
such that
	\begin{align*}
		&\|Q^\Op_{\hookrightarrow _K} a\|_{L_{p}(\real^d,\ell_2^c)}
\le D_1
\bigl\| \bigl (\sum_{\sigma,n}4^{c(\sigma)n}\bigl | \sum_{(\tau,\ell)\hookrightarrow_K (\sigma,n)}\MM_{F,f} a_{\tau,\ell} \bigr |^2\bigr )^{\frac 12}\bigr \|_{L_p}\, ,
\quad \forall f, F, a\, .
\end{align*}
Set $\lambda_{F,s,t}=\max\{|F|_{+}^t,|F^{-1}|_{-}^{|s|}\}$. By  Cauchy--Schwarz, and \eqref{eq:split:estimate:I}--\eqref{eq:split:estimate:I'}, we find   $D_2$, $D_3$ such that 
	\begin{align*}
		& 3 D_1 \sum_{(\tau,\sigma)\ne (0,0)}\|\bigl (\sum_{n}
\sum_{j:(\tau,j)\hookrightarrow_K (\sigma,n)}2^{c(\sigma)n-c(\tau)j} 
\sum_{\ell:(\tau,\ell)\hookrightarrow_K (\sigma,n)}2^{c(\sigma)n+c(\tau)\ell}|\MM_{F,f} a_{\tau,\ell}|^2\bigr )^{\frac 12}\|_{L_p}\\
&\le D_2\sum_{(\tau,\sigma)\ne (0,0)}
\bigl \|\bigl (
\lambda_{F,s,t}
\sum_{n}\sum_{\ell:(\tau,\ell)\hookrightarrow_K (\sigma,n)}2^{c(\sigma)n+c(\tau)\ell}
|\MM_{F,f} a_{\tau,\ell}|^2\bigr )^{\frac 12}\bigr \|_{L_p}\\
&= D_2\sum_{(\tau,\sigma)\ne (0,0)}
\bigl \|\bigl ( \lambda_{F,s,t}
\sum_{\ell}4^{c(\tau)\ell}|\MM_{F,f} a_{\tau,\ell}|^2
\sum_{n:(\tau,\ell)\hookrightarrow_K (\sigma,n)}2^{c(\sigma)n-c(\tau)\ell}\bigr )^{\frac 12}
\bigr \|_{L_p}\\
&\le D_3 \, \sup_K\bigl |f |\det DF|^{-1/p}\bigr |\sum_{\sigma,\tau}
\bigl \|\bigl (\lambda_{F,s,t}^2\, 
 \sum_{\ell}4^{c(\tau)\ell}
| a_{\tau,\ell}|^2\bigr )^{\frac 12} \|_{L_p}\, ,
\, \, \forall f, \, F, \, a \, .
	\end{align*}	
Finally, since  there exists $C_{00}<\infty$,  independent of $F$, such that
\begin{equation}
\label{hook00}
\sum_{n:(0,\ell)\hookrightarrow_K (0,n)}
2^{ c(0)(n-\ell)}\le C_{00}|F|_{0}^{q}\, , \,\,\forall \ell\, ,
\qquad\mbox{ and }
	\sum_{\ell: (0,\ell)\hookrightarrow_K (0,n)}
2^{ c(0)(n-\ell)}\le C_{00}|F|_{0}^{q}\, , \,\,\forall n\, ,
\end{equation}
 we find
 $D_4$ such that for all $F$, $f$, and $a$
$$
\bigl\| \bigl (\sum_{n}4^{c(0)n}\bigl | \sum_{(0,\ell)\hookrightarrow_K (0,n)}\MM_{F,f} a_{\tau,\ell} \bigr |^2\bigr )^{\frac 12}\bigr \|_{L_p}
\le D_4 \, \sup_K \bigl |f|\det DF|^{-1/p}|F|_{0}^q
\bigr |
\bigl\| (\sum_{\ell}4^{q\ell}| a_{0,\ell}|^2)^{\frac 12}\|_{L_p} \, .
$$
\end{proof}

We next bound the other term in the decomposition 
\eqref{eq:multop} of $\MM_{F,f}$. For this, we need the following strengthening of condition  \eqref{regul-sqt}:
\begin{equation}\label{regulS-sqt}
t-(r-1)<s<0<q< t \, .
\end{equation}

\begin{lemma}[The Compact Term]\label{l:bound:II}
	Fix $p\in (1,\infty)$, and fix $s$, $q$, $t$   as in \eqref{regulS-sqt}. Let
$s'<s$, $q'\le q$ and
 $t'<t$ satisfy $t-(r-1)<s'<0<q'\le t'$. Let $F$ be a $C^r$ cone-hyperbolic diffeomorphism
from $\Theta'$ to $\Theta$ on $K$,
let $\widetilde \Phi'$ be given by Lemma~\ref{forcomp}, and let $f\in C^{r-1}_0(K)$. Then there exists  $C(F,f)<\infty$ such that
	\begin{equation*}
\|Q_{\not\hookrightarrow_\KK}^{\Op}a\|_{L_{p}(\real^d,\ell_2^c)}
\le C(F,f)\cdot
\bigl \| \|a\|_{\ell_2^{c'}}  \bigr \|_{L_p}\, , \forall a\, .
\end{equation*}
\end{lemma}

\begin{proof}
We shall use Lemma~\ref{forcomp} and integration by parts, along the lines of \cite[pp. 144--147]{Baladi2007}.

Write, for  $x\in\real^d$ and $(\tau,\ell)\not\hookrightarrow_K(\sigma,n)$,
	\begin{align*}
		( J_{\sigma,n}^{\tau,\ell}a)(x)
&:=
\frac{(2\pi)^{2d}}{2^{(n+\ell)d}}{\Psi}^\Op_{\sigma,n}
\MM_{F,f}(\widetilde{\Psi}'^\Op_{\tau,\ell}b_{\tau,\ell})(x)\\
		&=\int_{y\in K}
\int_{\real^d\times\real^d\times\real^d }
e^{\iunit 2^n\tilde{\xi} (x-y)}e^{\iunit 2^\ell\tilde{\eta}(F(y)-w)} \Psi_{\sigma,1}(\tilde{\xi})
{\widetilde{\Psi}}' _{\tau,1}(\tilde{\eta})f(y) a_{\tau,\ell}(w) \dd w\, \dd\tilde{\xi}\, \dd\tilde{\eta}\,  \dd y\, ,
	\end{align*}
where we used the change of variables
$\tilde{\xi}= 2^{-n}\xi$ and  $\tilde{\eta}= 2^{-\ell}\eta$. 	Integrating  by parts $r-1$ times (see Lemmas~\ref{l:partint}--\ref{l:partint:reg}) in $y$, and using   \eqref{eq:l:multop:III}, we rewite
$( J_{\sigma,n}^{\tau,\ell}a)(x)$ as
\begin{align*}
		&\int_{y\in K}
\int_{\real^d\times\real^d\times\real^d }
e^{\iunit 2^n\tilde{\xi} (x-y)}
e^{\iunit 2^\ell\tilde{\eta}(F(y)-w)} \Psi_{\sigma,1}(\tilde{\xi})
{\widetilde{\Psi}}'_{\tau,1}(\tilde{\eta})
\frac{f_{r-1,n,\ell}(\tilde\eta,\tilde{\xi},y)}
{2^{\max\{n,\ell\}(r-1)}} 
a_{\tau,\ell}(w) \dd w\, \dd\tilde{\xi}\, \dd\tilde{\eta}\, \dd y\, ,
	\end{align*}
where  all partial derivatives of $f_{r-1,n,\ell}(\tilde\eta,\tilde{\xi},y)$ with respect to $\tilde{\eta}$ and $\tilde{\xi}$ are bounded  by a constant $C_2(F,f)$ uniformly in
$n$, $\ell$, and $(\tilde{\xi},\tilde{\eta},y)\in \supp\Psi_{\sigma, 1}\times\supp{\widetilde{\Psi}}' _{\tau, 1}\times K$. Define $b:\real^d\to [0,1]$ by
	\[
b(y)\coloneqq
1 \mbox{ if }  |y|\le 1 \, , \, b(y):=
|y|^{-d-1} \mbox{ if }  |y|>1 \, .
\]
If $|x-y|2^n>1$ we integrate $(d+1)$-times by parts in $\tilde{\xi}$, and if $|w-F(y)|2^\ell>1$ we integrate $(d+1)$-times by parts in $\tilde{\eta}$. Hence, we arrive at
the following formula for	$( J_{\sigma,n}^{\tau,\ell}a)(x)$:
\begin{align*}
\int_{y\in K}
\int_{\supp\Psi_{\sigma,1}\times\supp{\widetilde{\Psi}}'_{\tau, 1}}
\int_{\real^d }
\frac{\widetilde{f}_{r-1,n,\ell}(\tilde{\xi},\tilde{\eta},y)}
{2^{\max\{n,\ell\}(r-1)}}\,  b_n(x-y)b_\ell(w-F(y))a_{\tau,\ell}(w)\, \dd w\,  \dd\tilde{\xi}\, \dd\tilde{\eta} \,\dd y\, ,
\end{align*}
where $b_m(w)=b(2^m w)$ ($m\ge 0$), and $\widetilde{f}_{r-1,n,\ell}(\tilde{\xi},\tilde{\eta},y)$ is uniformly bounded by $C_2'(F,f)$.
Thus, there exists  $C_3<\infty$ such that for all $x\in \real^d$
	\begin{align}\label{eq:l:multop:V}
		|(J_{\sigma,n}^{\tau,\ell}a)(x)|\le {C_3}\, C_2(F,f)\,
 2^{-\max\{n,\ell\}(r-1)}\bigl (b_n\ast 
(b_\ell\circ F)\ast
|a_{\tau,\ell}|\bigr ) (x)\, , \, \mbox{ if } (\tau,\ell)\not\hookrightarrow_K(\sigma,n) \, .
\end{align}
Since 
$r-1>t-s'> 0$ and $c(\sigma)\le t$, $c'(\tau)\ge s'$,
there exists $\epsilon >0$ such that for
all $\sigma$ and $\tau$,   
\begin{align}\label{eq:l:multop:VI}
	2^{(c(\sigma)+\epsilon)n-c'(\tau)\ell-\max\{n,\ell\}
(r-1)}&\le 2^{(t+\epsilon)n-s'\ell-\max\{n,\ell\}(r-1)}\le 2^{-\epsilon\ell }\, , \forall n\ge 1\, , \forall \ell \ge 1 \, .
\end{align}
We can assume $n \cdot\ell\ne 0$ since if $n=0$ or $\ell=0$  then $\xi$ or $\eta$ is bounded. (By Footnote~\ref{sixteen}, we have  $n \cdot\ell\ne 0$ in our application.) 
Hence  starting with   the triangle inequality, then using  \eqref{eq:l:multop:V}, , we find $C(\epsilon)<\infty$ such that for all $a$
\begin{align*}
		&
\|Q_{\not\hookrightarrow_K}^{\Op}a\|_{L_p(\real^d, \ell_2^c)}
\le	C(\epsilon) \sup_{\sigma,n}2^{(c(\sigma)+\epsilon)n}
\bigl \|{\Psi}^\Op_{\sigma,n}
\sum_{(\tau,\ell)\not\hookrightarrow_K (\sigma,n)}
\MM_{F,f}({\widetilde{\Psi}'^\Op}_{\tau,\ell}
a_{\tau,\ell})\bigr \|_{L_p(\real^d)}\nonumber	\\
&\le C(\epsilon) \sup_{\sigma,n}
\sum_{(\tau,\ell)\not\hookrightarrow_K(\sigma,n)}
2^{(c(\sigma)+\epsilon)n-c'(\tau)\ell}2^{c'(\tau)\ell}
\|{\Psi}^\Op_{\sigma,n}
\MM_{F,f}(\widetilde{\Psi}'^\Op_{\tau,\ell}a_{\tau,\ell})\|_{L_p}\\
		&= C(\epsilon)
\frac{2^{d(n+\ell)}}{(2\pi)^{2d}}\, \sup_{\sigma,n}
\sum_{(\tau,\ell)\not\hookrightarrow_K (\sigma,n)}
2^{(c(\sigma)+\epsilon)n-c'(\tau)\ell}2^{c'(\tau)\ell}
\|J_{\sigma,n}^{\tau,\ell}(a)\|_{L_p}\, .
\end{align*}
Applying Young's inequality (with $\|b_m\|_{L_1}=2^{-dm}\|b\|_{L^1}$, for $m=\ell$ and $n$), and  \eqref{eq:l:multop:VI}  then yields finite constants $C_4$,  $C_5$, $C_6$ (depending on $\epsilon$, $F$ and $f$) such that, again for all $a$,
\begin{align*}
	&
\|Q_{\not\hookrightarrow_K}^{\Op}a\|_{L_p(\real^d, \ell_2^c)}\\
	&	\le		 C_4\sup_{\sigma,n}
\sum_{\tau,\ell}
2^{(c(\sigma)+\epsilon)n-c'(\tau)\ell
-\max\{n,\ell\}(r-1)}2^{(n+\ell)d}2^{c'(\tau)\ell}
\|b_n\ast
 (b_\ell\circ F)\ast a_{\tau,\ell}\|_{L_p}\\
	&	\le C_5 \sum_{\tau,\ell}2^{-\ell\epsilon}
2^{c'(\tau)\ell}\|a_{\tau,\ell}\|_{L_p}\\
&\le C_6\sup_{\tau,\ell}  2^{c'(\tau)\ell}\|a_{\tau,\ell}\|_{L_p}\le C(F,f)\, \|a\|_{L_p(\real^d, \ell_2^{c'})} \, .
	\end{align*}
\end{proof}

We end this subsection with  a bound on the transfer operator 
$\MM_{F,f}$  from \eqref{eq:to:loc}.

\begin{lemma}[Bounding the Local Transfer Operator]\label{l:to:loc:bound}
Let $K\subset\real^d$ be compact with nonempty interior.  Fix $p\in(1,\infty)$ and
fix $s$, $q$, $t$  as in \eqref{regulS-sqt}. Let $t-(r-1)<s'<0<q'\le t'$ satisfy $s'<s$, $q'\le q$ and $t'<t$.
Then there exists $C<\infty$ such that  for  any  cone-hyperbolic $C^r$-diffeomorphism $F$ from $\Theta'$ to $\Theta$ on $K$, taking the covering $\widetilde \Phi'$  given by Lemma~\ref{forcomp}, we have for 
any $f:\real^d\rightarrow\complex$ in $C^{r-1}_0(K)$ and all $v \in W^{s,t,q}_{p, \Theta,F(K)}$
\[	\|
\MM_{F, f}v\|_{W^{s,t,q}_{p, \Theta}}
\le C \bigl ( C(F,f)
\|v\|_{W^{s',t',q'}_{p, \Theta',F(K)}}
+C_0(F,f)
\|v\|_{W^{q}_{p, \Theta',F(K)}}
+C_\pm(F,f)\|v\|_{W^{s,t,q}_{p, \Theta',F(K)}}\bigr ) \, ,\]
where  $C(F,f)$ is  from Lemma~\ref{l:bound:II},
and,  recalling $|F|_{+}$, $|F|_{0}$, and $|F^{-1}|_{-}$ from \eqref{covv},
\begin{align*}
C_0(F,f)&= 
\sup_{K}\biggl |\frac{f}{|\det \D F|^{\frac 1p}}\biggr| \cdot
\max\{1,
|F|_{0}^q\}\, , \,\, 
	C_\pm(F,f)=  \sup_{K}
\biggl |\frac{f}{|\det \D F|^{\frac 1p}}\biggr | \cdot
\max\{|F|_{+}^t,|F^{-1}|_{-}^{|s|}\}\,  .
\end{align*}
\end{lemma}

\begin{proof}
Let $v\in W^{s,t,q}_{p,\Theta,F(K)}\subset W^{s',t',q'}_{p,\Theta,F(K)}$. For  $\tau\in\{\pm, 0\}$,  $\ell \ge 0$, set
$a_{\tau,\ell}=(\Psi'_{\tau,\ell})^\Op v$.  Recalling
\eqref{eq:reg:I},
we have $a=a(v)\in L_p(\real^d,\ell_2^c)\subseteq L_p(\real^d,\ell_2^{c'})$, and more precisely
	\begin{align*}
		\|a\|_{L_p(\real^d,\ell_2^c)}
= 
\|v\|_{W^{s,t,q}_{p,\Theta',F(K)}}\, ,\, 
\|a\|_{L_p(\real^d,\ell_2^{c'})}=\|v\|_{W^{s',t',q'}_{p,\Theta',F(K)}}\, , \,\,
 \bigl \| \bigl (\sum_{\ell}4^{q\ell}
\bigl | a_{0,\ell}\bigr |^2\bigr )^{\frac 12} \bigr \|_{L_p}=\|v\|_{W^{q}_{p,\Theta',F(K)}}\, .
\end{align*}
Letting  $\hookrightarrow_K$ be as in Definition~\ref{d:arrow}, the decomposition \eqref{eq:multop} gives
\begin{align*}
&\|\MM_{F,f}v\|_{W^{s,t,q}_{p,\Theta,K}}
=\|{Q}^\Op_{\not\hookrightarrow,K}a+{Q}^\Op _{\hookrightarrow,K}a\|_{L_p(\real^d, \ell_2^c)}
\le
\|{Q}^\Op_{\not\hookrightarrow,K}a\|_{L_p(\real^d, \ell_2^c)}+\|{Q}^\Op_{\hookrightarrow, K}a\|_{L_p(\real^d, \ell_2^{c})} \, .
\end{align*}
We conclude using Lemmas~\ref{l:bound:I} and ~\ref{l:bound:II}.
\end{proof}

\subsection{Lasota--Yorke Type Bounds on the Transfer Operator}

For $s<0<t$ and $\alpha>0$, set
\begin{align}\label{eq:expconst}
	\lambda^{(s,t,\alpha)}(x)\coloneqq\max
\bigl \{
\|(\D g_{-\alpha})^{\tr}|_{E_{+,x}^*}\|^{t},
\|(\D g_{\alpha})^{\tr}|_{E_{-,g_{-\alpha}(x)}^*}\|^{-s}\bigr \}\, ,
\quad x \in M\,  .
\end{align}
There exists $C'<\infty$ such that $\sup_x\lambda^{(s,t,\alpha)}(x)\le C' \theta^{\min \{t, |s|\}\alpha}$  by property \eqref{eq:split:ineq}.

\begin{lemma}[Bounding the Transfer Operator]\label{l:to:bound}
Fix $p\in(1,\infty)$ and fix $s$, $q$, $t$ as in \eqref{regulS-sqt}. Let
$t-(r-1)<s'<s<0<q< t'<t$.
There exist  $A=A(X,V)<\infty$ and $C=C(X,V)<\infty$, such that for all $\varphi\in W_p^{s,t,q} (M)$ 
\begin{align*}
\|\LL_{\alpha,V}\varphi\|_{W_p^{s,t,q}}
\le C e^{A\alpha}\|\varphi\|_{W_p^{s',t',q}}
 + C
\|\phi_{\alpha}|\det\D g_{-\alpha}|^{-\frac 1p}\cdot\lambda^{(s,t,\alpha)}\|_{L_\infty}
\|\varphi\|_{W_p^{s,t,q}}\, , \forall \alpha \ge 0\, .
\end{align*}
\end{lemma}

The bound in the above lemma shows that $\LL_{\alpha,V}$ is an operator semigroup on $W_p^{s,t,q} (M)$.
As usual for flows, however, it is not a true Lasota--Yorke bound since $W_p^{s,t,q}(M)$ is not compactly embedded in $W_p^{s',t',q'}(M)$ if $q'=q$.
However, using Lemma~\ref{l:norm:loc:derivative}, we will prove in Theorem~\ref{t:res:lasota} that the resolvent $(z-V-X)^{-1}$
satisfies a Lasota--Yorke bound.

\begin{proof}
If $\alpha<\alpha_0$, using $\int_0^{\alpha_0}=\int_0^{\alpha_0-\alpha}
+\int_{\alpha_0-\alpha}^{\alpha_0}$, we find $
\|\LL_{\alpha,V}\varphi\|_{W_p^{s,t,q}}\le \|\varphi\|_{W_p^{s,t,q}}
+\|\LL_{\alpha_0,\phi_{\alpha_0}}\varphi\|_{W_p^{s,t,q}}$.

We may assume from now on that $\alpha \ge \alpha_0$.
Recalling the charts $\kappa_\omega: V_\omega \to \real^d$, the partition of unity $\vartheta_\omega$, and the  cone systems
$\Theta_\omega$ (from Lemma~\ref{l:excone})  above  Definition~\ref{d:bs:ani}, write,  as before, 	
$$
V_{\alpha,\omega\omega'}=
V_\omega \cap g_\alpha(V_{\omega'})\, , \mbox{ and }
F_{-\alpha,\omega\omega'}(x)=\kappa_{\omega'} \circ g_{-\alpha} \circ \kappa_\omega^{-1}(x)\, ,
\,\, x\in \kappa_{\omega}(V_{\alpha,\omega\omega'})\, ,
\alpha \ge \alpha_0\, ,\, \omega, \omega'\in \Omega
\, .
$$
Since $\alpha \ge \alpha_0$,     each $F_{-\alpha,\omega\omega'}$ is cone-hyperbolic
from $\Theta_{\omega'}$ to
$\Theta_{\omega}$ on $\kappa_{\omega}(V_{\alpha,\omega\omega'})$. 

The intersection multiplicity of a family of sets is the maximal number of sets having nonempty intersection,
while the intersection multiplicity of a family of functions is the
intersection multiplicity of the family of the supports of the
functions.

We claim that there exists an integer $\nu_d\ge 2$, depending only on $d$, such that the following holds: There exist $C=C_{\alpha_0}<\infty$ and, for each 
$\alpha \ge \alpha_0$, a finite refinement $\WW_{\alpha}=\{W_{\alpha, \vec \omega}\}_{\vec \omega \in \Omega_{\alpha}}$
 of  $\VV_\alpha=\{V_{\alpha,\omega\omega'}\}_{(\omega,\omega') \in \Omega^2}$, of  intersection multiplicity at most $\nu_d$,  such that 
\begin{align}\label{eq:boundinf}
\sup_{W} 
|\phi_\alpha \cdot |\det\D g_{-\alpha}|^{-\frac 1p}|&
\le C \inf_{W}
|\phi_\alpha \cdot |\det\D g_{-\alpha}|^{-\frac 1p}|\, ,\,\,
\forall\,  W\in \WW_{\alpha}\, .
\end{align}
Indeed, since there exists $K_\alpha<\infty$
such that $\sup_{\beta \in [0,\alpha]}d(g_{-\beta}(x), g_{-\beta}(y))\le K_\alpha d(x,y)$,
while $\log \phi_\alpha=\int_0^\alpha V(g_{-\beta}(x)) \dd \beta$,
and (noting that $\alpha-\alpha_0 [\alpha/\alpha_0]\in [0, \alpha_0)$)
$$\log |\det\D g_{-\alpha}|=\log |\det\D g_{-(\alpha-\alpha_0 [\alpha/\alpha_0])}|
\circ g_{-([\alpha/\alpha_0]-1) \alpha_0}+\sum_{\ell=0}^{[\alpha/\alpha_0]-1}
\log |\det\D g_{-\alpha_0}| \circ g_{-\ell \alpha_0} \, ,
$$ 
where $V$ and $|\det \D g_{-\alpha_0}|$, $|\det \D g_{-(\alpha-\alpha_0 [\alpha/\alpha_0])}|$    are uniformly continuous on $M$
(they are in fact $\gamma$-H\"older for $\gamma=\min\{r-1,1\}$), there exists  a finite refinement\footnote{Using the bounded  distortion argument
 for hyperbolic maps  \cite[Proposition 20.2.6]{Katok1995}, one can construct
 $\widetilde \VV_\alpha$  by first iterating $[\alpha/\alpha_0]$ times the
 cover $\{V_\omega\}$ and then  refining to guarantee that the
 diameter in the last coordinate in charts is $<\alpha^{-1/\gamma}$. The cardinality of
such $\widetilde \VV_\alpha$  grows like $C\alpha^{1/\gamma}C^\alpha$, for
$C>1$. This is not needed.}
$\widetilde \VV_\alpha$ of  $\VV_{\alpha}$ such that \eqref{eq:boundinf} holds for all $W\in \widetilde \VV_\alpha$. 
 A finite refinement $\WW_\alpha$ of $\widetilde \VV_\alpha$ satisfying the claimed intersection multiplicity
bound can then be obtained e.g. by covering $M$ with $d$-dimensional balls of radius the Lebesgue number of $\widetilde \VV_\alpha$ centered on an appropriate lattice, see e.g. \cite[Footnote 19 p. 46]{Bbook}. 
(Note that the cardinality of $\widetilde \VV_\alpha$ or $\WW_\alpha$  is immaterial in view of the use of the reconstitution Lemma~\ref{partgroupdag} below.)
Finally,   fix   a $C^{r}$  partition\footnote{We can ensure  that $\alpha^{-1/\gamma}\|\vartheta_{\alpha,\vec \omega}\|_{C^{r-1}}$ is controlled by  the largest expansion  of $F_{-\alpha,\omega\omega'}$. This is
not needed.}  of unity
$\{\vartheta_{\alpha,\vec \omega}\}_{\vec \omega \in \Omega_{\alpha}}$ of $M$, 
subordinate to the cover $\WW_{\alpha}$, with intersection multiplicity
at most $\nu_d$.

After these preliminaries, we perform the estimate: Let $\varphi\in{W_p^{s,t,q}}(M)$, by definition, we have
$$
\|\LL_{\alpha,V}\varphi\|_{{W_p^{s,t,q}}}^2
\le \# \Omega \max_{\omega\in\Omega}
		\int_0^{\alpha_0}
\|(\vartheta_\omega\cdot(\LL_{\alpha',V}\circ \LL_{\alpha,V}\varphi))\circ  \kappa_\omega^{-1}\|_{W_{p,\Theta_\omega}^{s,t,q}}^2
\, \dd\alpha'\, .
$$
Next, using $\LL_{\alpha',V}\circ \LL_{\alpha,V}=\LL_{\alpha,V}\circ\LL_{\alpha',V}$,
and setting $\varphi_{\alpha'}=\LL_{\alpha',V}(\varphi)$,
we find
\begin{align}
\nonumber \|(\vartheta_\omega
&\cdot(\LL_{\alpha',V}\circ \LL_{\alpha,V}\varphi))\circ \kappa_\omega^{-1}\|_{W_{p,\Theta_\omega}^{s,t,q}}\\
&=
\|\sum_{\omega'\in\Omega}
\sum_{\vec \omega \in \Omega_{\alpha}} 
(\vartheta_\omega\vartheta_{\alpha,\vec \omega}\cdot \phi_{\alpha})\circ \kappa_\omega^{-1}\cdot
(\vartheta_{\omega'}\cdot\varphi_{\alpha'})\circ \kappa_{\omega'}^{-1}\circ F_{-\alpha,\omega\omega'}\|_{W_{p,\Theta_\omega}^{s,t,q}}\nonumber\\
&\nonumber \le C \nu_d^{(p-1)/p} \max_{\omega'\in\Omega}\biggl(\sum_{\vec \omega \in \Omega_{\alpha}}
\|
(\vartheta_\omega\vartheta_{\alpha,\vec \omega}\cdot \phi_{\alpha})\circ \kappa_\omega^{-1}\cdot
(\vartheta_{\omega'}\cdot\varphi_{\alpha'})\circ \kappa_{\omega'}^{-1}\circ F_{-\alpha,\omega\omega'}\|_{W_{p,\Theta_\omega}^{s,t,q}}^p\biggr )^{1/p}\\
\label{RHS2}&\qquad + \widetilde C_{\vartheta_\alpha}  \max_{\omega'\in\Omega}\sum_{\vec \omega \in \Omega_{\alpha}}
\|
(\vartheta_\omega\vartheta_{\alpha,\vec \omega}\cdot \phi_{\alpha})\circ \kappa_\omega^{-1}\cdot
(\vartheta_{\omega'}\cdot \varphi_{\alpha'})\circ \kappa_{\omega'}^{-1}\circ F_{-\alpha,\omega\omega'}\|_{W_{p,\Theta_\omega}^{s',t',q'}}\, ,
\end{align}
using the fragmentation Lemma~\ref{partexpdag0}.
By Lemma~\ref{l:to:loc:bound}  the term in the last line of \eqref{RHS2} is bounded by
\begin{align}\label{11}
\widetilde C_{0,\alpha} (X,V) \max_{\omega'\in \Omega} 
\|
(\vartheta_{\omega'}\varphi_{\alpha'})\circ \kappa_{\omega'}^{-1}\|_{W_{p,\Theta_{\omega'}}^{s',t',q'}}\, ,
\end{align}
for $\widetilde C_{0,\alpha}(X,V)<\infty$.
Remark~\ref{forLeibniz}, gives systems $\Theta'_{\omega}<\Theta_\omega$ (independent of
$\alpha$)
such that $F_{-\alpha,\omega\omega'}$ is  cone-hyperbolic from $\Theta'_{\omega'}$ to $\Theta_{\omega}$ on $\kappa_{\omega}(V_{\alpha,\omega\omega'})$. 
For  $\alpha\ge \alpha_0$ and $\vec \omega\in \Omega_{\alpha}$,
let $\tilde \vartheta_{\alpha,\vec \omega}:M\to [0,1]$ be  $C^{r-1}$,
 supported in $W_{\alpha,\vec \omega}$,
and such that 
$\tilde \vartheta_{\alpha,\vec \omega} \vartheta_{\alpha,\vec \omega}= \vartheta_{\alpha,\vec \omega}$, and  set 
\begin{equation}f_{\alpha, \vec\omega}=(\vartheta_\omega\tilde \vartheta_{\alpha,\vec \omega}\phi_\alpha )\circ \kappa_\omega^{-1}\, ,\,\,
\bar\vartheta_{\alpha, \vec \omega}=(\vartheta_\omega \vartheta_{\alpha,\vec \omega})\circ \kappa_\omega^{-1}\circ F_{-\alpha,\omega,\omega'}^{-1}
\end{equation}
 Then, Lemma~\ref{l:to:loc:bound} gives $C_p <\infty$ and $\widehat C_{0,\alpha} (X,V)<\infty$ such that   each term in the sum on the second-to-last line of
\eqref{RHS2} is bounded by
\begin{align}\label{recc}
& C_p C_\pm
(F_{-\alpha,\omega\omega'},f_{\alpha,\vec \omega})
\|\bar \vartheta_{\alpha,\vec \omega}\cdot (
(\vartheta_{\omega'}\cdot\varphi_{\alpha'})\circ \kappa_{\omega'}^{-1})
\|_{W_{p,\Theta'_{\omega'}}^{s,t,q}}^p
\\
\nonumber &\qquad\qquad\qquad\qquad\qquad\qquad\qquad
+\widehat C_{0,\alpha} (X,V)\|
(\vartheta_{\omega'}\cdot\varphi_{\alpha'})\circ \kappa_{\omega'}^{-1}
\|_{W_{p,\Theta'_{\omega'}}^{s',t',q}}^p\, .
\end{align}
Due to the strict inequality between cone ensembles, the reconstitution Lemma~\ref{partgroupdag}
bounds the $p$th root of the sum of  \eqref{recc} over $\vec \omega$, uniformly in $\alpha$. 
Recalling \eqref{11}, taking the square, the maximum over $\omega$, and integrating over $\alpha'$, we find  $\bar C<\infty$ and $C'_\alpha=C'_\alpha(X,V)<\infty$ such that
\begin{align}\label{44}
\|\LL_{\alpha,V} (\varphi)\|^2_{{W_p^{s,t,q}}} 
\le {C}'_{\alpha}\cdot 
\| \varphi\|^2_{W_p^{s',t',q}}+ \bar C  \max_{\omega,\omega',\vec \omega}
(C_\pm(F_{-\alpha,\omega\omega'},
f_{\alpha, \vec \omega}))^2  \cdot
\| \varphi\|^2_{W_p^{s,t,q}}\, ,\,\forall \alpha\ge \alpha_0\, . 
\end{align}

We next estimate $C_\pm(F_{-\alpha,\omega\omega'},f_{\alpha,\vec\omega})$.
By construction of  $\Theta_\omega$ in the proof of Lemma~\ref{l:excone},  
and since the covering $\widetilde \Phi'$ from Lemma~\ref{forcomp} used in
Lemma~\ref{l:to:loc:bound} can be taken  such that  $\supp\widetilde \Phi_\sigma$ is bounded away from $E^*_\tau$ (in charts)
if $\tau\ne \sigma$,  there exists $C<\infty$ such that,
recalling  \eqref{covv}, we have for all $\alpha\ge \alpha_0$,
all $\omega$, $\omega'$, and all $\vec \omega \in \Omega_{\alpha}$, setting $K_{\alpha, \vec \omega}=\kappa_\omega(\supp \tilde \vartheta_{\alpha, \vec \omega})$,
\begin{align*}
|F_{-\alpha,\omega\omega'}|_{+,K_{\alpha, \vec \omega}}
&\le C \sup_{x\in K_{\alpha, \vec \omega}}
\bigl \|(\D g_{-\alpha})^{\tr}|_{E_{+,x}^*} \bigr \| \, ,\,\,\,
|F_{-\alpha,\omega\omega'}^{-1}|_{-,K_{\alpha, \vec \omega}}
\le C \sup_{x\in K_{\alpha, \vec \omega}}
\bigl \|(\D g_{\alpha})^{\tr}|_{E_{-,g_{-\alpha}(x)}^*} \bigr \|\, .
\end{align*}
Thus, using \eqref{eq:boundinf} and  $\inf |\psi_1|\sup |\psi_2|\le \sup |\psi_1 \psi_2|$ for continuous $\psi_1$, $\psi_2$, we find $C<\infty$ such that
\begin{align*}
\nonumber &C_\pm(F_{-\alpha,\omega\omega'},f_{\alpha,\vec\omega})\le C\max_{W\in \WW_{\alpha}}
\bigl (
\sup_W|\phi_\alpha \cdot |\det\D g_{-\alpha}|^{-\frac 1p}|
\cdot \sup_W |\lambda^{(s,t,\alpha)}|\bigr )\\
&\qquad\le \bar C\max_{W\in \WW_{\alpha}}
\bigl (\inf_W \biggl  |\frac{\phi_{\alpha}}{ |\det\D g_{-\alpha}|^{\frac 1p}}\biggr |
\cdot \sup_W|\lambda^{(s,t,\alpha)}|\bigr )
\le  \widehat C
\sup_M \biggl |\frac{\phi_\alpha}{|\det\D g_{-\alpha}|^{\frac 1p}}\lambda^{(s,t,\alpha)}\biggr |\, , \, \forall \alpha \ge \alpha_0\, .		
\end{align*}
In view of \eqref{44}, we have proved the lemma.
\end{proof}

\noindent Strong continuity  suffices (it is not necessary
\cite{Butterley_2016}) to show that $X+V$ is the generator  of $\LL_{\alpha,V}$:

\begin{lemma}[Strong Continuity. Domain of the Generator $X+V$]\label{l:to:C0}
Let $p\in(1,\infty)$ and  fix $s$, $q$, $t$ as in \eqref{regulS-sqt}.
The  family 
$\{\LL_{\alpha,V}\}_{\alpha \ge 0}$ of bounded operators on
$W_p^{s,t,q}(M)$
forms a strongly continuous semigroup.
The generator of this semigroup is
$X+V\colon D(X+V)\rightarrow W_p^{s,t,q}(M)$,
which is  closed  on its (dense) domain $D(X+V)\subseteq W_p^{s,t,q}(M)$.  
Moreover,
if $q<r-2$ or if $\phi_\alpha$ is $C^{r}$ in the flow direction, setting $\DD_{r-1}:={C^{r-1}(M)}$
if $q<r-2$, and otherwise 
$$\DD_{r-1}:= {C^{r-1, r}(M)}=\{\varphi \in C^{r-1}(M)\mid
\varphi \mbox{ is } C^r \mbox{ in the flow direction} \}\, ,
$$
$D(X+V)$ contains\footnote{As observed in \cite{Butterley_2016}   strong continuity implies  that the completion $\DD$ of
$\DD_0=\{\int_0^\beta \LL_{\alpha,V} \varphi \, \dd \alpha \mid \varphi\in W_p^{s,t,q}(M)\, ,\,\,\beta >0\}$ 
under $\|\cdot \|_{W_p^{s,t,q}(M)}$ is a dense subset of
$W_p^{s,t,q}(M)$, so that $\DD=W_p^{s,t,q}(M)$. Clearly, $\LL_{\alpha,V}(\DD_0) \subset \DD_0$ and $\DD_0\subset D(X+V)$. Thus, $\DD_0$ is a dense subset of $D(X+V)$ for the graph norm, without any conditions on $q$ or $\phi_\alpha$.}
$\DD_{r-1}$ as a dense subset
for the graph norm
$\|\cdot\|_{W_p^{s,t,q}(M)}+
\|(X+V)(\cdot)\|_{W_p^{s,t,q}(M)}$.
\end{lemma}

\begin{proof} 
To establish strong continuity, it suffices to show  $\lim_{\alpha\downarrow 0}
\|\LL_{\alpha}\varphi-\varphi\|_{W_p^{s,t,q}(M)}=0$ for all $\varphi\in W_p^{s,t,q}(M)$ (\cite[Proposition I.1.3]{Engel_2006}).
 Lemmas~\ref{l:to:bound} and \ref{l:contembed} give $C<\infty$ such that $\|\LL_{\alpha,V}\varphi\|_{W_p^{s,t,q}}
\le C\|\varphi\|_{W_p^{s,t,q}}$ for all  $\alpha\in [0,1]$.
By density of $C^\infty(M)$, for every $\epsilon>0$ there is $\tilde{\varphi}=\tilde \varphi_\epsilon\in{C^{r-1,r}(M)}$ 
 such that
$\|\varphi-\tilde{\varphi}\|_{W_p^{s,t,q}}\le\epsilon$.
Therefore,
\begin{align}\label{eq:to:C0:I}
\|\LL_{\alpha,V}\varphi-\varphi\|_{W_p^{s,t,q}}&\le 
\|\LL_{\alpha,V}(\varphi-\tilde{\varphi})\|_{W_p^{s,t,q}}
+\|\varphi-\tilde{\varphi}\|_{W_p^{s,t,q}}
+
\|\LL_{\alpha,V}\tilde{\varphi}-\tilde{\varphi}\|_{W_p^{s,t,q}}\nonumber\\
	&\le (C+1)\epsilon 
+ \|\LL_{\alpha,V}\tilde{\varphi}-\tilde{\varphi}\|_{W_p^{s,t,q}}\, ,\, \, \forall \epsilon>0\, , \forall  \alpha\in[ 0,1]\, .
\end{align}
Since $\tilde \varphi\in{C^{r-1,r}(M)}$ (if $0<q<r-2$ then  the argument  can be adapted to $\tilde \varphi \in C^{r-1}$) we have
 $\partial_{\alpha'}\LL_{\alpha',V}
{\tilde \varphi}\in C^{r-1}(M)$. Thus,  there exists $C(\tilde \varphi_\epsilon)<\infty$ such that
\begin{align}\label{eq:to:C0:II}
\|\LL_{\alpha,V}\tilde{\varphi}
-\tilde{\varphi}\|_{W_p^{s,t,q}}
=\bigl \|\int_0^\alpha
\partial_{\tilde \alpha}\LL_{\tilde \alpha,V}{\tilde \varphi}\, 
\dd \tilde \alpha\bigr \|_{W_p^{s,t,q}}
\le \alpha\sup_{0\le \tilde \alpha\le \alpha}
\|\partial_{\tilde \alpha}
\LL_{\tilde \alpha,V}{\tilde \varphi}\|_{W_p^{s,t,q}}
\le C(\tilde \varphi_\epsilon)\alpha\, .
	\end{align}
Combining \eqref{eq:to:C0:I} and\eqref{eq:to:C0:II} gives $\lim_{\alpha\downarrow 0}
\|\LL_{\alpha}\varphi-\varphi\|_{W_p^{s,t,q}(M)}=0$. Also, $\lim_{\alpha\downarrow 0}\alpha^{-1}
(\LL_{\alpha}\varphi-\varphi)=X\tilde \varphi + V\tilde \varphi$ for $\tilde \varphi\in{C^{r-1,r}(M)}$.
Strong continuity and \cite[Thm II.1.4]{Engel_2006} then imply that $X+V$ is the generator of the semigroup, and that it is   closed with domain dense in $W_p^{s,t,q}(M)$. 
Clearly, 
$\LL_{\alpha,V}({C^{r-1}(M)})\subseteq{C^{r-1}(M)}
$, and, if $\phi_\alpha$ is $C^r$ in the flow direction, then $\LL_{\alpha,V}({C^{r-1,r}(M)})\subseteq{C^{r-1,r}(M)}
$.
The final  claim thus follows from  \cite[Proposition II.1.7]{Engel_2006}, since $\LL_{\alpha, V}(\DD_{r-1})\subset \DD_{r-1}$
and 
$\lim_{\alpha \to 0}\alpha^{-1} ( \LL_{\alpha, V}\varphi -\varphi)$ exists for $\varphi \in \DD_{r-1}$
in the two cases considered.
\end{proof}

\subsection{Lasota--Yorke Bounds for the Resolvent. Discrete Spectrum of $X+V$}\label{unn}
Recall $\lambda^{(s,t,\alpha)}$ from \eqref{eq:expconst}. We set
(the limit below exists and is finite by superadditivity)
\begin{align}\label{eq:res:lasota:A1}
	\lambda_{\min}=\lambda_{\min}^{s,t,p}(X,V)&:= \lim_{\alpha\to\infty}
\frac 1\alpha\log
\|\phi_\alpha
|\det\D g_{-\alpha}|^{-\frac1p}
\lambda^{(s,t,\alpha)}\|_{L_\infty(M)}\, .
\end{align}

 Recalling $A(X,V)$ from Lemma~\ref{l:to:bound},
we may and shall replace  $A(X,V)$ by $\max \{A(X,V), \lambda_{\min}\}$ 
from now on. 
The following theorem will furnish an essential spectral bound for $X+V$:

\begin{theorem}[Lasota--Yorke Inequality for the Resolvent]\label{t:res:lasota}
Let $p\in (1,\infty)$ and let $s$, $q$, $t$ be as in \eqref{regulS-sqt}.
Let $t-(r-1)<s'<0<q'\le t'$ satisfy $s'<s$, $t'<t$, and $q-1\le q'<q$.
For any $\epsilon>0$, there exists $C<\infty$ such that for  all
$\delta>0$, all $z\in\complex$ with $\Re z>A(X,V)+\delta$,  all $n\in\naturall$,  and all $\varphi\in{W_p^{s,t,q}}(M)$,
recalling our notation $\RR_z=(z-X-V)^{-1}$,
\[
\delta \|\RR_z^{n+1}\varphi\|_{{W_p^{s,t,q}}}
\le \frac{C(|z|+1)}{(\Re z-A(X,V))^{n}}
\|\varphi\|_{{W_p^{s',t',q'}}}+ \frac{C}
{(\Re z-\epsilon-\lambda_{\min}^{s,t,p}(X,V))^{n}}\|\varphi\|_{{W_p^{s,t,q}}}\, .
\]
\end{theorem}

The above theorem implies that the  spectral radius of the resolvent $\RR_z$ on $W_p^{s,t,q}(M)$ is bounded by $|\Re z-A(X,V)|^{-1}$ if $\Re z>A(X,V)$ (a very rough bound). In addition, we have:

\begin{corollary}[Essential Spectral Radius]\label{c:res:spectrum}
For all $z\in\complex$ with $\Re z>A(X,V)$,
the essential spectral radius of  $\RR_z$ on $W_p^{s,t,q}(M)$ is bounded by
 $|\Re z-\lambda_{\min}^{s,t,p}(X,V)|^{-1}$.  Moreover, the set
$\{\lambda\in\sigma (X+V)|_{W_p^{s,t,q}(M)} 
\mid \Re \lambda>\lambda_{\min}^{s,t,p}(X,V)\}$
consists of isolated eigenvalues of finite multiplicity.	
\end{corollary}

\begin{proof}
Since the inclusion $W_p^{s,t,q}(M)\subset W_p^{s',t',q'}(M)$ is compact by Lemma~\ref{l:compact},
the first claim  follows from a result of Hennion \cite[Corollaire 1]{Hennion_1993} and Theorem~\ref{t:res:lasota}.
 The second claim then follows 
 from the spectral mapping theorem\footnote{This spectral mapping theorem says that if
$z \notin \sigma(X+V)$ then $\sigma(\RR_z)\setminus \{0\}=\{(z-\lambda)^{-1}\mid \lambda \in \sigma(X+V)\}$, where $(z-\lambda)^{-1}$ is an eigenvalue if and only if $\lambda$ is
an eigenvalue.} for the resolvent \cite[Thm V.1.13]{Engel_2006}.
\end{proof}

If $\lambda_{\max}^{s,t,q,p}(X,V)>\lambda_{\min}^{s,t,p}(X,V)$, where
\begin{align*}
\lambda_{\max}^{s,t,q,p}(X,V):=\sup \Re \, \sigma (X+V)|_{W_p^{s,t,q}(M)}
 \, , \end{align*}
then
the isolated eigenvalues furnished by Corollary~\ref{c:res:spectrum}
are called the {\it Ruelle--Pollicott resonances} of $X+V$
on $W_p^{s,t,q}(M)$. We will apply the following theorem to our scale $W_{p}^{s,t,q}(M)$
and, in Lemma~\ref{l:invmeas}, to the scale from \cite{Giulietti_2013}:

\begin{theorem}[Intrinsicness of Ruelle--Pollicott Resonances]\label{intr}
Let $\BB_1$ and $\BB_2$ be two Banach spaces of distributions on $M$ on which
$\{\LL_{\alpha, V}\}$  is a strongly continuous semigroup with generator $X+V$.
Assume that both $\BB_1$ and $\BB_2$ 
 contain $C^{r-1}(M)$ as a dense subset and are  continuously embedded in the dual of
$C^{r-1}(M)$. If  there exists $\lambda_{\min}>-\infty$ such that the sets
$D_i=\{\lambda\in\sigma (X+V)|_{\BB_i} 
\mid \Re \lambda>\lambda_{\min}\}$, $i=1,2$, 
consist of isolated eigenvalues of finite multiplicity, then  $D_1=D_2$,
including multiplicities. In particular, the corresponding generalised eigenvectors belong to
$\BB_1\cap\BB_2$ (in fact, to the intersection of the domains of $X+V$ on
$\BB_1$ and on $\BB_2$).
\end{theorem}

\begin{proof}
If $r=\infty$, this is  a special case of \cite[Thm 2.3]{Go2}, which refers to \cite[Thm 1.5]{FauSj}.
If $r<\infty$\footnote{We expect that intrinsicness can also be proved by using
dynamical determinants.}  the proof of \cite[Thm 1.5]{FauSj} using meromorphic extensions
of suitable resolvents applies,
replacing $L_2(M)$ by the dual of $C^{r-1}(M)$ and
using that $C^{r-1}(M)$ is a dense subset of both $\BB_1$ and $\BB_2$.
 \end{proof}
 
\noindent Lemma~\ref{l:invmeas} says  that $\lambda_{\max}^{s,t,q,p}(X,V)=\htop$ (for
suitable $s,t,q$,  $p$) for
$V$ as in 
Section~\ref{s:horo}.

The remainder of \S\ref{unn} is devoted to the proof of Theorem~\ref{t:res:lasota}.
Since the resolvent can be written as a Laplace transform (integrating along the flow),
this proof   will follow from 
the flow box condition \eqref{eq:chart:flowbox},  Lemma~\ref{l:to:bound}, and the lemma 
below:
\begin{lemma}[Integration Along the Flow]\label{l:norm:loc:derivative}
Fix $p\in (1,\infty)$. There exists $C<\infty$ such that 
\begin{align*}
\|(\sum_{n=0}^\infty4^{\tilde{r}n}
|{\Psi}^\Op_{0,n}v|^2)^{\frac 12}\|_{L_p}\le C
\|(\sum_{n=0}^\infty4^{(\tilde{r}-1)n}
|{\Psi}^\Op_{0,n}\partial_{x_d}v|^2)^{\frac 12}\|_{L_p} \, ,
\forall \tilde{r}>0\, .
	\end{align*}	
\end{lemma}

\noindent (Adapting the proof  gives $
\|(\sum_{n=0}^\infty4^{(\tilde{r}-1)n}
|{\Psi}^\Op_{0,n}\partial_{x_d}v|^2)^{\frac 12}\|_{L_p}\le C \|(\sum_{n=0}^\infty4^{\tilde{r}n}
|{\Psi}^\Op_{0,n}v|^2)^{\frac 12}\|_{L_p}$.)
	
\begin{proof}
It is enough to consider the terms with $n>0$.
The starting point is 
	$${\Psi}^\Op_{0,n}(\partial_{x_d}v)=
(\FFF^{-1}{\Psi}_{0,n})*(\partial_{x_d}v)=
(\partial_{x_d}\FFF^{-1}{\Psi}_{0,n})* v=2^n({\DD_d}^\Op v)_n
\, , \forall v\in L_p(\real^d)\, , 
$$
where 
$(\DD_d(\xi)b)_n:= \iunit\frac {\xi_d}{2^n}\Psi_{0,n}(\xi)b$, for $n\in\naturall$, $\xi\in\real^d$, and $b\in\complex$.

For  a sequence $a$ of complex numbers with 
$\|a\|_{\ell_2(\tilde r)}:=\bigl (\sum_{n=1}^\infty 4^{\tilde{r}n}|a_n|^2\bigr )^{1/2}<\infty$,  we put 
\[(Q_d(\xi)a)_n:=-\iunit\frac {2^n}{\xi_d}\widetilde{\Psi'}_{0,n}(\xi)a_n\, ,\quad \xi\in\real^d\, ,\, n\in\naturall\, ,\]	
where $\widetilde{\Psi'}_{0,n}$ is associated via \eqref{eq:suppext}
to a covering $\widetilde{\Psi'}$ of $\Theta$ with $\supp\widetilde{\Psi'}_0$
contained in a cone around $\xi_d$. Then $({Q_d}^\Op {\DD_d}^\Op v)_n={\Psi}^\Op_{0,n}v$. Since  there exists $\tilde \gamma_0<\infty$ such that
$2^{n-1}\le |\xi| \le \tilde \gamma_0 |\xi_d|$ for any $\xi\in \supp \widetilde \Psi'_{0,n}$, the map $Q_d$ satisfies the decay condition  in Theorem~\ref{l:multiplier}. Hence, taking  $\HH_1=\HH_2=\{ a\mid \|a\|_{\ell_2(\tilde r)}<\infty\}$,
 the map ${Q_d}^\Op$ is a bounded linear operator  on $L_p(\real^d,\ell_2(\tilde r))$. This concludes the proof, since it gives $C<\infty$ such that
$
\|\|{Q_d}^\Op {\DD_d}^\Op v\|_{\ell_2(\tilde r)}\|_{L_p}\le 
C\|\|{\DD_d}^\Op v\|_{\ell_2(\tilde r)}\|_{L_p}$.
\end{proof}

\begin{proof}[Proof of Theorem~\ref{t:res:lasota}]
By Lemma~\ref{l:to:bound}, for  $z\in\complex$ such that $\Re z>A(X,V)$
(\cite[Cor. II.1.11]{Engel_2006}),
\begin{align}\label{eq:res:int}
\RR^n_{z}\varphi= \int_{0}^\infty \frac{\alpha^{n-1}e^{-z\alpha}}{(n-1)!}\LL_{\alpha,V}\varphi\, \dd\alpha\, ,  \, \forall n \in \naturall \, ,
\end{align}
for all $\varphi\in{W_p^{s,t,q}}(M)$.
Introducing the truncated iterated resolvent
\begin{equation}\label{truncres}
\RR_{tr,z}^{n}\varphi:=\int_0^{\alpha_0} \frac{\alpha^{n-1}e^{-z\alpha}}{(n-1)!}\LL_{\alpha,V}\varphi\, \dd\alpha \, ,
\end{equation}
we claim that
\begin{equation}\label{truncbound}
\|\RR_{tr,z}^{n}\varphi\|_{W_p^{s,t,q}}\le \frac{C}{
(\Re z+\Delta)^{n}}\|\varphi\|_{W_{p}^{s,t,q}} \,,\,\,  \forall \Delta \ge 0
\, ,\, \forall \Re z>0 \, , \,\forall n>e\cdot \alpha_0\cdot(\Re z+\Delta)\, .
\end{equation}
This bound  holds because, using Lemma~\ref{l:to:bound},
$\sup_{\alpha \in [0, \alpha_0]}
e^{-\Re z\alpha}\le 1$, and  $\int_0^{\alpha_0} \frac{\alpha^{n-1}}{(n-1)!} \dd\alpha=\frac {\alpha_0^n}{n!}
\le \frac{1 }{(\Re z+\Delta)^n} $
if $n>e\cdot \alpha_0\cdot(\Re z+\Delta)$
(recall that $n! \ge n^n e^{-n}$),
we find
$$
\int_0^{\alpha_0} \frac{\alpha^{n-1}e^{-\Re z\alpha}}{(n-1)!} \|\LL_{\alpha,V}\varphi \|_{W_{p}^{s,t,q}}\, \dd\alpha
\le 
C \frac{\alpha_0^n}{n!}\|\varphi\|_{W_{p}^{s,t,q}}
\le \frac{C\|\varphi\|_{W_{p}^{s,t,q}}}{(\Re z+\Delta)^n} \, ,  \, \forall \Re z>0\, , \, \forall n >e \alpha_0 (\Re z+\Delta) \, .
$$
We can therefore focus on times $\alpha\ge \alpha_0$ in \eqref{eq:res:int} and invoke Remark~\ref{forLeibniz}.

 Lemma~\ref{l:to:bound} gives $C_1=C_1(\epsilon)<\infty$ such that for all $n\in \naturall$
	\begin{align}
		\|\RR^{n+1}_z\varphi\|_{W_p^{s,t,q}}&
\le\int_0^\infty 
\frac{\alpha^{n-1}e^{-\Re z\alpha}}{(n-1)!}
\|\LL_{\alpha,V}\RR_z\varphi\|_{W_p^{s,t,q}}\dd\alpha\nonumber\\
		\label{eq:res:bound:III}&
\le\frac{C_1}{(\Re z-A(X,V))^n}
\|\RR_z\varphi\|_{W_{p,\Theta'}^{s',t',q}}+ 
\frac{C_1}
{(\Re z				-\epsilon-\lambda_{\min})^{n}}
\|\RR_z\varphi\|_{W_p^{s,t,q}}\, ,
	\end{align}
where,  for $\Theta'_\omega<\Theta_\omega$  as in Remark~\ref{forLeibniz}, we replaced $\Theta$ by $\Theta'$
in the first term of \eqref{eq:res:bound:III}.
Lemma~\ref{l:to:bound} also gives  $C_2<\infty$ such that 
	\begin{align}\label{eq:res:bound}
\|\RR_z\varphi\|_{W_p^{s,t,q}}\le \frac {C_2}{\Re z-A(X,V)}\|\varphi\|_{W_p^{s,t,q}}\, ,
	\end{align}
so that the second term of \eqref{eq:res:bound:III} is bounded as claimed. The starting point to bound the first term is the fact that the flow box condition \eqref{eq:chart:flowbox} gives
$(\D\kappa_{\omega}^{-1})(\partial_{x_d}|_{\kappa_\omega(V_\omega)})=X|_{V_\omega}$, and hence
\begin{align}\label{eq:res:bound:0}
\partial_{x_d}((\vartheta_\omega\cdot\tilde \varphi)\circ \kappa_\omega^{-1})
&=((X\vartheta_\omega)\cdot\tilde \varphi+\vartheta_\omega
\cdot(X\tilde \varphi))\circ\kappa_{\omega}^{-1}\, .
\end{align}

Using the triangle inequality 
(and $-\infty <q'$) to separate the contribution
of $\Theta'_{\omega,0}$, and 
applying Lemma~\ref{l:norm:loc:derivative}  with \eqref{eq:res:bound:0} (for $\tilde \varphi=\LL_{\alpha',V}\RR_{z}\varphi$) to bound this term, we find $C_3<\infty$ such that
\begin{align*}
\|(\vartheta_\omega\cdot&\LL_{\alpha',V}\RR_{z}\varphi)\circ \kappa_\omega^{-1}\|_{W_{p,\Theta'_\omega}^{s',t',q}(K_\omega)}
\le 
\|(\vartheta_\omega\cdot \LL_{\alpha',V} \RR_z \varphi)\circ \kappa_\omega^{-1}\|_{W_{p,\Theta'_\omega}^{s',t',q'}(K_\omega)}\\
\nonumber&+C_3
\|(( X \vartheta_\omega)\cdot {\LL_{\alpha',V}\RR_{z}\varphi})\circ \kappa_\omega^{-1}\|_{W_{p,\Theta'_\omega}^{q-1}(K_\omega)}+C_3
\|(\vartheta_\omega\cdot X  \RR_{z}\LL_{\alpha',V}\varphi)\circ \kappa_\omega^{-1}\|_{W_{p,\Theta'_\omega}^{q-1}(K_\omega)}\, ,
\nonumber
\end{align*}	
where we used for the  last  term that \eqref{eq:res:int} implies $\LL_{\alpha',V}\RR_{z}\varphi=\RR_{z}\LL_{\alpha',V}\varphi$.
Since 
$( X \vartheta_\omega)\circ\kappa_\omega^{-1}=\partial_{x_d} (\vartheta_\omega \circ \kappa_\omega^{-1})\in  C^{r-1}_0(\kappa_\omega(V_\omega))$ 
(using that $\vartheta_\omega$ and $\kappa_\omega$ are $C^r$, with $\vartheta_\omega$ is compactly supported
in $V_\omega$) 
and $q-1\le q'$,  Lemma~\ref{l:to:loc:bound} for the identity map gives  $C_4<\infty$ such that
\begin{align*}
\|(( X \vartheta_\omega)\cdot \LL_{\alpha',V}\RR_{z}\varphi)\circ \kappa_\omega^{-1}\|_{W_{p,\Theta'_\omega}^{q-1}(K_\omega)}
\le C_4\sup_{\omega\in\Omega}
\|(\vartheta_\omega\cdot {\LL_{\alpha',V}\varphi})\circ \kappa_\omega^{-1}\|_{W_{p,\Theta_\omega}^{s',t',q'}(K_\omega)}\, .
\end{align*}
Using  $X \RR_z\varphi=z\RR_z \varphi-V\RR_z\varphi-\varphi
$, and, again, $\RR_{z}\LL_{\alpha',V}\varphi=\LL_{\alpha',V}\RR_{z}\varphi$, we find
\begin{align}\nonumber
\|(\vartheta_\omega\cdot  X\RR_z&\LL_{\alpha',V}\varphi) \circ\kappa_\omega^{-1}\|_{W_{p,\Theta'_\omega}^{q-1}(K_\omega)}
\le |z|
\|(\vartheta_\omega\cdot \LL_{\alpha',V}\RR_z\varphi)\circ \kappa_\omega^{-1}\|_{W_{p,\Theta'_\omega}^{s',t',q-1}(K_\omega)}\\
\label{eq:res:bound:II}		
&+
\|(\vartheta_\omega\cdot V\LL_{\alpha',V}\RR_z\varphi)\circ \kappa_\omega^{-1}\|_{W_{p,\Theta'_\omega}^{s',t',q-1}(K_\omega)}
+
\|(\vartheta_\omega\cdot \LL_{\alpha',V}\varphi)\circ \kappa_\omega^{-1}\|_{W_{p,\Theta'_\omega}^{s',t',q-1}(K_\omega)}
\, .
\end{align}
Since $V\in C^{r-1}(M)$, we may bound the first term  in \eqref{eq:res:bound:II}	 with Lemma~\ref{l:to:loc:bound} for the identity map.
Using the definition of $\|\RR_z\varphi\|_{{W_p^{s',t',q}}}$
as an integral of local norms over $\alpha'\in[0,\alpha_0]$, this bounds the first  term of \eqref{eq:res:bound:III} as claimed
(we use \eqref{eq:res:bound} for the terms with $\RR_z\varphi$).
\end{proof}


\subsection{Dolgopyat Bounds for the Resolvent of  Weighted Transfer Operators}

In the contact Anosov case  and for the potential
$V$ introduced in the next section, the spectrum
of $X+V$ has already been studied \cite{Giulietti_2013},
on a different Banach space.  We will use in the proof of Lemma~\ref{l:invmeas} that
the discrete spectra of $X+V$ on our Banach spaces and the spaces of \cite{Giulietti_2013} coincide
in a big enough half-plane (``intrinsicness''), but it is not clear how to exploit this to obtain the bounds on the resolvent needed in Section~\ref{s:horo}.  Indeed, 
 in the self-adjoint
case, there exist good bounds on the iterated resolvent $\RR_{z}^n$ in terms of the distance between
$z$ and the spectrum. 
However, even when $W_p^{s,t,q}(M)$ is a Hilbert space, the operator $X+V$ is not self-adjoint a priori,
so the existence of a spectral gap for $X+V$ does {\emph not} imply good
 bounds on the resolvent in general
(see \cite{Tsujii_2010, FT3,CG} for special cases where such bounds are known).
For this reason, we introduce the following condition:

\begin{condition}[Weak Dolgopyat Bounds on the Resolvent]\label{cnd:wA}
There exist  $p\in (1, \infty)$, $s$, $q$, $t$ as in \eqref{regulS-sqt}, constants
$$
 s''\in \real\, ,\quad \delta'\in (\lambda_{\min}^{s,t,p},\lambda^{s,t,q,p}_{\max})\, ,
$$
and constants $a_0>0$, $b'_0>1$,  $c_1<1<C_1$,  such that,
for all $a\ge a_0$ and 
$\gamma'$ in the range
\begin{equation}
aC_1<\gamma' < \frac{c_1}{\log \bigl (1+(\lambda^{s,t,q,p}_{\max}-\delta')/a\bigr )}\, , 
\end{equation}
we have 
\[
\|\RR_{a+\iunit b+\lambda_{\max}^{s,t,q,p}}^{n}\varphi \|_{W^{s''}_p}\le C_1|a+(\lambda_{\max}^{s,t,q,p}-\delta')|^{-n}
\|\varphi\|_{C^1}\, ,\,\,
\forall |b|\ge b'_0\, ,  
\mbox{ where }n=\lceil\gamma'  \log |b|\rceil\, .\]
\end{condition}

Using  the mollification
ideas introduced by Liverani   in \cite[\S5,\S7]{Baladi2012} and \cite[\S 9]{Baladi_2018}, we will show that Condition~\ref{cnd:wA}  implies norm estimates on the resolvent (see \cite[Remark 2.6]{Butterley_2016}):

\begin{proposition}[Strong Dolgopyat Bounds on the Resolvent]\label{cnd:A} If there exists $C_0< \infty$ with 
\begin{equation}\label{exact}
\|\LL_{\alpha,V}\|_{W_p^{s',t',q'}}\le C_0 e^{\lambda^{s,t,q,p}_{\max}\alpha}\, ,\,\, 
\forall \alpha\ge 0\, ,
\end{equation}  
for  $t-(r-1)<s'<0<q'<t'$, with $t-t'=q-q'=s-s'>0$, for  some $s$ , $q$ , $t$ and $p>\max\{d/t,d/(r-1+s)\}$ such that  Condition~\ref{cnd:wA} holds
for $c_1$, $C_1$, then  there exist
$\delta\in (\lambda_{\min}^{s,t,p},\lambda^{s,t,q,p}_{\max})$, $a>0$, $b_0>1$, $C<\infty$,   and
$\gamma\in (aC_1,c_1/\log \bigl (1+(\lambda^{s,t,q,p}_{\max}-\delta)/a \bigr ))$, 
with
\begin{equation}\label{cnd:sA}
\|\RR_{a+\iunit b+\lambda_{\max}}^{n}\|_{{W_p^{s,t,q}}}\le C|a+(\lambda_{\max}^{s,t,q,p}-\delta)|^{-n}\, ,\,\,
\forall |b|\ge b_0\, , 
\mbox{ where }n=\lceil\gamma  \log |b|\rceil\, .
\end{equation}
\end{proposition}

Before proving Proposition~\ref{cnd:A}, we make further remarks.
Bounds~\eqref{cnd:sA} on suitable anisotropic Banach spaces are used in many places in the literature,
starting with Liverani's breakthrough paper \cite{Liverani_2004} (see e.g. \cite{Baladi2012, Giulietti_2013, Baladi_2018}). This has been axiomatized by Butterley in \cite{Butterley_2016}:  Together with a weak Lipschitz control on $\alpha \mapsto \LL_{\alpha,V}$, the bounds~\eqref{exact} and ~\eqref{cnd:sA}  imply
the 
\textit{spectral gap}  property\footnote{\label{fqc}Beware that this property
does not imply a spectral gap (quasi-compactness) for the time-one transfer operator:
we do not expect  $\LL_{\alpha, V}$ to be eventually norm continuous \cite[Thm \S II.5.3]{Engel_2006}, so a priori we only have $\sigma (\LL_{\alpha, V})\subset \exp(\alpha \sigma (X+V))$ for  $\alpha\ge 0$ (equality
holds for eigenvalues and residual spectrum), see \cite[\S V.2.b]{Engel_2006}. A spectral gap
for the time-one transfer operator is only known in special cases,
\cite{Tsujii_2010, FT3}.}
\begin{align}
\nonumber
	\sigma (X+V)|_{W_p^{s,t,q}(M)} \cap \{ \Re \lambda > \delta\} \mbox{ is a finite set.}
\end{align}
The above spectral gap   can be used to get exponential decay of
correlations, but   the implication in the other direction 
is not known in general. (Dolgopyat \cite{Dolgopyat_1998}
obtained exponential decay of correlations for Gibbs measures
with arbitrary H\"older potentials for geodesic flows on surfaces of strictly negative curvature or, more generally, $C^5$ Anosov flows such that $E_-$ and $E_+$ are $C^1$ and not jointly integrable,  using symbolic dynamics. His ideas  led to results of Liverani on the SRB measure
of contact Anosov flows \cite{Liverani_2004}.
See \cite{GSt} and \cite{TsZh}, and references therein, for recent sufficient conditions  ensuring
 exponential mixing for Gibbs measures and Anosov flows.)

The bounds~\eqref{cnd:sA} have been established
\cite{Liverani_2004, Tsujii_2010, Baladi2012, Giulietti_2013} 
for the generator $X$ associated to contact Anosov flows
and  the  potential $V=0$, replacing  our spaces $W_p^{s,t,q}$ by other
anisotropic Banach spaces. 
For the potential $V$ used in Section~\ref{s:horo}, Dolgopyat bounds are shown  in \cite[\S 7]{Giulietti_2013}
(see also the argument sketched by Faure and Guillarmou before \cite[Proposition 3.4]{Faure_2017_2}). 
We expect that \eqref{cnd:sA} or Condition~\ref{cnd:wA}
can be proved directly in our setting. For our purposes it is sufficient
instead to refer to  \cite[\S 7]{Giulietti_2013} in Section~\ref{s:horo} to establish Condition~\ref{cnd:wA}, and then invoke Proposition~\ref{cnd:A}.
We thereby illustrate how to  build bridges
between results for different anisotropic spaces (once the essential radius
is controlled, exact growth is obtained, and, for the Dolgopyat estimate, mollification bounds are known).

\begin{proof}[Proof of Proposition~\ref{cnd:A}]
Let $\{\Theta'_\omega\}$ and $\{\Theta_\omega\}$ form an adapted pair for $\AAA$ and $g_\alpha$ in the sense of Remark~\ref{forLeibniz}. 
Denote by 
$\|\varphi\|_{W_p^{s',t',q'}(\Theta')}$ the norm constructed with $\Theta'_\omega$ instead of $\Theta_\omega$.
We start with three trivial but useful observations:
First, for any  $\beta>0$,  $\delta_2 >0$, and $\delta_1\ge 0$, we have for all $|b|\ge 1$
and $a >\delta_1$ that
\begin{equation}\label{triviall}
|a-\delta_1|^{-\lceil\gamma'  \log |b|\rceil} |b|^{-\beta} \le  |a+\delta_2|^{-\lceil\gamma'  \log |b|\rceil}\, ,\,\,   
\forall \gamma' \in \bigl (0,
\frac{\beta}{\log \bigl (1+\delta_2/a)-\log \bigl (1-\delta_1/a)}\bigr )\, .
\end{equation}
(If $\delta_1=0$ and $\beta>0$, taking $\delta_2>0$ small enough, we can choose $\gamma'$ arbitrarily large in \eqref{triviall}.)

Second, for any  $\beta'>0$,   and $\delta_3>\delta_2> 0$, we have for all $|b|\ge 1$
and $a >0$ 
\begin{equation}\label{triviall'}
|a+\delta_3|^{-\lceil\gamma'  \log |b|\rceil} |b|^{\beta'} \le  |a+\delta_2|^{-\lceil\gamma'  \log |b|\rceil}\, ,\,\,   
\forall \gamma' >\frac{\beta'}{\log  (1+\delta_3/a)-\log  (1+\delta_2/a)}\, .
\end{equation}

Third, for any $0<\delta_0<\delta_2$
 and  $\delta_1\ge 0$, 
we have for all $a>\delta_1$ and $m_1, m_2\in \naturall$,
\begin{equation}\label{triviall2}
|a+\delta_2|^{-m_1} 
|a-\delta_1|^{-m_2}  
\le  |a+\delta_0|^{-m_1-m_2}\,\,\,
\,\,
\mbox{ if }\,\,
\frac{m_1}{m_2}
\ge \frac{\log (1+\delta_0/a)
-\log(1-\delta_1/a)
}
{\log (1+\delta_2/a)-\log(1+\delta_0/a)} \, .
\end{equation}
(If $\delta_1=0$, for fixed $\delta_2>0$, taking $\delta_0>0$ small enough, we can choose $m/n$ arbitrarily small.)

Set $\lambda_{\max}=\lambda^{s,t,q,p}_{\max}$.
To deduce ~\eqref{cnd:sA} from Condition \ref{cnd:wA}, we
use the Lasota--Yorke estimate: We may assume that  $s''<\min\{ -d-1,s\}$. Then,
for any
$
\tilde \delta_2\in (0, \lambda_{\max}-\lambda_{\min}^{s,t,p})$,
Theorem~\ref{t:res:lasota} with \eqref{trick} give  $C$, $C(s'',m)$,  and $A(X,V)\ge \lambda_{\max}$,
such that for all $\varphi \in W_p^{s,t,q}(M)$ and all $m, n\in \naturall$, 
\begin{align}
\nonumber
\|\RR_{a+\iunit b+\lambda_{\max}}^{m+1+n}&\varphi\|_{{W_p^{s,t,q}}}
\le \frac{2C(s'',m)|b|}{(a+\lambda_{\max}-A(X,V))^{m}}
\|\RR_{a+\iunit b+\lambda_{\max}}^{n}\varphi\|_{{W_p^{s'',s'',s''}}}
\\ &
\label{line0}+ \frac{C}
{(a+\tilde \delta_2)^{m}}\|\RR_{a+\iunit b+\lambda_{\max}}^{n}\varphi\|_{{W_p^{s,t,q}}}\, ,
 \forall |b|\ge 1 \,,\forall a>A(X,V)-\lambda_{\max}+1 \, .
\end{align}
Then,  we proceed as  in \cite[\S5,\S7]{Baladi2012} or \cite[\S 9]{Baladi_2018}: First,
since \eqref{exact} 
gives $\|\RR_{a+\iunit b+\lambda_{\max}}^{n}\varphi\|_{{W_p^{s,t,q}}}\le C a^{-n}\|\varphi\|_{{W_p^{s,t,q}}}$,   for any  $\epsilon_0>0$ (to be fixed in the next paragraph and by \eqref{cheps}) there exist
 $m$ and $n$ with $m\le \epsilon_0 n$, 
 and
such that the last term  on the right-hand side of \eqref{line0} satisfies the required condition:
Indeed, apply \eqref{triviall2} 
for $m_1=m$, $m_2=n$, $\delta_1=0$, $\delta_2=\tilde \delta_2$,
taking  $\delta_0>0$ small enough such that $m\le \epsilon_0 n$ is allowed,
and choose $\delta \ge \min(\epsilon_1,\lambda_{\max}-\delta_0)$. 

For the first term on the right-hand side of
\eqref{line0}, it is enough
to bound  $\frac{2C(s'',m)|b|}{(a+\lambda_{\max}-A(X,V))^{m}}
\|\RR_{z}^{n}\varphi-\RR^n_{tr,z}\varphi\|_{{W_p^{s'',s'',s''}}}$. Indeed, it is not hard to see that there exists $\bar C(s'')\ge 1$ such that $C(s'',m)\le C(s'')^m$. Let $\delta_1(s'')>A(X,V)-\lambda_{\max}>0$ be such that 
$$
\frac{2C(s'',m)}{(a-A(X,V)-\lambda_{\max})^{m}}
\le \frac{2C(s'')^m}{(a-A(X,V)-\lambda_{\max})^{m}}\le \frac C{(a-\delta_1(s''))^{m}}\, ,\,\,
\forall m\ge 1\,,\,\,\,
\forall  a>\delta_1(s'')
\, .
$$
Then, the contribution
of   $\|\RR^n_{tr,z}\varphi\|_{{W_p^{s'',s'',s''}}}\le \|\RR^n_{tr,z}\varphi\|_{{W_p^{s,t,q}}}$
is controlled by
\eqref{truncbound}, applying  \eqref{triviall2} for
$m_1=n$, $m_2=m$, 
$\delta_1=\delta_1(s'')> 0$, and  $\delta_2=\Delta>\delta_0>0$,  for large enough $\Delta$, 
taking $\epsilon_0$ small enough so that
\begin{equation}
n\ge m \frac{\log (1+\delta_0/a)-\log(1-\delta_1(s'')/a)}
{\log (1+\Delta/a)-\log(1+\delta_0/a)}\, ,
\end{equation}
 and taking $b_0$ large enough to ensure
$n=\lceil\gamma'  \log |b|\rceil>e\alpha_0(a+\lambda_{\max}+\Delta)$ if $|b|\ge b_0$,  for $\gamma'\ge C_1a$ determined below.
(Again,   choose $\delta \ge \min(\epsilon_1,\lambda_{\max}-\delta_0)$.)

Set $\RR^n_{*,z}:=\RR_{z}^{n}\varphi-\RR^n_{tr,z}$. 
Fixing $s'$, $q'$, $t'$ with  $q'-q=t'-t=s'-s<0$, 
and $t-(r-1)<s'<0<q'<t'$, we decompose,  for any $s''\le s'$,
\begin{equation}
\label{lafin}
|b|\|\RR_{*,z}^{n}\varphi\|_{{W_p^{s''}}} 
\le |b| \|\RR_{*,z}^{n}(\MMM_\epsilon \varphi)\|_{{W_p^{s''}}}+ C|b|
\|\RR_{*,z}^{n}(\varphi-\MMM_\epsilon\varphi)\|_{{W_p^{s',t',q'}}}\, ,
\end{equation}
where  $\MMM_\epsilon$ is the mollification
operator in charts defined by \eqref{defmoll}, for $\epsilon=|b|^{-\kappa}$, with
$\kappa>1$ to be chosen later.
Let $\{\Theta'_\omega\}$ form an adapted pair with $\{\Theta_\omega\}$ 
By ~\eqref{exact}  we have
\begin{equation}
\label{ani0}
\|\RR_{*,a+\iunit b+\lambda_{\max}}^{n}(\varphi-\MMM_\epsilon\varphi)\|_{{W_p^{s',t',q'}}}\le 
\frac{C}{ a^n} \|\varphi-\MMM_\epsilon\varphi
\|_{W_p^{s',t',q'}(\Theta')}
\, .
\end{equation}

Then the mollification estimate  Lemma~\ref{moll} gives
\begin{equation}
\label{ani}
\frac{C}{ a^n}\|\varphi-\MMM_\epsilon\varphi\|_{{W_p^{s',t',q'}(\Theta')}}\le 
\frac{C}{ a^n} \epsilon^{s-s'}\|\varphi\|_{{W_p^{s,t,q}}}\le
\frac{\bar  C}{ a^n} |b|^{-\kappa (s-s')}\|\varphi\|_{{W_p^{s,t,q}}}\, ,
\end{equation}

If $\kappa>1/(s-s')$, applying \eqref{triviall} with  $\beta=\kappa(s-s')-1>0$ and $a>\delta_1=0$, $\delta_2=\lambda_{\max}-\delta$, the bounds \eqref{ani0} and \eqref{ani} take care of the second term in the right-hand side of \eqref{lafin}, assuming
\begin{equation}\label{uba}
\gamma' <\frac{\kappa(s-s') -1}{\log \bigl (1+\frac { \lambda_{\max}-\delta}a)}
\, .
\end{equation}
Note that  this inequality is compatible with $\gamma' >aC_1$ if $\kappa$
is large enough. 

Fix $\eta_0\in (0,\min\{t, r-1+s\})$, small. By the Sobolev embeddings for $W^{1+\eta_0}_p=F^{1+\eta_0}_{p,2}$ and
$B_{\infty,\infty}^{1_+}$ \cite[Thm 2.2.3(i)]{RS} in dimension $d$, we have 
$$\|\tilde \varphi \|_{C^1}\le \widehat C  \|\tilde \varphi \|_{W_p^{1+\eta_0}}\, ,\,\mbox{ if } p>\frac{d}{\eta_0}\, .
$$
Thus, Condition \ref{cnd:wA} 
 bounds  the first term in the right-hand side of \eqref{lafin} by 
$$
\frac{C_1 |b|}{ |a+\lambda_{\max}-\delta'|^{n}}\|\MMM_\epsilon\varphi \|_{C^1}
\le \frac{\bar C |b|}{ |a+\lambda_{\max}-\delta'|^{n}}\|\MMM_\epsilon\varphi \|_{W_p^{1+\eta_0}} \, .
$$
Since the charts in $\AAA$ are $C^r$,
the classical  isotropic
 mollification estimate of \cite[Lemma 5.3]{Baladi2012} (replacing $X_0$  by
 $M$ and $2$  by $r$ there) becomes:
For each $p\in (1,\infty)$ and all $-r+1<s \le s'<r+s\le r$,
there exists $C_\#$ so that
for all small
enough $\epsilon >0$ and every $\varphi\in W^{s}_p(M)$, we have 
\begin{equation*}
\|\MMM_\epsilon (\varphi)\|_{W^{s'}_p(M)} 
\le C_\# \epsilon^{s-s'} \|\varphi\|_{W^{s}_p(M)}  \, .
\end{equation*}
Therefore, since $-r+1<s<0<1+\eta_0<r + s$, taking $s'=1+\eta$ and recalling
\eqref{compareiso}, we have
\begin{equation*}
\frac{\bar C |b|}{ |a+\lambda_{\max}-\delta'|^{n}}\|\MMM_\epsilon\varphi \|_{W_p^{1+\eta_0}}
\le \frac{\bar C |b|\, \epsilon^{s-1-\eta_0}}{ |a+\lambda_{\max}-\delta'|^{n}}\|\varphi \|_{W_p^{s,s,s}}
=\frac{\bar C\, |b|^{1+\kappa r}}{ |a+\lambda_{\max}-\delta'|^{n}}
\|\varphi \|_{W_p^{s,s,s}}\, .
\end{equation*}
We need to multiply the above by $C(a-\delta_1(s''))^{-m}$.
For this, we use  that 
\begin{equation}\label{conf}
\gamma'>\frac{1+\kappa r}{\log  (1+ (\lambda_{\max}-\delta')/a)-\log  (1+\delta_2/a)}
\end{equation} 
is compatible with the upper bound \eqref{uba} on $\gamma'$, up to taking small enough
 $\delta_2\in (0,\lambda_{\max}- \delta')$ in the right-hand side of \eqref{triviall'} for $\beta'=1+\kappa r$
and $\delta_3=\lambda_{\max}-\delta'$.

We conclude the proof of the
proposition by applying \eqref{triviall2} 
for $a>\delta_1=\delta_1(s'')$, $\delta_2$ as in the previous paragraph,  $m_1=n$, $m_2=m$, 
for $\delta_0\in(0,\delta_2)$,
and $\epsilon_0\in (0,1)$ such that 
\begin{equation}\label{cheps}
\frac{\log (1+\delta_0/a)
-\log(1-\delta_1(s'')/a)
}
{\log (1+\delta_2/a)-\log(1+\delta_0/a)}
<\frac 1 {\epsilon_0} \, .
\end{equation}
Indeed,  taking  $\delta>\epsilon_1$ closer to $\lambda_{\max}$ if necessary to ensure $\delta\le \lambda_{\max}-\delta_0$, and for $\gamma'>aC_1$ satisfying \eqref{uba}--\eqref{conf}, 
take $\gamma>0$ such that 
(using $ m\le \epsilon_0 n$)
$$ 
\gamma \lceil \log |b|\rceil=m+n \le (\epsilon_0+1)\gamma'\lceil \log |b|\rceil
\, . $$ 
Then  $\gamma <1/\log \bigl (1+\frac{\lambda_{\max}-\delta}{a}\bigr )$
follows from \eqref{uba}, up to taking $\delta<\lambda_{\max}$  closer to $\lambda_{\max}$.
\end{proof}

Remark~\ref{C.3} explains why using mollifiers  through isotropic spaces as in \cite[Lemma 5.4, (7.5)--(7.6)]{Baladi2012} does not allow to carry out successfully the bounds in the previous
proof.

\section{Asymptotics of Horocycle Integrals}\label{s:horo}
In this section,  we assume throughout that $r\ge 2$,   the $C^r$ Anosov flow $g_\alpha$ 
on $M$ is topologically mixing
with stable dimension $d_-=1$, and that the strong-stable distribution $E_-$ is orientable. 
\subsection{Horocycle Flow $h_\rho$. Horocycle Integral $\gamma_x(\varphi,T)$. Renormalisation Time $\tau(\rho,\alpha,x)$}

We shall focus on stable horocycle flows. Analogous results exist for unstable horocycle flows.

\begin{definition}[(Stable) Horocycle Flow]
A   (stable) horocycle flow 
for a topologically mixing $C^r$ Anosov flow $g_\alpha$ on $M$
with $d_-=1$ and $E_-$ orientable
is a $C^0$ flow $h_\rho$ on $M$  such that 
$\partial_\rho h_\rho\in E_-\setminus\{0\}$
for all $\rho \in \real$.
\end{definition}

\begin{remark}[Unit Speed Parametrisation]
The stable manifolds of  the flow $g_\alpha$ are the 
submanifolds  tangent to the  bundle $E_-$
(this bundle is in general only H\"older, existence is
ensured by the stable manifold theorem, see e.g. \cite[Thm 17.4.3]{Katok1995}).  
We can
 parametrise  stable manifolds by the
arc-length induced by the Riemannian metric on $M$. Since we assumed that $E_-$
is orientable, this defines uniquely a horocycle flow  with $|\partial_\rho h_\rho|\equiv 1$,  called the unit speed horocycle flow. All other horocycle flows are obtained by time reparametrisations.
Topological mixing of $g_\alpha$ implies that each  stable
 manifold is dense in $M$ \cite[p. 84]{Marcus_1977} so any horocycle flow is minimal. 
\end{remark}

Our main object of interest is the following (stable) horocycle integral:

\begin{definition}[Horocycle Integral]
The horocycle integral of the  horocycle flow $h_\rho$ for the observable $\varphi\in{C^{0}(M)}$ at $x\in M$ is defined by
\begin{align}\label{eq:cocycle:current}
\gamma_{x}(\varphi,T)=
\int_{0}^T \varphi\circ h_\rho(x)\dd\rho\, .
\end{align}
\end{definition}

Writing $\mu(\varphi)=\int \varphi \dd \mu$,  where
$\mu$ is the unique\footnote{See \cite{Bowen_1977} for a proof of unique ergodicity. If $g_\alpha$ preserves a smooth measure
see also \cite{Marcus_1975}. See also Remark~\ref{MME}.} $h_\rho$-invariant probability measure,  
we have
\begin{equation}\label{thm2.1}
\gamma_{x}(\varphi,T)= T \cdot \mu(\varphi)  +\EE_{T,x}(\varphi)\, ,\quad
\lim_{T\to \infty} \frac{\EE_{T,x}(\varphi)}{T}=0\, ,\forall x\in M\,,\,\,\,
\forall \varphi\in{C^{0}(M)}\, .
\end{equation}

Our main result, Theorem~\ref{t:decomposition} in \S\ref{55b}, gives a more precise asymptotic expansion, 
involving the spectrum and eigendistributions of a suitably weighted transfer operator $\LL_{\alpha,V}$. 
A crucial ingredient in our analysis is the renormalisation time (first introduced by Marcus   \cite[p.83]{Marcus_1977} to study ergodic properties of the horocycle flow):

\begin{definition}[(Pointwise) renormalisation time]
A map $\tau\colon \real^2\times M\rightarrow \real$ which satisfies 
\begin{align}\label{eq:lr}
		g_{\alpha}\circ h_{\rho}(x)=h_{\tau(\rho,\alpha,x)}\circ g_{\alpha}(x)\, ,
		\qquad \forall \rho,\alpha\in\real\, , \forall x\in M,
\end{align}
is called a (pointwise) renormalisation time for the stable horocycle flow $h_\rho$.
\end{definition}

For the unit speed horocycle flow of the geodesic flow on a compact surface of constant negative curvature, the renormalisation time is 
 $\tau(\rho,\alpha,x)=\rho\cdot \exp( -\alpha\cdot \htop)$. More  generally:

\begin{lemma}[Properties of $\tau(\rho,\alpha,x)$]\label{l:existlr}
There exists a	unique solution $\tau(\rho,\alpha,x)$ to \eqref{eq:lr}. In addition 
$\tau(\rho,\alpha,x)$ is differentiable in $\rho$, and we have\footnote{In \eqref{p:proplr:XII},  we denote by $(\partial_\rho h_0)^*\in E_-^*$ 
 the canonical dual of $\partial_\rho h_0:=\partial_\rho h_\rho{}|_{\rho=0}.$}
\begin{align}
\label{p:proplr:XI}&\tau(\rho,\alpha,x)
=\gamma_{x}(\partial_\rho\tau(0,\alpha,\cdot),\rho)\, ,\,\,\forall x\in M\,,\,\, \forall \alpha\in \real\,,\, \forall  \rho \in\real\, ,
	\\
&\label{p:proplr:XII}
\partial_{\rho}\tau(0,\alpha,x)=
\det \D g_{\alpha}|_{E_-}(x)  \cdot
\frac{(\partial_{\rho}h_{0}(x))^*
(\partial_{\rho}h_{0}(x))}{(\partial_{\rho}h_{0}\circ g_{\alpha}(x))^*(\partial_{\rho}h_{0}\circ g_{\alpha}(x))}\, ,\,\,\forall x\in M\,,\,\, \forall \alpha\in\real\, .
\end{align}
In particular, $\partial_{\rho}\tau(0,\alpha,x)>0$,
$\tau(0,\alpha,x)=0$,  $\tau(\rho,0,x)=\rho$. Moreover,
there exists $C<\infty$ with
\begin{align}
\label{p:proplr:IX}&
\frac{1}{C} \le \frac{\tau(\rho,-\alpha,x)}{\rho}e^{-\htop\alpha}\le C  \,,\, \forall x\in M\, ,\,
\forall \rho \in \real\,, \forall \alpha \ge 0  \mbox{ such that } | \rho|\ge 1 \, ,\\
\label{p:proplr:VIII}&  
\frac{1}{C}\le\frac{ \rho}{\tau(\rho,\alpha,x)}e^{-\htop\alpha}\le C \,,\, \forall x\in M\, ,\,
\forall \rho \in \real\,,\forall \alpha \ge 0 \mbox{ such that } |\tau(\rho,\alpha,x)|\ge 1\, .
\end{align}
\end{lemma}

The bounds \eqref{p:proplr:IX}--\eqref{p:proplr:VIII}
will come from \cite[App. C]{Giulietti_2013}. That
$\lim_{\rho \to  \infty}(\tau(\rho,\alpha,x)/\rho)=e^{-\alpha\htop}$ for all $\alpha\ge 0$
follows from \cite{Marcus_1977}, see the proof of Lemma~\ref{wow}.

The key fact behind our main result (Theorem~\ref{t:decomposition}) is the following 
consequence\footnote{The proof of \eqref{thekey} is a simplification of that of
Sublemma~\ref{sublemmarenorm} below.} of \eqref{eq:lr} 
\begin{equation}\label{thekey}
\int_{0}^T \varphi\circ h_\rho(x)\dd\rho=
\int_0^{\tau(T,\alpha,x)}
(\LL_{\alpha, V}\varphi)\circ h_\rho\circ g_\alpha(x)\dd \rho\, ,
\end{equation}
(``renormalisation'') where  the transfer
operator $\LL_{\alpha,V}$   is defined by \eqref{eq:to}, 
choosing 
\begin{align}
\label{eq:weight}
	V\equiv -\partial_{\alpha}\partial_{\rho}\tau(0,0,\cdot)\, ,\,\,\, \mbox{i.e. }	\phi_\alpha= \partial_\rho\tau(0,-\alpha,\cdot)\, ,
\end{align}
and assuming that $\phi_\alpha$ is $C^{r-1}$.
The underlying idea will be to take $\alpha=O(\log T)$ so that
$\tau(T,\alpha,x)=O(1)$, and then exploit the information on the spectrum of the semigroup
$\LL_{\alpha,V}$ obtained in the previous section.

Since the derivative of the Jacobian is the divergence and $\partial_\alpha g_{-\alpha}|_{\alpha=0}=X$,
we find for the unit speed horocycle flow  that \eqref{p:proplr:XII} implies
\begin{align*}
V= \divv ( X|_{E_-})
\mbox{ and  } \phi_\alpha=\det\D g_{-\alpha}{}|_{E^-}\, .
\end{align*}
Hence, if  $E_-$ is $C^{r-1}$  then  $\phi_\alpha\in{C^{r-1}(M)}$. 
More generally, if  $E_-$ is $C^{r-1}$,  for any $C^{r}$ time reparametrisation of the unit speed horocycle flow, the weight $\phi_\alpha$ is $C^{r-1}$ by \eqref{p:proplr:XII}.
  (Cf. \cite[Remark 2.4]{giul_liv_2017}.)
In order to fit in the Banach norm setting of Sections~\ref{sec2} and \ref{s:prop},
we will need $\phi_\alpha$ to be $C^{r-1}$
 for $r\ge 2$, and we will have to introduce the horocycle integrals
\eqref{eq:cocycle:current2} localised by smooth cutoff functions (following \cite{giul_liv_2017},
see the  proof of Lemma~\ref{l:horo:decomp:local}), replacing thus
\eqref{thekey} by the more involved version 
of ``renormalisation'' in Sublemma~\ref{sublemmarenorm}.

Before proving Lemma~\ref{l:existlr}, we state and prove a consequence
of \eqref{thekey} and classical results:

\begin{lemma}[The Invariant Measure $\mu$ as an Eigenvector]\label{wow}
If $\phi_\alpha$ from \eqref{eq:weight} is $C^{r-1}$, we have 
\begin{equation}
\label{l:invmeas:I} \mu(\LL_{\alpha,V} \varphi)=e^{\htop\alpha}\mu(\varphi)\, , \, \, \forall \alpha \ge 0 \,,\,\,\forall \varphi \in C^0(M)\, .
\end{equation}
\end{lemma}

\begin{remark}[Spectrum of $X+V$ on $L^1(\mu)$]\label{Rwow}Lemma~\ref{wow} gives
$
\mu(|\LL_{\alpha,V}\varphi|)\le
\mu(\LL_{\alpha,V}|\varphi|)=\mu(|\varphi|)$ for all $\alpha \ge 0$ and any $\varphi\in{C^{0}(M)}$.
Therefore, 
since $\mu$ is a Radon measure,  for each $\alpha \ge 0$, the
operator $\LL_{\alpha,V}$ is  bounded 
on the Banach space $L^1(\mu)$, with  norm equal to $e^{\alpha\htop}$. Hence,
using \cite[Cor. II.1.11]{Engel_2006} and  \eqref{eq:res:int}, the spectral radius of  $\RR_z=(z-(X+V))^{-1}$ on
 $L^1(\mu)$ is bounded by
$|\Re z-\htop|$ if $\Re z > \htop$.
The spectrum of $X+V$ on $L^1(\mu)$ thus lies
in the half-plane $\Re z \le \htop$.
\end{remark}

\begin{proof}[Proof of Lemma~\ref{wow}]
Unique ergodicity  \eqref{thm2.1} (twice),  renormalisation \eqref{thekey}, and
a result of Marcus \cite[Lemma~3.1, p~84]{Marcus_1977}  give
that for all $\alpha \ge 0$ and $\varphi\in{C^{0}(M)}$
\begin{align*}	
\mu(\varphi)&=\lim_{T\to\infty}\frac 1T\gamma_{x}(\varphi,T)=\lim_{T\to\infty}\frac{\tau(T,\alpha,x)}{T}\frac 1{\tau(T,\alpha,x)}
\gamma_{g_{\alpha(x)}} (\LL_{\alpha,V}\varphi, \tau(T,\alpha,x))=e^{-\alpha\htop}\mu(\LL_{\alpha,V}\varphi) .
\end{align*}
\end{proof}

\begin{proof}[Proof of Lemma~\ref{l:existlr}]
Since stable leaves are dense and the flow $h_\rho$ is non-singular, this flow does not admit any  periodic orbits.
For $x\in M$ and  $\rho, \alpha\in\real$,  set $h_{\alpha,\rho}(x)\coloneqq g_{\alpha}\circ h_{\rho}\circ g_{-\alpha}(x)$. Then $\partial_\rho h_{\alpha,\rho}\in E_{-,x}\setminus\{0\}$. Hence $h_{\alpha,\rho}(x)$ parametrizes the same stable manifold  as $h_{\rho}(x)$. If there were two different pointwise 
times $\tau$, there would exist $\rho_1<\rho_2\in\real$ such that
$h_{\alpha,\rho}(x)=h_{\rho_1}(x)=h_{\rho_2}(x)$, and this would contradict the absence of periodic orbits.
Thus,  $\tau(\rho,\alpha,x)$ is uniquely defined and differentiable in $\rho$.
 We deduce from \eqref{eq:lr} that
\begin{align*}
h_{\tau(\rho,\alpha_1+\alpha_2,\cdot)}(g_{\alpha_1+\alpha_2}(\cdot))&=g_{\alpha_1+\alpha_2}( h_\rho(\cdot))=g_{\alpha_1}( h_{\tau(\rho,\alpha_2,\cdot)}( g_{\alpha_2}(\cdot)))
=h_{\tau(\tau(\rho,\alpha_2,\cdot),\alpha_1,g_{\alpha_2}(x))}( g_{\alpha_1+\alpha_2}(\cdot)),
\end{align*}
and
\begin{align*}
h_{\tau(\rho_1+\rho_2,\alpha,x)}( g_\alpha(x))&=g_\alpha( h_{\rho_1+\rho_2}(x))
=g_\alpha( h_{\rho_1}\circ g_{-\alpha}\circ g_\alpha( h_{\rho_2}\circ g_{-\alpha}\circ g_{\alpha}(x))
\\&
=h_{\tau(\rho_1,\alpha, h_{\rho_2}(x))}( h_{\tau(\rho_2,\alpha,x)}( g_{\alpha}(x))\, .
\end{align*}
This implies that  for all $ \alpha_1,\alpha_2, \rho_1,\rho_2 \in\real$, 
\begin{equation}\label{p:proplr:II}
\tau(\rho,\alpha_1+\alpha_2,\cdot)=\tau(\tau(\rho,\alpha_2,\cdot),\alpha_1,g_{\alpha_2}(\cdot)) \, , \,\,
\tau(\rho_1+\rho_2,\cdot,x)=
\tau(\rho_1,\cdot,h_{\rho_2}(x))+\tau(\rho_2,\cdot,x)\, .
\end{equation}
Then, using  $\tau(0,\alpha,x)=0$,  and differentiating  the identities in \eqref{p:proplr:II}  at $\rho=0$ and $\rho_1=0$, we find 
\begin{equation}\label{p:proplr:V}
\partial_\rho \tau(\rho,\alpha,x)=
\partial_\rho \tau(0,\alpha,h_\rho(x))\, ,\,\,\,
\partial_\rho\tau(0,\alpha_1,g_{\alpha_2}(x))\partial_\rho\tau(0,\alpha_2,x)=
\partial_\rho\tau(0,\alpha_1+\alpha_2,x)\,, \,\forall  \alpha_i\, ,
\end{equation}
Next \eqref{p:proplr:XI} follows from the definition \eqref{eq:cocycle:current} of $\gamma_x$ and the first claim of \eqref{p:proplr:V}.

To show	\eqref{p:proplr:XII}, we take derivatives on both sides of (\ref{eq:lr}) with respect to $\rho$
\begin{align}\label{eq:p:proplr:I}
\D g_\alpha\partial_\rho h_\rho(x)=\partial_\rho\tau(\rho,\alpha,x)\cdot(\partial_\rho h_0)\circ h_{\tau(\rho,\alpha,x)}\circ g_\alpha(x)\, .
\end{align}
We have 
\begin{align}\label{eq:p:proplr:II}
(\partial_\rho h_0\circ g_\alpha)^*&(\D g_{\alpha}\partial_\rho h_0)=(\partial_\rho h_0\circ g_\alpha)^*((g_{\alpha})_\ast\partial_\rho h_0)=(g_{\alpha})^\ast(\partial_\rho h_0\circ g_\alpha)^*(\partial_\rho h_0)\nonumber\\
&={\det (\D g_{\alpha}{}|_{E_-})^\ast}(\partial_\rho h_0)^*(\partial_\rho h_0)={\det \D g_{\alpha}{}|_{E_-}}(\partial_\rho h_0)^*(\partial_\rho h_0)\, .
\end{align}
Setting $\rho=0$ in \eqref{eq:p:proplr:I}, we obtain \eqref{p:proplr:XII}, using \eqref{eq:p:proplr:II} and the non-singularity of the horocycle flow.
That $\partial_{\rho}\tau(0,\alpha,x)>0$ follows  from \eqref{p:proplr:XII}, for instance.
 
Next, since  $h_\rho$ is non-singular, the stable manifold $W_y:=h_{[0,1]}(y)$ has length bounded from above and below uniformly in $y\in M$.  Using \eqref{p:proplr:XII}  and\footnote{There is a typo in \cite[Section C]{Giulietti_2013}
and the set  $W$ there should actually be unstable.} \cite[Lemma~C.3, Remark~C.4]{Giulietti_2013}  (recalling that $g_\alpha$ is transitive), we find  $C_3,C_4, C_5<\infty$
such that
\begin{align*}
\tau(\rho,-\alpha,x)&\le C_3\int_0^\rho {\det \D g_{-\alpha|E_-}}\circ h_{\rho}(x)\dd\rho\le C_4\rho \sup_{y \in M} \int {\det \D g_{-\alpha|E_-}}\dd W_y\\
&\le C_5\rho \sup_{y \in M} \voll (g_{-\alpha}(W_y))\le C_6\rho e^{\htop\alpha}\, , \,\,
\forall \rho\ge 1\, ,\, \alpha\ge 0\, ,\, x\in M\, .
\end{align*}
A lower bound for $\tau(\rho,-\alpha,x)$ is obtained analogously, using  \cite[Lemma~C.1]{Giulietti_2013}. This shows \eqref{p:proplr:IX} for $\rho\ge  1$.
We get \eqref{p:proplr:IX} for all $\rho\le - 1$ since \eqref{p:proplr:XI} implies $\tau(-\rho,\alpha,x)=-\tau(\rho,\alpha,h_{-\rho}(x))$.  

Finally,  \eqref{p:proplr:VIII} follows from \eqref{p:proplr:XII} and the following consequence of the first claim of \eqref{p:proplr:II} 
\[
\rho=\tau(\tau(\rho,\alpha,x),-\alpha,g_{\alpha}(x))=\tau(|\tau(\rho,\alpha,x)|,-\alpha,g_{\alpha}(x)) \, .
\]
\end{proof}

\subsection{Main Result: Asymptotic Expansion for the Horocycle Integral (Theorem~\ref{t:decomposition})}\label{55b}

We need some notation to state our main result: Denote by $(X+V)'$ the dual of $X+V$ (acting on the dual of
$W_p^{s,t,q}(M)$). Recall that $\sigma((X+V)')|_{({W_p^{s,t,q}(M)})'}=\sigma (X+V)|_{W_p^{s,t,q}(M)}$ (by \cite[Section II.2.5]{Engel_2006},   strong continuity of the dual semigroup is not needed for this).
Therefore, by Corollary~\ref{c:res:spectrum}, each $\lambda\in\sigma (X+V)|_{W_p^{s,t,q}(M)}$ with $\Re\lambda>\lambda_{\min}^{s,t,p}$ is an eigenvalue of finite geometric multiplicity
 $n_\lambda$ and finite algebraic multiplicities $m_{\lambda,i}$, $1\le i\le n_\lambda$,
of  $(X+V)'$, with generalised eigenstates
$\OO_{(\lambda,i,j)}$ in the domain of  $(X+V)'$, for $ 1\le j\le m_{\lambda,i}$,
 with 
\begin{equation}\label{defOO}
((X+V)'-\lambda)^j\OO_{(\lambda,i,j)}=0\, , \quad
((X+V)'-\lambda)^{j-1}\OO_{(\lambda,i,j)}
\ne 0\,,\,\,1\le j\le m_{\lambda,i}\, .
\end{equation}

We may now state our main theorem:
 
\begin{theorem}[An Expansion for Horocycle Integrals]\label{t:decomposition}
Let $r>2$ and let $h_\rho$ be a $C^r$ reparametrisation of the unit speed horocycle flow of  a topologically mixing $C^r$ Anosov flow $g_\alpha$,  such that $d_-=1$, with $E_-$  orientable and\footnote{Hence 
$\partial_\rho\tau(0,-\alpha,\cdot)\in C^{r-1}(M)$ for all $\alpha\ge 0$.} $C^{r-1}$.  
Assume that there exist $t-(r-1)<s<0 <q<t<r-2$ and
$p\in(d/\min\{t,r-1+s\},\infty)$   with   
$$\lambda^{s,t,p}_{\min}< \htop\, , 
\mbox{ and Condition~\ref{cnd:wA} holds for some $\delta$ with }
\lambda_{\min}^{s,t,p}< \delta < \htop\, .
$$
Then  $
\Sigma_\delta:=\sigma (X+V)|_{W_p^{s,t,q}(M)} \cap 
\{\lambda\in\complex\mid \Re \lambda >\max\{0, \delta\}\}
$ is finite, and there exists $T_0>1$ such that for each $(\lambda ,i,j)$ with $\lambda\in \Sigma_\delta$, 
$1\le i \le n_\lambda$, and
$1\le j \le m_{\lambda,i}$,  there are\footnote{Note that \eqref{defc} gives a formula
for $c_{(\lambda,i,j)}$ using the generalised eigenvector of $(\lambda ,i,j)$.
We do not show $\inf_{T>T_0 ,\, x} |c_{(\lambda,i,j)}(T,x)|>0$.} functions
$$
c_{(\lambda,i,j)}:(T_0,\infty)\times M \to \complex \, ,\mbox { with } \sup_{T>T_0 ,\, x\in M} |c_{(\lambda,i,j)}(T,x)|<\infty \, ,
$$
and  for any $\tilde \delta >\max\{0, \delta\}$ there exists  
 $C_{\tilde \delta}<\infty$ such that for all $\varphi\in C^{r}(M)$ and all $T\ge T_0$ 
\begin{align*}
\int_0^T&\varphi\circ h_{\rho}(x)\dd\rho=
		T\mu(\varphi)+{\EE}_{T,x,\Sigma_\delta}(\varphi)
+\sum_{\substack{\lambda\in \Sigma_\delta \\\lambda \ne \htop}}
\sum_{\substack {1\le i \le n_\lambda\\ 1\le j\le m_{\lambda,i}}}
T^{\frac{\lambda}{\htop}}(\log T)^{j-1}\,
c_{(\lambda,i,j)}(T,x)\cdot
\OO_{(\lambda,i,j)}(\varphi)\, ,
\end{align*}
where 
\begin{align}\label{boundEE}
\sup_{x\in M}|{\EE}_{T,x,\Sigma_\delta}(\varphi)|
\le C_{\tilde \delta} (T^{\tilde \delta/\htop}\|\varphi\|_{C^r}+\|\varphi\|_{C^0}) \, .
\end{align}
\end{theorem}

The proof of Theorem~\ref{t:decomposition}  is given at the end of Section~\ref{s:lastt}.
We record an immediate corollary:

\begin{corollary}[Power Law Convergence]
\label{c:cor0}
Under the assumptions of Theorem~\ref{t:decomposition},
  there exist
$\epsilon\in (0, \min\{1,1-\delta/\htop\})$ and $C_\epsilon<\infty$ such that for all $\varphi\in{C^{r}(M)}$ 
\[|
\frac 1T\int_0^T\varphi\circ h_{\rho}(x)\dd\rho
	-\mu(\varphi)|\le \frac {C_\epsilon} {T^{\epsilon}} \|\varphi\|_{C^{r}}+\frac{C_\epsilon}{T} \|\varphi\|_{C^0}\, ,
\,\, \forall T > 0\, .
\]
\end{corollary}

A {\it contact form} is a $1$-form $\upsilon\in T^\ast M$ such that $\upsilon\wedge(\wedge_{n=1}^{\frac{d-1}2}\dd\upsilon)$ vanishes nowhere, where $\dd$ denotes the exterior derivative. (A contact form can only exist if $d$ is odd.) 
A flow $g_\alpha$ on $M$ is a {\it contact flow} if there exists a $C^1$ contact form $\upsilon$ which is preserved by the pullback of $g_{\alpha}$.
Geodesic flows on (the unit tangent bundle of) negatively curved compact manifolds  are examples of an Anosov contact flow.
Contact Anosov flows are  topologically mixing \cite[Thm 3.6]{Katok_1994}. 

\begin{proposition}\label{p:Anosov}
Let $g_\alpha$ be a $C^3$ contact Anosov flow on a compact manifold $M$ of dimension $d=3$.  Assume the strong-stable distribution $E_-$ is orientable. Then for any $\epsilon_1 >0$, 
we may choose $r \in (2, 3)$ and $t-(r-1)<s<0<t<r-2$ such that $E_-$ is $C^{r-1}$ and such that,
for any $C^{r}$ reparametrisation of the unit speed horocycle flow, we have
  $\lambda_{\min}^{s,t,p}<\epsilon_1$ for all $p\in (1,\infty)$.
In addition, if the flow satisfies the bunching condition \eqref{smalleta}, there exists  $p_0>1$ such that   Condition~\ref{cnd:wA} holds for some $\delta'\in
[\max\{\lambda_{\min}^{s,t,p},0\}, \htop)$ if $p> p_0$.
\end{proposition}

The proposition above is proved in~\S~\ref{finalS}.
Its assumptions  hold  if $g_\alpha$ is the geodesic
flow on a $C^3$ surface of strictly negative curvature 
($E_-$ is $C^{2-\tilde \eta}$ for any $\tilde \eta\in (0,1)$ 
by \cite[Thm 3.1]{Hurder_1990}; for the orientability of $E_-$, see
\cite[Lemma~B.1]{Giulietti_2013}) with $\hat\varpi$ satisfying \eqref{smalleta}. In particular, they hold
if $g_\alpha$ is the geodesic
flow on a  $C^3$ compact surface of constant negative curvature, where $\hat\varpi=2$.

We compare our main theorem and Proposition~\ref{p:Anosov} with
the results of Flaminio and Forni \cite{Flaminio_2003}:
Let $M$ be the unit tangent bundle of a compact 
surface
of constant negative curvature, and let $g_\alpha$ be its unit speed geodesic flow. Then the canonical volume form $\voll $ on $M$
is  the measure of maximal entropy for $g_\alpha$.  The unit speed horocycle flow leaves $\voll$ invariant as well, so that $\mu=\voll $. Also, the vector fields $X$ and $V=\htop$  are constant, $r=\infty$, and $\htop=1$ because $\tau(\rho,\alpha,x)=\rho\exp(-\alpha)$. In this setting \cite[Thm~1.5]{Flaminio_2003} the noninteger
obstructions to convergence, corresponding to our eigenvalues $\lambda$, are connected to
the nonzero eigenvalues $\sigma$ of the Laplacian via $\lambda=\frac 1 2 \pm \sqrt{\frac 1 4 -\sigma}$.
The   eigenvalue $\htop=1$ is simple, there are no other eigenvalues of real part equal to one, all eigenvalues with $\Re \lambda>1/2$ are semi-simple,  
and there are only finitely many eigenvalues with $\Re \lambda>\frac 12$. Moreover, since $r=\infty$, for any $p\in (1,\infty)$
(including  $p=2$)
 the  parameters $-s$, $t$ can be taken large enough to ensure $\lambda_{\min}^{s,t,p}=\lambda_{\min}^{s,t}<0$.
 Since Condition~\ref{cnd:wA} holds for some $\delta>\frac 12$  (see Proposition~\ref{p:Anosov}), we find $c_{(\lambda,i,1)}$ and ${\EE}_{T,x,\Sigma_\delta}$  as in Theorem~\ref{t:decomposition} such that
\begin{align*}
\int_0^T\varphi\circ h_{\rho}(x)\dd\rho=&T\voll (\varphi)+
\sum_{\lambda \in \Sigma_{\delta}\setminus \{1\}}
\sum_{i=1}^{n_\lambda}
T^\lambda c_{(\lambda,i,1)}(T,x)\OO_{(\lambda,i,1)}(\varphi)+{\EE}_{T,x,\Sigma_\delta}
(\varphi)\, .
\end{align*}

\subsection{Localised Horocycle Integrals, Properties of the Renormalisation Time $\tau$}

In view of the  smooth cutoff  decomposition of $\gamma_x(\cdot, T)$
in Lemma~\ref{l:horo:decomp:local} below, we introduce
localised horocycle integrals as follows.
For any  bounded compactly
supported  $w: \real\to \real$ let
\begin{align}\label{eq:cocycle:current2}
\gamma_{w,x}(\varphi)\coloneqq\int_{\real} w(\rho)\cdot(\varphi\circ h_\rho(x))\dd\rho\, ,\,\, x\in M \, , \, \, \varphi\in{C^{0}(M)}\, .
\end{align}

To show Theorem~\ref{t:decomposition}, it will be  useful\footnote{The corresponding result  for the anisotropic norms of \cite{giul_liv_2017} is slightly more intuitive.} to view $\gamma_{w,x}$ as an element of the dual  of $W_p^{s,t,q}(M)$:

\begin{lemma}\label{l:bound:dual}
There exists $\bar C<\infty$, depending on 
$\max _{\alpha \in [0,\alpha_0]} \|\phi_{-\alpha}\circ g_{-\alpha}\|_{C^{r-1}}$, 
the partition of unity $\vartheta_{\omega}$, and the charts $\kappa_{\omega}$  (Definition~\ref{d:bs:ani}), such that for any 
 $w\in C_0^{r-1}(\real)$ and $\varphi\in W_p^{s,t,q'}(M)$
\[
\sup_{x\in M}|\gamma_{w,x}(\varphi)|\le \bar C|\supp{w}|\cdot \|w\|_{C^{|s|}}\cdot \|\varphi\|_{W_p^{s,t,q'}}\, , \,\, \forall p\in(1,\infty)\,,  \, \forall t-(r-1)<s<0<q'\le t\, .\]
\end{lemma}

Before proving it,  we  show an easy consequence
of Lemma~\ref{l:bound:dual}:
The unique $h_\rho$ invariant measure $\mu$  belongs to the
dual space of $W_p^{s,t,q}(M)$.

\begin{corollary}\label{l:invmeas:II} 
Let $p\in (1,\infty)$ and  let $s$, $q$, $t$ be as in \eqref{regulS-sqt}. If
  $\phi_\alpha$ from \eqref{eq:weight} is $C^{r-1}$
then $\mu\in(W_p^{s,t,q}(M))'$.
Also, by Lemma \ref{wow}, we have
$\lambda_{\max}^{s,t,q,p}\ge \htop\in\sigma (X+V)|_{W_p^{s,t,q}(M)}$.
\end{corollary}

If $\lambda_{\min}^{s,t,q}<\lambda_{\max}^{s,t,q,p}$,
the statement of this corollary  could alternatively be obtained from \cite[\S4]{Giulietti_2013}, see the proof
of Lemma~\ref{l:invmeas}.

\begin{proof}
Fix $\epsilon>0$ small (much smaller than the diameter of $M$).
For $x\in M$  denote by $C^{r-1}_{x,\epsilon}(M)$
the set of $\varphi\in{C^{r-1}(M)}$ which vanish in an $\epsilon$ neighbourhood of $x$.
Then there exist $\delta(\epsilon)>0$ and $C(\epsilon)$  such that
for any $T>1$ with $d(h_T(x),x)<\delta$ there exists $w^{T,\epsilon}\in C^r_0(\real,[0,1])$
with $|\supp(w^{T,\epsilon})|\le T+2$ and $\|w^{T,\epsilon}\|_{C^r}\le C(\epsilon)$, such that
$$
1_{[0,T]}(\rho)\varphi(h_\rho(x))=w^{T,\epsilon}(\rho)\varphi(h_\rho(x))\, ,
\forall \varphi\in{C^{r-1}_{x,\epsilon}(M)}\,,\,\,\forall \rho\in \real\, .
$$

For any $x\in M$, since $h_\rho(x)$ is dense, there is a sequence $T_n=T_n(x,\epsilon)$ such that
$d(h_{T_n}(x),x)<n^{-1}$ and $T_n\to \infty$.
By unique ergodicity \eqref{thm2.1} and Lemma~ \ref{l:bound:dual},  we have
\begin{align}\label{eq:l:invmeas:II}
|\mu(\varphi)|\le \lim_{n\to\infty}\bigl |\frac 1{T_n}\gamma_{x}(\varphi,T_n) \bigr|= \lim_{n\to\infty}\bigl |\frac 1{T_n}\gamma_{w^{T_n,\epsilon},x}(\varphi) \bigr |
&\le 2 C(\epsilon) \bar C\|\varphi\|_{{W_p^{s,t,q}}}\, , \,\forall \varphi\in{C^{r-1}_{x,\epsilon}(M)}\, .
\end{align}
Next, using a $C^\infty$ function $\psi=\psi_{x,y,\epsilon}:M\to[0,1]$, vanishing
in an $\epsilon$ neighbourhood of $x$, and $\equiv 1$ in an $\epsilon$ neighbourhood of
some  $y\ne x$, we can write any $\varphi\in C^{r-1}(M)$  as
$\psi \varphi+(1-\psi)\varphi$, where $\psi \varphi\in{C^{r-1}_{x,\epsilon}(M)}$,
and $(1-\psi )\varphi\in{C^{r-1}_{y,\epsilon}(M)}$.
Applying \eqref{eq:l:invmeas:II} at $x$ and $y$ gives
\begin{align*}
|\mu(\varphi)|
&\le 2 C(\epsilon) \bar C\bigl ( \|\psi \varphi\|_{\widetilde W_p^{s,t,q}}
+ \|(1-\psi) \varphi\|_{\widetilde W_p^{s,t,q}}\bigr ) \, , \,\forall \varphi\in{C^{r-1}(M)}\, ,
\end{align*}
where $\widetilde W_p^{s,t,q}$ is defined like $W_p^{s,t,q}$, but using
systems of cones $\widetilde \Theta_\omega$ (see Remark~\ref{forLeibniz}) ensuring that
$\|\psi \varphi\|_{\widetilde W_p^{s,t,q}}\le C \|\psi\|_{C^r} \|\varphi\|_{W_p^{s,t,q}}$
for some $C<\infty$ and all $\varphi$, $\psi$.
We conclude by density of $C^{r-1}$ functions in $W_p^{s,t,q}(M)$.
\end{proof}

\begin{proof}[Proof of Lemma~\ref{l:bound:dual}]
Let $\delta_*$ denote the Dirac distribution, fix
$w\in C_0^{r-1}(\real)$, and  set
\begin{align*}
w_{x,\phi,\omega_1,\omega_2,\alpha}(z)&:=
(\vartheta_{\omega_2}\cdot\phi_{-\alpha}\circ g_{-\alpha})\circ \kappa_{\omega_1}^{-1}(z)\cdot\int_{-\infty}^\infty w( \rho)\delta_*(z-\kappa_{\omega_1}\circ g_{\alpha}\circ h_{ \rho}(x))\dd \rho\, ,\\
\varphi_{\omega,\alpha}(z)&:=(\vartheta_{\omega}\cdot\LL_{\alpha,V}\varphi)\circ\kappa_{\omega}^{-1}(z)\,,\,\,  z\in \real^d\,,\,\,x\in M\, ,\, \alpha \ge 0\, ,\, \omega\, ,\, \omega_{1}\,, \omega_2\in \Omega\, , \varphi\in{C^{r-1}(M)}
\, .\nonumber
\end{align*}
Since $h_\rho$ has no periodic orbits and $M$ is compact,  $w_{x,\phi,\omega_1,\omega_2,\alpha}$ is  a  bounded function 
supported in the interior of a subset $J$ of the (one-dimensional) stable leaf at $g_{\alpha}(x)$ (using \eqref{eq:lr}) in charts. In addition,  there exists $\bar C_0$ 
such that  
 $$
|J|\le\bar C_0 |\supp w|\quad \mbox{ and } \quad
\sup_{\alpha \in [0,\alpha_0], \, x\in M, \, \omega_i\in \Omega}  
\|w_{x,\phi,\omega_1,\omega_2,\alpha}\|_{L_\infty} \le \bar C_0 \|w\|_{L_\infty}\, .
$$

Next, since 
$\phi_{-\alpha}\circ g_{-\alpha} = 1 /{\phi_\alpha}$,
 we have, exchanging the integrals with respect to $z$ and $\rho$
(so that $z=\kappa_{\omega_1}\circ g_\alpha\circ h_\rho(x)$)
\begin{align*}
\gamma_{w,x}(\varphi)=\sum_{\omega_1,\omega_2\in \Omega}\int_{\real^d}w_{x,\phi,\omega_1,\omega_2,\alpha}(z)
\cdot\varphi_{\omega_1,\alpha}(z)\dd z\, ,\,\, \forall \alpha \ge 0 \, .
\end{align*}
Recalling $\widetilde{\Psi'}_{\sigma,n}$ from \eqref{eq:suppext}, we find, using 
Plancherel's theorem for the inner product of two functions, $\widetilde{\Psi'}_{\sigma,n}{\Psi}_{\sigma,n}={\Psi}_{\sigma,n}$, Cauchy--Schwarz for the sum in
$\sigma$ and $n$,   H\"older's inequality,  
a constant $C<\infty$ such\footnote{The bounds below can be viewed as yet another avatar
of integration by parts.} that, for all $p\in(1,\infty)$ and all $t-(r-1)<s<0<q'\le t$, 
\begin{align}
\alpha_0\cdot |\gamma_{w,x}(\varphi)|&=\int_0^{\alpha_0}
|\sum_{\omega_1,\omega_2}
\int_{\real^d}w_{x,\phi,\omega_1,\omega_2,\alpha}(z)\cdot\varphi_{\omega_1,\alpha}(z)\dd z|\dd\alpha\nonumber\\
		\le &\int_0^{\alpha_0}\sum_{\omega_1,\omega_2}
|\int_{\real^d}\sum_{\sigma,n} 2^{-c(\sigma)n}{\widetilde{\Psi'}}^\Op_{\sigma,n}
(w_{x,\phi,\omega_1,\omega_2,\alpha})(z) 2^{c(\sigma)n}{\Psi}^\Op_{\sigma,n}
(\varphi_{\omega_1,\alpha})(z)\dd z|\dd\alpha\nonumber\\
		\label{eq:bound:dual:pre:II}&\le C\sup_{\alpha \in[0,\alpha_0]}\sum_{\omega_1,\omega_2}
\|(\sum_{\sigma,n}4^{-c(\sigma)n}
|{\widetilde{\Psi'}}^\Op_{\sigma,n}
(w_{x,\phi,\omega_1,\omega_2,\alpha})|^2)^{\frac 12}\|_{L_{1-1/p}}  \cdot \|\varphi\|_{{W_p^{s,t,q'}}}\, ,
\end{align}
using the definition \eqref{eq:norm} of the norm in \eqref{eq:bound:dual:pre:II}. To conclude, it suffices to find $C_0<\infty$ such that   
		\begin{equation}
\label{l:bound:horo:first:I}
\max_{\substack{\sigma\in\{+, 0\}\\\omega_1,\omega_2\in\Omega}}  \,\,\, \sup_{x \in M}
\|{\widetilde{\Psi'}}^\Op_{\sigma,n}
(w_{x,\phi,\omega_1,\omega_2,\alpha})\|_{L_{1-1/p}}\le C_0 |\supp w|\cdot \|w\|_{L_\infty}, \,
\,  \forall 0\le \alpha\le \alpha_0\, , \,  \forall n\in \naturall \, ,
\end{equation}
and (since $c(+)>0$ and $c(0)>0$, it is enough to consider $\sigma=-$)
\begin{equation}\label{l:bound:horo:first:II}
\max_{\omega_1,\omega_2\in\Omega}\sup_{x \in M}\|{\widetilde{\Psi'}}^\Op_{-,n}
(w_{x,\phi,\omega_1,\omega_2,\alpha})\|_{L_{1-1/p}}\le \frac{C_0}{ 2^{(r-1)n}}|\supp w|\cdot \|w\|_{C^{|s|}}\, , \,
\forall 0\le \alpha\le \alpha_0\, ,   \forall n\in \naturall\, .
\end{equation}
Young's inequality for $\|\FFF^{-1}(\widetilde{\Psi'}_{\sigma,n})*
w_{x,\phi,\omega_1,\omega_2,\alpha}\|_
{L_{1-1/p}}$  gives \eqref{l:bound:horo:first:I} for $C_0=\bar C_0 \max \diamm V_{\omega}$.

Finally, we show  \eqref{l:bound:horo:first:II}. There are 
$\widetilde C_0$ (depending on
$\max _{\alpha \in [0,\alpha_0]} \|\phi_{-\alpha}\circ g_{-\alpha}\|_{C^{r-1}}$, 
 $\vartheta_{\omega}$,  and $\kappa_{\omega}$), a subset $\tilde J\subset \real$, and a $C^{r-1}$ diffeomorphism $\tilde y:\tilde J\to J\subset\real^d$ with\footnote{As a warmup, the reader is invited to think
of the case when $J$ is a subset of a coordinate axis in $\real^d$.}
$|\tilde J|\le C|J|\le \widetilde C_0|\supp w|$,  with $\max\{\|\tilde y\|_{C^{r-1}},\|\tilde y^{- 1}\|_{C^{r-1}}\}\le \widetilde C_0$, 
and 
\begin{align*}
{\widetilde{\Psi'}}^\Op_{-,n}w_{x,\phi,\omega_1,\omega_2,\alpha}(z)&=
\frac {1}{(2\pi)^d}
\int_{\real^d}\int_{\real^d} 
\widetilde{\Psi'}_{-,n}(\xi)e^{\iunit \xi (z-\tilde z)}
w_{x,\phi,\omega_1,\omega_2,\alpha}(\tilde z)\dd\xi\dd \tilde z\\
&=
\frac {1}{(2\pi)^d}
\int_{\real^d}\int_{\tilde J} 
\widetilde{\Psi'}_{-,n}(\xi)e^{\iunit \xi (z-\tilde y(y))}
w_{x,\phi,\omega_1,\omega_2,\alpha}(\tilde y(y))\dd\xi\dd y\, ,
\end{align*}
with
 $w_{x,\phi,\omega_1,\omega_2,\alpha} \circ \tilde y$ a $C^{r-1}$ function supported in the interior of $\tilde J$ such that
$$
\sup_{\alpha \in [0,\alpha_0], \, x\in M, \, \omega_i\in \Omega}  
\|w_{x,\phi,\omega_1,\omega_2,\alpha} \circ \tilde y\|_{C^{\tilde r}(\tilde J)} \le \widetilde C_0 \|w\|_{C^{\tilde r}}\, ,\,\,
\forall \tilde r\le r-1 \, .
$$
Note that $J$  lies in a stable cone in charts. Thus,
there exists $C_1>0$ such that $ |\partial_{y} (\xi \tilde{y}(y))|\ge C_1 2^n$ for any $\xi$ in
the support of $\widetilde{\Psi'}_{-,n}$ (which lies inside $E^*_-$
 in charts).  Finally, integrating ${\lfloor|r-1|\rfloor}$ times  by parts  with respect to $y$,  following by a regularised integration by parts  if $|r-1|$ is not an integer (Lemmas~\ref{l:partint} and~\ref{l:partint:reg}),  and  ending with Young's inequality, we get  \eqref{l:bound:horo:first:II}.
\end{proof}

The next two lemmas use the following  version of the  renormalisation equation
\eqref{thekey} for the localised horocycle integral \eqref{eq:cocycle:current2}.

\begin{sublemma}[Renormalisation and Smooth Localisation]\label{sublemmarenorm} Fix  $x\in M$
and $\varphi\in C^{0}$, then 
\begin{align*}
\gamma_{w,x}(\varphi)&=\int_\real w( \tau(\rho,-\alpha,g_{\alpha}(x)))
\cdot \LL_{\alpha,V} \varphi(h_\rho (g_{\alpha}(x))) \dd \rho \, ,\, \forall \alpha \ge 0 \, .
\end{align*}
\end{sublemma}

\begin{proof} 
By definition and our choice $\phi_\alpha=\partial_{\rho}\tau(0,-\alpha,\cdot)$,
$$ 
\int_\real w( \tau(\rho,-\alpha,g_{\alpha}(x)))
\cdot \LL_{\alpha,V} \varphi(h_\rho (g_{\alpha}(x))) \dd \rho=\gamma_{w(\tau(\cdot,-\alpha,g_{\alpha}(x)),g_{\alpha}(x)}
(\LL_{\alpha,\partial_{\rho}\tau(0,-\alpha,\cdot)} \varphi)\, .
$$
Thus, the sublemma follows from  \eqref{eq:lr} and the first claims of \eqref{p:proplr:II} and \eqref{p:proplr:V}, since
\begin{align}
\gamma_{w,x}(\varphi)
&=\int_{-\infty}^\infty w(\rho)\cdot\varphi\circ g_{-\alpha}\circ h_{\tau(\rho,\alpha,x)}\circ g_{\alpha}(x)\dd\rho\nonumber\\
&=\int_{-\infty}^\infty  w(\tau(\rho,-\alpha,g_{\alpha}(x)))\cdot\varphi\circ g_{-\alpha}\circ h_{\rho}\circ g_{\alpha}(x)\cdot\partial_{\rho}\tau(\rho,-\alpha,g_\alpha(x))\dd\rho\nonumber\\
&=\int_{-\infty}^\infty  w(\tau(\rho,-\alpha,g_{\alpha}(x)))\cdot (\partial_{\rho}\tau(0,-\alpha,\cdot)\cdot \varphi\circ g_{-\alpha})\circ h_{\rho}\circ g_{\alpha}(x)\dd\rho\nonumber\\
\nonumber \label{eq:horo:renorm:II}
&=\gamma_{w\circ \tau(\cdot,-\alpha,g_{\alpha}(x)),g_{\alpha}(x)}
(\partial_{\rho}\tau(0,-\alpha,\cdot)\cdot \varphi\circ g_{-\alpha})
\, ,\,\,\forall\alpha\ge 0\, .
\end{align}
\end{proof}

Taking $w=w_T$ to be the characteristic function $w_T=1_{[0,T]}$, we have $\gamma_{w_T,x}(\varphi)=
\gamma_x(\varphi,T)$, and by choosing $\alpha=O(\log T)$ we can ensure (in view of 
Lemma~\ref{l:existlr})
that the support of $w_T\circ \tau(\cdot,-\alpha,g_{\alpha}(x))$ has size $O(1)$,
uniformly in $x$. In order to apply Lemma~\ref{l:bound:dual}, 
some regularity of $w$ is required: we thus need a more clever choice of localisation function $w_T=w_{T,x}$. We state the corresponding result,  similar
to \cite[Lemma 5.16]{Flaminio_2003}, \cite[Lemma 3.19]{giul_liv_2017}.

\begin{lemma}[Bounds for Localised Horocycle Integrals]\label{l:horo:decomp:local} 
Let $C> 1$ be as in \eqref{p:proplr:IX}--\eqref{p:proplr:VIII} and fix $\bar C>\max\{C,4\}$. 
If $\phi_\alpha$ from \eqref{eq:weight} is $C^{r-1}$, then 
for every $T> C \bar C$ and $x\in M$ there exists a compactly supported $C^{r}$ function
 $w=w_{T,x}:\real\to [0,1]$ with
\begin{equation}\label{l:horo:decomp:local:I} 
|\gamma_{x}
(\varphi, T)-\gamma_{w_{T,x},x}(\varphi)|\le 2 C \bar C \|\varphi\|_{C^0}\, ,\,\,\forall
\varphi\in  C^0(M)\, .
\end{equation}
Moreover, for $p\in (1,\infty)$ and $t-(r-1)<s<0<q'\le q\le t$,  there exists 
$\widetilde C<\infty$ 
such that, if $\tilde \varphi \in W_p^{s,t,q}(M)$ 
 satisfies
\begin{align}\label{eq:l:horo:decomp:local:bound}
\|\LL_{\alpha,V}\tilde{\varphi}\|_{{W_p^{s,t,q'}}}
\le  \exp( \alpha \cdot a)\max\{1,|\alpha|^{j-1}\} C_{\tilde \varphi}
\,,\,\,\forall \alpha \ge 0\, ,
\end{align}
for some\footnote{The proof gives \eqref{l:horo:decomp:local:IIa}, replacing
$\widetilde C \cdot C(a)\cdot T^{\frac a{\htop}}(\log T)^{j-1}$ by
$\widetilde C\cdot (\log T)^{j}$ if $a=0$, by
$\widetilde C \cdot (\log T)^{j-1}$ if $a<0$.}  $a>0$,  $j\ge 1$, and $C_{\tilde \varphi}<\infty$,  then, setting $C(a)=1/(1-(C/\bar C)^{a/\htop})$, we have
\begin{equation}\label{l:horo:decomp:local:IIa}
\sup_{x\in M}|\gamma_{w_{T,x},x}(\tilde{\varphi})|
\le \widetilde C \cdot C(a) 
T^{\frac a{\htop}}(\log T)^{j-1} 
 C_{\tilde \varphi} \, ,\quad \forall  T>C \bar C\, .
\end{equation}
\end{lemma}

\begin{proof}
For $x\in M$ and $T>C \bar C$, define inductively sequences $\alpha^\pm_k=\alpha^\pm_k(x,T)\in\real$,  $k\ge 1$,  by
\begin{align*}
	&\bar C=\tau(T,\alpha^+_1,x)\, ,\,
1=\tau(\tau(\bar C,-\alpha^+_{k},g_{\alpha^+_{k}}(x)),\alpha^+_{k-1},x)=
\tau(\bar C,\alpha^+_{k-1}-\alpha^+_{k},g_{\alpha^+_{k}}(x))\, ,\nonumber\\
	&	\alpha^-_1= \alpha^+_1\, ,
\,\,-1=\tau(\tau(-\bar C,-\alpha^-_{k},g_{\alpha^-_{k}}( h_T(x))),\alpha^-_{k-1},h_T(x))
=\tau(-\bar C,\alpha^+_{k-1}-\alpha^+_{k},g_{\alpha^-_{k}}(h_T(x)) )\, ,
\end{align*}
where  we used the first claim of \eqref{p:proplr:II}. In the special case when $\tau(\rho,\alpha,x)=\rho e^{-\alpha \htop}$ we find 
$$\alpha^-_1= \alpha^+_1=\frac{\log(T/\bar C)}{\htop}>0\,,\qquad \alpha_{k}^+-\alpha_{k-1}^+=\alpha_{k}^--\alpha_{k-1}^-=\frac{\log(1/\bar C)}{\htop}<0\,,\,\,k\ge 2\, .
$$
More generally, 
since $\tau(T,0,x)=T$
and $\tau(T,\alpha,x)$ is continuous in $\alpha$, the bounds
 \eqref{p:proplr:IX}--\eqref{p:proplr:VIII} give 
$0<\log (T/(\bar CC))\le \htop \cdot \alpha^\pm_1\le \log (TC/\bar C)$.
It is also easy to check that
$\alpha^+_{k}<\alpha^+_{k-1}$ and $\alpha^-_{k}<\alpha^-_{k-1}$
for all $k\ge 2$, and that \eqref{p:proplr:IX} gives 
\begin{align}\label{eq:l:horo:decomp:local:I}
\{e^{\htop(\alpha_{k}^--\alpha_{k-1}^-)}, e^{\htop(\alpha_{k}^+-\alpha_{k-1}^+)}\}\in \bigl [\frac 1{C \bar C},  \frac C {\bar C}\bigr ]\, ,
\, \, \forall k\ge 2 \, ,
\,  \forall x \in M\, , \forall T>  C\bar C\, .
\end{align}
Thus, we have
\begin{align}\label{eq:l:horo:decomp:local:0}
 T /(C\bar C)^{k} \le e^{\htop\alpha_k^\pm}\le T C^{k}/\bar C^k\, , \quad
\forall k \ge 1 \, .
	\end{align}
Fixing   a $C^\infty$ function $\chi:\real\to [0,1]$ such that
$\chi|_{[ 1,\infty)}\equiv 1$ and
$\chi|_{(-\infty, 0]}\equiv 0$,  we put
\[
w_1(\rho)= \chi(\tau(\rho,\alpha^+_1,x))\cdot \chi(- \tau(\rho-T,\alpha^-_1,h_T(x)))\, ,
\]
and, for $k\ge 2$ (note that $w_1$  and the  $w^\pm_k$ depend on $x$ and $T$), 
\begin{align*}
w_k^+(\rho)&= \chi(\tau(\rho,\alpha^+_{k},x))-\chi(\tau(\rho,\alpha^+_{k-1},x))\, ,\\
w_k^-(\rho)&= \chi (-\tau(\rho-T,\alpha^-_{k},h_T(x)))
-\chi(- \tau(\rho-T,\alpha^-_{k-1}, h_T(x)))\, .
\end{align*}
Then, for any $n_\pm\ge 1$, we have
\begin{align}\label{nN}
w_1(\rho)+ \sum_{k=2}^{n_+} w^+_k(\rho)+\sum_{k=2}^{n_-}w_k^-(\rho)=
&w_1(\rho)+\chi (\tau(\rho,\alpha^+_{n},x))-\chi( \tau(\rho,\alpha^+_1,x))\\
\nonumber &+\chi( -\tau(\rho-T,\alpha^-_{n},h_T(x))-\chi(-\tau(\rho-T,\alpha^-_1, h_T(x)))\, ,
\end{align}
so that  \eqref{p:proplr:VIII} and  \eqref{eq:l:horo:decomp:local:0} imply
$w_1+\sum_{k=1}^\infty( w^+_k+w_k^-)=1|_{(0,T)}$.
Define $N_\pm\ge 1$ by
  $\min\{\alpha_{N_+}^+,\alpha_{N_-}^-\}\ge 0$
and $\max\{\alpha_{N_++1}^+,\alpha_{N_-+1}^-\}\le 0$ 
(note that $N_\pm\in[ 1, \log (T)/\log (\bar C/C)]$ by \eqref{eq:l:horo:decomp:local:0}). Then put
 $$w_{T,x}(\rho)=w_1(\rho)+\sum_{k=2}^{N_+} w^+_k(\rho)+\sum_{k=2}^{N_-}w_k^-(\rho)\, .
$$
Thus,  setting $n_\pm=N_\pm$ in \eqref{nN}, the definitions of $\chi$ and $\alpha_{N+1}^\pm$ give
\begin{align}
\label{thus}
\supp(1|_{{[0,T]}}-w_{T,x})\subset&[0,\tau(\bar C,-\alpha_{N_++1}^+,g_{\alpha_{N_++1}^+}(x))]\\
\nonumber &\qquad \cup[T+\tau(-\bar C,-\alpha_{N_-+1}^-,g_{\alpha_{N_-+1}^-}( h_T(x)),T) ]\,  . 
\end{align} 
Using  $\max\{\alpha_{N_++1}^+,\alpha_{N_-+1}^-\}\le 0$, the claim \eqref{l:horo:decomp:local:I} then follows from \eqref{p:proplr:VIII}.

Next, for $\tilde \varphi\in W_p^{s,t,q}(M)$, using $\gamma_{v,x}(\tilde \varphi)=\gamma_{v\circ(\cdot+T),h_T(x)}(\tilde \varphi)$, Sublemma~\ref{sublemmarenorm} gives\footnote{This is analogous to the decomposition in \cite[Lemma 3.1]{giul_liv_2017}. We use   more explicit  smoothing functions.}  
\begin{align}\label{keydec}
\gamma_{w_{T,x},x}(\tilde \varphi)=&\gamma_{\tilde{w}_1,g_{\alpha_1^+}(x)}
(\LL_{\alpha_1^+,V}\tilde \varphi)
+\sum_{k=2}^{N_+}
 \gamma_{\tilde{w}_k^+,g_{\alpha^+_{k}(x)}}
(\LL_{\alpha_{k}^+,V}\tilde \varphi)
+\sum_{k=2}^{N_-}\gamma_{\tilde{w}_k^-,g_{\alpha^-_{k}}(h_T(x))}
(\LL_{\alpha_{k}^-,V}\tilde \varphi)\, \, ,
\end{align}
where, recalling \eqref{p:proplr:II}, we put $\tilde{w}_1= w_1( \tau(\cdot,-\alpha^+_1,g_{\alpha^+_1}(x)))
= \chi\cdot \chi(\bar C-\cdot)$, and
\begin{align*}
\tilde{w}_k^+(\rho)=& w_k^+(\tau(\rho,-\alpha^+_{k},g_{\alpha^+_{k}}(x)))
\, ,\,\,
\tilde{w}_k^-(\rho)= w_k^-(T+\tau(\cdot,-\alpha^-_{k},g_{\alpha^-_{k}}(h_T(x)))\,,\,\, k\ge 2\, .	
\end{align*}
Since  $\tilde{w}_k^+(\rho)
=\chi(\rho)-\chi( \tau(\rho,\alpha^+_{k-1}-\alpha^+_{k},g_{\alpha^+_{k}}(x)))$ and also  $\tilde{w}_k^-(\rho)
=\chi(-\rho)-\chi(-\tau(\rho,\alpha^-_{k-1}-\alpha^-_{k},g_{\alpha^-_{k}}( h_T(x))))$, we find
\begin{align}\label{eq:l:horo:decomp:local:support}
\supp\tilde{w}_1
\subseteq\bigl [0, C \bar C\bigr ]\,,\,\, \supp\tilde{w}_k^+\cup -\supp\tilde{w}_k^-
\subseteq\bigl [0, C \bar C\bigr ]\,,\,\,  \forall k\in\naturall\, .
\end{align}
Since $\phi_\alpha\in C^{r-1}$, \eqref{p:proplr:XI}--\eqref{p:proplr:XII} imply 
$\sup_{\alpha \ge 0, x \in M}\|\partial_\rho\tau(\cdot,\alpha,x)\|_{ C^{r-1}}<\infty$,  so
 \eqref{eq:l:horo:decomp:local:I}
 gives
\begin{align}\label{eq:l:horo:decomp:local:III:a}
\sup_{x \in M} \max\{\|\tilde{w}_1\|_{C^{r}},
\sup_{k\ge 2}\, 	 \, \sup_{T>\bar C C} \|\tilde{w}_k^\pm\|_{C^{r}}
\} <\infty  \, .
\end{align}
Thus, if \eqref{eq:l:horo:decomp:local:bound}
holds for $\tilde \varphi$, applying Lemma~\ref{l:bound:dual}, and   \eqref{eq:l:horo:decomp:local:0} to \eqref{keydec},  we find   $\hat C<\infty$ such that
\begin{align}\label{eq:l:horo:decomp:local:V}
\sup_x |\gamma_{w_{T,x},x}(\tilde \varphi)|\le \hat C C_{\tilde \varphi}
T^{\frac a{\htop}} (\log T)^{j-1}\,\sum_{k=1}^{\max \{N_-,N_+\}}  (C/\bar C )^{k \frac a \htop}
 \, ,
\end{align}	
for all $T > C\bar C$. Since $4<C<\bar C$, summing  
$\sum_{k=1}^\infty ( C/{\bar C} )^{k\frac a{\htop}}$ gives \eqref{l:horo:decomp:local:IIa}.
\end{proof}

\subsection{Exact Bounds ($\sup_{\alpha\ge 0} \|e^{-\alpha \htop}\LL_{\alpha,V}\|_{W^{s,t,q}_p}<\infty$). Proof of Theorem~\ref{t:decomposition}}\label{s:lastt}

We saw in Lemma~\ref{wow} that $\mu$, the  unique  $h_\rho$ invariant probability,
is a fixed point
of   $e^{-\htop\alpha}\LL_{\alpha,V}'$ acting on Radon measures,
in Remark \ref{Rwow} that $\sup_{\alpha \ge 0}\|e^{-\htop\alpha}\LL_{\alpha, V}\|_{L^1(\mu)}\le 1$,
and in Corollary~\ref{l:invmeas:II} that $\mu\in (W_p^{s,t,q}(M))'$ so that $\lambda_{\max}^{s,t,q,p}\ge\htop$.
If  $\lambda_{\min}^{s,t,p}<\htop$,
we  get more.\footnote{In fact,  only  the exact growth claim is needed
from Lemma \ref{l:invmeas}:
The rest of the information about the peripheral spectrum could be obtained by 
an ad hoc argument based on  \eqref{l:invmeas:I}, the identity in the proof of
Lemma~\ref{wow}, \eqref{l:horo:decomp:local:IIa}, and   Sublemma~\ref{sublemmarenorm}.
See \cite[Lemma 5.18 (v), and last claim of Lemma 5.14]{adam}.}  

\begin{lemma}[Peripheral Spectrum and Exact Growth]\label{l:invmeas}
If $E_-$ is $C^{r-1}$ and $\lambda_{\min}^{s,t,p}<\htop$
for some $p\in(1,\infty)$ and  $s$, $q$, $t$ as in \eqref{regulS-sqt}, then 
for all $0<q\le t$ we have
$\lambda_{\max}^{s,t,q,p}= \htop$.
Moreover,  $\htop$ is a simple 
eigenvalue and the only element of $\{ \lambda\in \sigma(X+V)|_{W_p^{s,t,q}(M)}\,,\,
\Re \lambda= \htop\}$.
In particular, there are no maximal Jordan blocks, and\footnote{This exact growth estimate
is a key ingredient, e.g.,  for \cite[Assumption~1]{Butterley_2016}
used in the proof of Theorem~\ref{t:decomposition}. See e.g. the inverse Laplace transform in \cite[Lemma~4.3]{Butterley_2016}.} 
 $\sup_{\alpha\ge 0}\| e^{-\alpha \htop}\LL_{\alpha,V}\|_{W^{s,t,q}_p}<\infty$.
\end{lemma}

For the potential $V$  associated to the SRB measure, 
where $\lambda^{s,t,q,p}_{\max}=0$, the results above are well known,  see
\cite[Lemma~5.1]{Butterley_2007} and \cite{Butterley_2013} (the claims
there are for other Banach spaces, but
intrinsicness can be applied as in our proof of Lemma~\ref{l:invmeas}).

\begin{remark}[MME and Bypassing Unique Ergodicity]\label{MME}
Exploiting the results of \cite{Giulietti_2013} as in the proof
of Lemma~\ref{l:invmeas}, it can be shown (without using Corollary~\ref{l:invmeas:II}), that the unique fixed point of $e^{-\htop \alpha}\LL_{\alpha, V}'$ in the dual of $W^{s,t,q}_p(M)$
is a Radon measure $\mu$, and letting $\nu\in W^{s,t,q}_p(M)$
be the unique fixed point of $e^{-\htop \alpha}\LL_{\alpha, V}$, that the distribution formally defined by
$\mu_{*}(\varphi)=\mu( \varphi\nu)$
is a
Radon measure, and it is the unique measure of maximal entropy (MME) of $g_\alpha$, which 
is (exponentially) mixing. (See e.g. \cite{Gouezel_2008} for the discrete-time analogue.)
In fact, unique ergodicity of $h_\rho$ (that is \eqref{thm2.1}) could be obtained from the information
 on the peripheral spectrum of $\LL_{\alpha, V}$, bypassing the results of Bowen and Marcus from \cite{Bowen_1977}. To keep the paper short, we refer to \cite{Bowen_1977}.
\end{remark}
Before showing  Lemma~\ref{l:invmeas}
we state and prove consequences of the exact
growth. 

\begin{corollary}[Exact Growth for the Resolvent]\label{c:exact}
Assume that $E_-$ is $C^{r-1}$ and $\lambda_{\min}^{s,t,p}<\htop$
for $p\in(1,\infty)$ and  $s$, $q$, $t$ as in \eqref{regulS-sqt}. Fix $0<q\le t$. 
There exists $C<\infty$ such that 
\begin{align*}
\|\RR_{z}^{n}\varphi\|_{{W_p^{s,t,q}}}&\le
\frac{ C}{(\Re z-\htop)^n}\|\varphi\|_{{W_p^{s,t,q}}}\, ,\, \forall \Re z>\htop\,, 
\,\, \forall n\ge 1\, .
\end{align*}
Moreover, recalling \eqref{truncres}, there exist $C<\infty$ and a  system $ \Theta=\{ \Theta'_\omega\}$ with $ \Theta'_\omega<\Theta_\omega$
such that 
$$
\|\RR^n_{z}\varphi-\RR^n_{tr, z}\varphi \|_{W_{p}^{s,t,q}}
\le\frac C  {(\Re z-\htop)^{n}}\|\varphi\|_{W_{p, \Theta'}^{s,t,q}}
\, ,\,\,\forall \Re z>\htop\, , \, 
\forall n \ge 1\, .
 $$ 
\end{corollary} 

\begin{proof}
The first bound follows from  exact growth
($\sup_{\alpha\ge 0}\| e^{-\alpha \htop}\LL_{\alpha,V}\|_{W^{s,t,q}_p}<\infty$),
simplifying \eqref{eq:res:bound:III}.
The second claim follows from exact growth,
using Remark~\ref{forLeibniz} (for $\alpha\ge\alpha_0$) with  \eqref{eq:res:int}.
\end{proof}

\begin{proof}[Proof of Lemma~\ref{l:invmeas}]
By Corollary~\ref{c:res:spectrum}, since  $\lambda_{\min}^{s,t,p}<\htop$,  we claim that we can  exploit
Theorem~\ref{intr} about
intrinsicness to transfer\footnote{An independent proof should exist.
To  exclude maximal Jordan blocks, a geometric argument is needed see \cite[\S4.3]{Giulietti_2013}.
 Maybe \eqref{p:proplr:IX}--\eqref{p:proplr:VIII} can help, see \cite[Lemma 5.17]{adam}.} the results of Giulietti, Liverani, and Pollicott in \cite{Giulietti_2013} to our spaces.

Indeed, first recall  that a $C^{r-1}$ one-form is a $C^{r-1}$ section of the cotangent bundle $T^* M$, or equivalently, a $C^{r-1}$ map from the tangent space $TM$ to $\mathbb {R}$ whose restriction to each fibre $T_xM$ is a linear functional on $T_xM$.
Using that the Anosov flow $g_\alpha$ is topologically mixing (and $E_-$ is orientable),
they showed  \cite[(4.5), Lemma 4.7, Proposition 4.9, for $\ell=d_-=1$]{Giulietti_2013} that $\htop$ 
(denoted $\sigma_{d_s}$ there, with $d_-s=d_-$) is a simple eigenvalue
and the only  element $\lambda$ of the spectrum with $\Re \lambda \ge \htop$ for the generator
$Y^{(d_-)}$
of the pullback semigroup $\LL_\alpha^{(d_-)}$ of $g_{-\alpha}$ \cite[(2.9)]{Giulietti_2013}
acting on the closure   $\widetilde \BB^{1,|s|,d_-}$
of $C^{r-1}$ one-forms on $M$ vanishing in the flow direction, 
for an anisotropic Banach norm
(see   \cite[Def. 3.6 and (4.6)]{Giulietti_2013} for $\ell=d_-=1$,
$p=1$, and $q=t$, and note that this is equivalent to letting the pullback semigroup
act on the Grassmannian of line bundles in $TM$ as in \cite{Gouezel_2008,giul_liv_2017}).
A key step for this   is the fact that, setting 
 $\tilde \lambda_{\min}^{1,|s|}:= \htop +\min\{1,|s|\}\log\theta <\htop$, the intersection of the spectrum
of $Y^{(d_-)}$  on $\widetilde \BB^{1,|s|,d_-}$
with the half-plane $\Re \lambda >\tilde \lambda_{\min}^{1,|s|}$ contains only isolated eigenvalues
of finite multiplicity (this is shown  by establishing the
corresponding result for the resolvent $\RR^{(d_-)}_z$ \cite[Def.~4.4, Lemma~4.8]{Giulietti_2013}).

Next, recall that, by our assumptions,
$r\ge 2$ and $E_-$ is $C^{r-1}$ (so that $E_-^*$ is $C^{r-1}$ too).
The (closed) subspace  $\overline \BB^{1,|s|}$ 
of $\widetilde \BB^{1,|s|,d_-}$ obtained by taking the closure
(for the norm of $\widetilde \BB^{1,|s|,d_-}$) of the space  $\Omega^{r-1}_{E_-}$
of those $C^{r-1}$ one-forms taking values in $E_-^*$,
is invariant under $\LL_\alpha^{(d_-)}$.
Using the natural bijection 
$\varphi \mapsto (\varphi(\cdot), E^*_-(\cdot))$ from $C^{r-1}(M)$ to  $\Omega^{r-1}_{E_-}$, 
we see that
the restriction of $\LL_\alpha^{(d_-)}$
to  $\Omega^{r-1}_{E_-}$  coincides with our  operator $\LL_{\alpha,V}$ on $C^{r-1}(M)$.
 It is well-known\footnote{For example, the spectrum of the
two-sided shift restricted to sequences vanishing on one side 
is the whole disc, while the original spectrum is the circle.} that restricting  a bounded operator
$\RR$ to a closed invariant subspace $\overline \BB \subset \widetilde \BB$  can 
fill up the holes (a hole in a compact set of $\complex$ is a bounded connected component of
its complement)
in the original spectrum, but the spectrum of the restriction does not
intersect the unbounded connected component of 
$\sigma(\RR|_{\widetilde\BB})$ (see \cite[Corollary~4.1]{Scro}). Hence,
the intersection of the spectrum
of $Y^{(d_-)}|_{\overline \BB^{1,|s|}}$
with the half-plane $\Re \lambda >\tilde \lambda_{\min}^{1,|s|}$ still contains only isolated eigenvalues
of finite multiplicity.
Some of the eigenvalues of $Y^{(d_-)}$ on ${\widetilde \BB^{1,|s|,d_-}}$
can disappear for the restricted operator, but we already established in Corollary~\ref{l:invmeas:II}  that $\htop$ is an eigenvalue in our space.
Finally, since $\max \{\lambda_{\min}^{s,t,p}, \tilde \lambda_{\min}^{1,|s|}\}<\htop$, we can apply 
Theorem~\ref{intr}, using also
that $C^{r-1}$ functions are dense in  $W^{s,t,q}_p$ and in $\overline \BB^{1,|s|}$, 
and that $W^{s,t,q}_p$ and $\overline \BB^{1,|s|}$  are both continuously embedded into
the dual of  $C^{|s|+\max\{1,t\}}$ \cite[Lemma 3.10]{Giulietti_2013}.
\end{proof}

\begin{proof}[Proof of Theorem~\ref{t:decomposition}]
The starting point of the proof is 
$$\gamma_x(1,T)\mu(\varphi)=\int_0^T \mu(\varphi) d\rho=T\mu(\varphi)$$
(this is trivial for the unit speed horocycle flow, an easy computation otherwise). Thus, we may and shall assume
$\mu(\varphi)=0$, replacing $\varphi$ by $\varphi -\mu(\varphi)$ (constants belong to our Banach space).

Fix $0< q\le t$. 
By Lemma~\ref{l:invmeas}, $\lambda^{s,t,p}_{\max}=\htop$, it is a simple eigenvalue and the only maximal eigenvalue of $X+V$ on $W_p^{s,t,q}(M)$ . Hence, $\OO_{\htop}:=\OO_{\htop,1,1}=\mu$, so that  $\OO_{\htop}(\varphi)=0$. 

For a general Ruelle--Pollicott resonance $\lambda\in\sigma (X+V)|_{W_p^{s,t,q}(M)}$ with $\Re\lambda>\lambda_{\min}^{s,t,p}$, recalling the
notation introduced above \eqref{defOO},  we denote
by  $\DD_{(\lambda,i,j)}\in D(X+V)$, for $1\le i \le n_\lambda$ and $1\le j\le m_{\lambda,i}$,
its  generalised eigenstates, i.e.,
\[
(X+V-\lambda)^j\DD_{(\lambda,i,j)}=0\, , \,\,\,\,
(X+V-\lambda)^{j-1}\DD_{(\lambda,i,j)}
\ne 0\, .
\]
We  write $\DD_{\htop}:=\DD_{\htop,1,1}$.
There is a curve
$\Gamma_\lambda$ around $\lambda$ with 
$\frac{1}{2 \iunit \pi} \oint_{\Gamma_\lambda}(z-(X+V))^{-1}\varphi \, \dd z=\sum_{i=1}^{n_\lambda} \Pi_{\lambda,i}\varphi$, 
where $\Pi_{\lambda,i}$  is a projector  of rank $m_{\lambda,i}$, with
\begin{align*}
	\Pi_{\lambda,i}=\sum_{j=1}^{m_{\lambda,i}}\DD_{(\lambda,i,j)}\otimes
\OO_{(\lambda,i,j)}\, , \, 1\le i\le n_\lambda
\, ,\, \,\,\OO_{(\lambda,i,j)}\in D((X+V)')\, ,
\end{align*}
where   
$
\OO_{(\lambda_1,i_1,j_1)}(\DD_{(\lambda_2,i_2,j_2)})=1$ if $(\lambda_1,i_1,j_1)=(\lambda_2,i_2,j_2)$ and 
$\OO_{(\lambda_1,i_1,j_1)}(\DD_{(\lambda_2,i_2,j_2)})=0$, otherwise.
In addition, there are finite rank
 nilpotent operators $\NN_{\lambda,i}$, for $1\le i\le n_\lambda$, such that
\begin{align*}
	&\Pi_{\lambda_1,i_1}\Pi_{\lambda_2,i_2}\equiv 0\mbox { and } \NN_{\lambda_1,i_1}\NN_{\lambda_2,i_2}\equiv 0\,  \mbox{ if }\lambda_1\ne \lambda_2\mbox{ or } i_1\ne i_2 \, ,
\\
\nonumber
&\NN_{\lambda,i}^{m_{\lambda,i}-1}\equiv 0\, ,\qquad\,	\Pi_{\lambda_1,i_1}\NN_{\lambda_2,i_2}=\NN_{\lambda_2,i_2}\Pi_{\lambda_1,i_1}=
\begin{cases}
\NN_{\lambda_2,i_2}&\mbox{if } \lambda_1=\lambda_2\mbox{ and }i_1=i_2\\0&\mbox{if } \lambda_1\ne \lambda_2\mbox{ or }i_1\ne i_2\, ,
\end{cases} 
\end{align*}
and,  using the surjection from eigenvalues of $X+V$  to those of the semi-group \cite[V (2.3)]{Engel_2006}
 $$\LL_{\alpha,V}\Pi_{\lambda,i}=\exp(\alpha \lambda)
\exp(\alpha \NN_{\lambda,i})\Pi_{\lambda,i}\,,\,\, \forall \alpha\ge 0\, .
$$
(See also e.g. \cite{Butterley_2016}.)
Therefore,   for each $(\lambda,i,j)$ there exists $C_{i,j}<\infty$ such that
	\begin{equation}\label{that}
\|\LL_{\alpha,V}\DD_{(\lambda,i,j)}
\|_{W_p^{s,t,q}}
\le C_{i,j} \exp(\alpha \cdot \Re \lambda )\max\{1,|\alpha|^{j-1}\}
\|\DD_{(\lambda,i,j)}\|_{W_p^{s,t,q}}\, ,\,
\forall \alpha \in \real
\, .\end{equation}
In other words, $\DD_{(\lambda,i,j)}$ satisfies  \eqref{eq:l:horo:decomp:local:bound} for all $\alpha\in\real$ if $a=\Re \lambda>
 \lambda_{\min}^{s,t,p}
$.

Assume that $\delta >0$, let 
$$\lambda\in\Sigma_\delta=\sigma (X+V)|_{W_p^{s,t,q}(M)}\cap \{z\in\complex\mid\Re z\ge \delta\}$$
 and fix $x\in M$.
Let  $T\ge T_0=\bar C C>1$ and $w_{T,x}\in C^{r-1}_0$ be given by Lemma~\ref{l:horo:decomp:local},
and
define
\begin{equation}\label{defc}
c_{(\lambda,i,j)}(T,x)
:= T^{-\frac{\lambda}{\htop}}
(\log T)^{1-j}
\gamma_{w_{T,x},x}(\DD_{(\lambda,i,j)})
\in \complex\,,\,\,
1\le i \le n_\lambda\, ,\,\,1\le j\le m_{\lambda,i}\, .
\end{equation} 
Then \eqref{l:horo:decomp:local:IIa} from Lemma~\ref{l:horo:decomp:local}
implies that   $\sup_{x, T\ge T_0}|c_{(\lambda,i,j)}(T,x)|<\infty$. 
Decomposing
\begin{align*}
\gamma_{w_{T,x},x}(\Pi_{\lambda,i}\varphi)
&=\sum_{j=1}^{m_{\lambda,i}}\OO_{(\lambda,i,j)}(\varphi) 
\gamma_{w_{T,x},x}(\DD_{(\lambda,i,j)})
=\sum_{j=1}^{m_{\lambda,i}}c_{(\lambda,i,j)}
T^{\frac{\lambda}{\htop}}
(\log T)^{j-1}\OO_{(\lambda,i,j)}(\varphi)\, ,
\end{align*}
and using Lemma~\ref{l:invmeas},  we find for any finite subset $\Lambda_\delta\subset \Sigma_\delta$ that
\[
\gamma_x(\varphi, T)
=\gamma_{w_{T,x},x}(\DD_{\htop})\mu(\varphi)+
\sum_{\substack{\lambda\in \Lambda_\delta\\
\lambda\ne\htop}}
\sum_{i=1}^{n_\lambda}
\sum_{j=1}^{m_{\lambda,i}} 
c_{(\lambda,i,j)}T^{\frac{\lambda}{\htop}}
(\log T)^{j-1}
\OO_{(\lambda,i,j)}(\varphi)
+{\EE}_{T,x,\Lambda_\delta}(\varphi)\, ,
\]
where $\gamma_{w_{T,x},x}(\DD_{\htop})\mu(\varphi)=0$,
and 
\begin{align*}
{\EE}_{T,x,\Lambda_\delta}(\varphi)=& 
\gamma_{w_{T,x},x}\bigl (\varphi 
-\sum_{ \lambda\in \Lambda_\delta}
		\sum_{i=1}^{n_\lambda}\Pi_{\lambda,i} \varphi\bigr )
+\gamma_{x}(\varphi,T)-\gamma_{w_{T,x},x}(\varphi)\,  .
\end{align*}

To conclude, we show that finiteness of
$\Sigma_\delta$ and the claimed bound on  ${\EE}_{T,x,\Sigma_\delta}$  follow from  Condition~\ref{cnd:wA}. 
We first check that  Assumptions 1, 2, and 3A from \cite{Butterley_2016} hold for the semigroup $e^{-\htop\alpha}\LL_{\alpha,V}$ on $W_p^{s,t,q}(M)$ (with generator $X+V-\htop$, resolvent $\RR_{z+\htop}$):
Note that $\RR_{\htop}=(\htop-(X+V))^{-1}$ is bounded on the codimension one subset $W(\htop)$
of $W^{s,t,q}_p$ formed of those $\tilde \varphi$ such that $\mu(\tilde \varphi)=0$. Therefore,
the norm on $W(\htop)$ defined by
$$
\|\tilde \varphi\|_{weak}=\frac
{\|(\htop-(X+V))^{-1}(\tilde \varphi)\|_{W^{s,t,q}_p}}
{\|\RR_{\htop}\|}
$$ 
satisfies  $\|\tilde \varphi\|_{weak}\le \|\tilde \varphi\|_{W^{s,t,q}_p}$. The  identity
$\tilde \varphi-e^{-\htop\alpha}\LL_{\alpha,V}\tilde \varphi
=(\htop-(X+V))\int_0^{\alpha} e^{-\htop\tilde \alpha}\LL_{\tilde \alpha,V}\tilde \varphi \dd \tilde \alpha
$
thus implies  Assumption~1 in \cite{Butterley_2016}, i.e.
\begin{align}
\label{eq:AS1}
	\sup_{\alpha>0}\frac 1\alpha
\|\id -e^{-\htop\alpha}\LL_{\alpha,V}\|_{W_p^{s,t,q}(M)
\rightarrow{weak}} <\infty\, ,\, \forall \tilde \varphi \in W(\htop)\, .
\end{align}
(Indeed, it is enough to consider $\alpha \in (0,1]$ in
\eqref{eq:AS1} due to the exact growth.)
Since $\htop-\lambda_{\min}^{s,t,p}>0$, 
the essential spectral radius of $\RR_{z+\htop}$ is not larger than   $|\Re z +\htop-\lambda_{\min}^{s,t,p}|^{-1}$ by Corollary~\ref{c:res:spectrum}, giving Assumption~2 in \cite{Butterley_2016}. 
Finally,  since $p>d/\min\{t,r-1+s\}$,  Proposition~\ref{cnd:A} and Corollary~\ref{c:exact} imply that Condition~\ref{cnd:wA}
gives \eqref{cnd:sA}, i.e., Assumption~3A from
 \cite{Butterley_2016} for   $\RR_{z+\htop}$.

Thus \cite[Thm~1]{Butterley_2016} gives   $\#\Sigma_\delta<\infty$ and furnishes, for 
  $\tilde \delta>\delta$, a constant $C_B(\tilde \delta)<\infty$ with 
\begin{align}\label{this}
\bigl \|e^{-\htop \alpha}\LL_{\alpha,V}\bigl (\psi-\sum_{ \lambda\in \Sigma_\delta}\sum_{i=1}^{n_\lambda}
\Pi_{\lambda,i}\psi\bigr )
\bigr\|_{weak}
&\le C_Be^{ \tilde \delta \alpha}
\|(\htop-(X+V))\psi\|_{W_p^{s,t,q}}\, , \,\, \forall\alpha\ge 0\, ,
	\end{align}
for all $\psi \in W(\htop)$. 

Finally, Lemma~\ref{l:horo:decomp:local} gives the bound \eqref{boundEE} for $\sup_x |{\EE}_{T,x,\Sigma_\delta}(\varphi)|$:
First,  \eqref{l:horo:decomp:local:I}  implies that
$\sup_{x,T}|\gamma_{x}(\varphi,T)-\gamma_{w_{T,x},x}(\varphi)|\le C \|\varphi\|_{C^0}$. Second,  
setting $\psi_j=(\htop-(X+V))^j\varphi$, for $j=1,2$,  we have
 $\max_{j=1,2}\|\psi_j\|_{W_p^{s,t,q}}\le C \|\varphi\|_{C^r}$, 
 because  $V\in C^{r-1}$, and $t<r-2$. 
Then  \eqref{l:horo:decomp:local:IIa} (with 
$\tilde \varphi=\psi_1\in W_p^{s,t,q}(M)$ and
$C_{\tilde \varphi}=C_B\|(\htop-(X+V))\tilde \varphi\|_{W_p^{s,t,q}}$)  gives $T_0<\infty$ such that
$$
\sup_x |
\gamma_{w_{T,x},x}
(\varphi-\sum_{ \lambda\in \Sigma_\delta}\sum_{i=1}^{n_\lambda}\Pi_{\lambda,i}\varphi)|
\le  \widetilde C C(\tilde \delta) C_B(\tilde \delta) \cdot T^{\tilde \delta /\htop}\|(\htop-(X+V)) \psi_1\|_{W_p^{s,t,q}}\, ,\,\forall T \ge T_0 \, .
$$
\end{proof}

\subsection{Proof of Proposition~\ref{p:Anosov}}
\label{finalS}

We assumed the flow fixes  a $C^1$  contact $1$-form $\upsilon\in T^\ast M$.  In particular, $\upsilon$ is annihilated on $E_++E_-$ and the volume  in $\wedge^3 T^\ast M$ is preserved by the flow.  Then, since $d=3$, we already mentioned that \cite[Thm 3.1]{Hurder_1990} 
gives that   $E_-$ is $C^{2-\tilde \eta}$ for any $\tilde \eta\in (0,1)$. Taking $r=3-\tilde \eta$,  we find $\partial_\rho \tau(0,-\alpha,\cdot)\in C^{r-1}$ for any  $C^{r}$ reparametrisation of the unit speed  horocycle flow.

It follows from \eqref{p:proplr:XII} 
and Lemma~\ref{l:to:loc:bound} that  the  transfer operators associated to $C^{r}$ reparametrisations  are conjugate to each other and thus have the same spectrum (using Remark~\ref{forLeibniz}). For the unit speed parametrisation we have
 $\phi_\alpha=\partial_\rho\tau(0,-\alpha,0)=\det \D g_{-\alpha}|_{E_-}$. We claim that
\begin{align}\label{gdbd}
\lambda_{\min}^{s,t,p}(X,V)=
\lim_{\alpha\to\infty}\frac 1\alpha\log
\|\phi_\alpha
|\det(\D g_{\alpha})^{\tr}{}|_{E_-^\ast}|^{\min\{t,-s\}}\, \|_{L_\infty(M)}\, .
\end{align}
Indeed, $d_-=d_+=1$, and
since  the flow $g_\alpha$ preserves volume, i.e.  $|\det\D g_{\alpha}|\equiv 1$, we find
\[
|\det(\D g_{-\alpha})^{\tr}|_{E_+^\ast}|=
|\det(\D _{g_{-\alpha}}g_{\alpha})^{\tr}|_{E_0^\ast}|
|\det(\D _{g_{-\alpha}}g_{\alpha})^{\tr}|_{E_-^\ast}|
\, .\]
Using \eqref{eq:expconst} and the upper and lower bounds on
$|\det(\D _{g_{-\alpha}}g_{\alpha})^{\tr}{}|_{E_0^\ast}|$, we get
\eqref{gdbd}.

Then,  taking 
$
-s=t=\frac {r-1}2-\frac{\tilde \eta} 2=1-\tilde \eta 
$,  we have $t-r+1\le s<0<t<r-2$, and formula \eqref{gdbd}  together with  \eqref{p:proplr:XII} 
give that
$\lambda^{s,t,p}_{\min}<\epsilon_1$, if $\tilde \eta>0$ is small enough.

It remains to discuss Condition~\ref{cnd:wA}. 
(This condition is stable under  reparametrizations of the horocycle flows, using
\eqref{eq:res:int} and the conjugacy mentioned in the previous paragraph.)
Since we assumed \eqref{smalleta}, the second\footnote{The proof of \cite[Proposition 7.5]{Giulietti_2013}  has a gap
since a factor $e^{z\tau_W\circ H_{\beta,i,W}}$ is missing from 
\cite[(7.14)]{Giulietti_2013}. However, the statement is correct \cite[Theorem 1 and its proof]{Livv},
replacing the condition $\min \{1,\hat\varpi\}>2/3$
in \cite{Giulietti_2013} by: ``$\lambda_+-\lambda_-< \vartheta_0 \lambda_-$
  with $\vartheta_0\in(0, \frac{1}{4})$
 and  $\varpi\ge 1$,'' which hold since we assumed \eqref{smalleta}.} claim of  \cite[Proposition 7.5]{Giulietti_2013}  holds for $\ell=d_-$.
Therefore,  since our operator $\RR_z$ coincides with the operator denoted $\RR^{(d_s)}(z)$ 
in \cite{Giulietti_2013}
restricted to those one-forms which take their images in $E_-^*$ (as in the proof of Lemma~\ref{l:invmeas}),  the second claim of  \cite[Proposition~7.5]{Giulietti_2013} and \cite[Lemma~7.4]{Giulietti_2013}
combined with \cite[(7.1)]{Giulietti_2013} (which holds due to the exact growth bounds)
and \eqref{triviall} (for $\delta_1=0$, $\delta_2>0$ and $\beta=3\gamma_0$)
 gives 
$\eta \in (0,1)$, 
  $a_0\ge 1$, $b'_0\ge 1$, $\delta_2\in (0,\htop)$,  $C<\infty$, $\gamma_0\in (0,1)$, $C_1>1$,  and,
for any $a\ge a_0$ and
$\gamma'\ge a C_1$
such that\footnote{We may take  $\delta_2>0$  small enough 
in \eqref{triviall} to ensure $aC_1<3\gamma_0/\log \bigl (1+\delta_2/a)$.}
 $\gamma'< 3\gamma_0/\log \bigl (1+\delta_2/a\bigr )$, we have
\begin{equation}\label{wrongnorm}
\|\RR_{a+\iunit b+\htop}^{n}\varphi \|_{\bar \BB^{0,1+\eta}}
\le C|a+\delta_2|^{-n}\|\varphi \|_{\bar \BB^{1,\eta}}\, ,\,\,
\forall |b|\ge b'_0\, , 
\mbox{ where }n=\lceil\gamma'  \log |b|\rceil\, ,
\end{equation}
for the anisotropic Banach
spaces $\bar \BB^{j,\eta}$, $j=0,1$, in the scale
from the  proof of Lemma~\ref{l:invmeas}. 
(The statement of \cite[Proposition~7.5]{Giulietti_2013} is for $a\in [a_0,2a_0]$ and 
$\gamma'\ge C_1$, the proof gives \eqref{wrongnorm}.)

By
\cite[Lemma 3.10]{Giulietti_2013} the space $\bar \BB^{1,\eta}$ 
lies in the dual of 
 $C^{1+\eta}(M)$. Thus, if $s''<-1-\eta -3-\frac{3}{p}$,
we have $\|\tilde \varphi\|_{W_p^{s'',s'',s''}}\le C\|\tilde \varphi\|_{\bar \BB^{0,1+\eta}}$. 
(Use Sobolev embeddings \cite[Thm 2.2.3(i)]{RS} for
the dual $B^{-1-\eta}_{1,1}$  of  $b_{\infty,\infty}^{1+\eta}\supset C^{1+\eta}(M)$  \cite[Def. 2.1.3.1(ii), Remark 2.1.5.1]{RS} and
$W^{s''}_p=F^{s''}_{p,2}$ in dimension $d=3$.)
The last bound of \cite[Remark 3.8]{Giulietti_2013}  gives
$\|\tilde \varphi \|_{\bar \BB^{1,\eta}}\le C\|\tilde \varphi \|_{C^1}
$, ending the proof.

\begin{appendix}

\section{Integration by Parts}\label{intparts}

\begin{lemma}[Integration by Parts (cf. text after {\cite[Remark 3.3]{Baladi2007}})]\label{l:partint}
Let $f\colon\real^d\rightarrow\real$ be $C^1$ and 
compactly supported.
For any $C^2$ function $ G\colon\real^d\rightarrow \real$ 
such that $\inf_{\supp f}\sum_{j=1}^d(\partial_{j} G)^2>0$, 
\begin{align*}
\int_{\real^d}e^{i  G(z)}f(z)\dd z=i\int_{\real^d}e^{i G(z)}\sum_{k=1}^d\partial_{k} \frac{(\partial_k G(z))f(z)}{\sum_{j=1}^d(\partial_{j} G(z))^2}\dd z\, .
\end{align*}
\end{lemma}

For $z\in\real^d$ we write  $\nabla_z G=(\partial_{j}G(z))_{j=1,\ldots ,d}$ for the gradient   and $\nabla_z^{\tr}G=\sum_{j=1}^d\partial_{j} G(z)$  for the divergence of $G:\real^d\to \real$.
 Let $\nu:\real^d\rightarrow\real_+$ be $C^\infty$, supported in the unit ball, and such that $\int_{\real^d}\nu(x)\dd x=1$.  Then we have:

\begin{lemma}[Regularised Integration by Parts  {\cite[(3.4)]{Baladi2007}}]\label{l:partint:reg}
Fix $1<r<2$. Let $f\colon\real^d\rightarrow\complex$ be a compactly supported $C^{r-1}$-map, let $ G\colon\real^d\rightarrow \real$ be $C^{r}$ and such that $|\nabla_z G(z)|^2=\sum_{j=1}^d(\partial_{j} G)^2>0$ on $\supp f$. Set\footnote{In particular, there exists $\bar C\ge 1$ such that $\|\nabla_z h_{\delta}\|_{L_\infty}\le \bar C \|h\|_{C^{r-1}}\delta^{r-1}$ and $\|h-h_{\delta}\|_{L_\infty}\le \bar C \|h\|_{C^{r-1}}\delta^{r-1}$.} 
$$h(z):=i\frac{\nabla_z G(z)f(z)}{|\nabla_z G(z)|^2}\,,
\quad h_\delta:= \delta^{-d} \cdot h* \nu\biggl (\frac{\cdot}{\delta}\biggr )\,, \,\,\delta >0\, .
$$
Then, for every $L\ge 1$, 
\begin{align*}
	\int_{\real^d}e^{\iunit L G(z)}f(z)\dd z=&
\frac 1 L\int_{\real^d}e^{\iunit L G(z)}\nabla_z^{\tr} h_\delta(z)\dd z
- \iunit\int_{\real^d}e^{\iunit L  G(z)}\nabla_z^{\tr} G(z)\bigl(h(z)-h_{\delta}(z)\bigr )\dd z\nonumber\, .
\end{align*}
\end{lemma}

\section{Fragmentation and Reconstitution}
\label{fragrec}
 
A finite set of $C^{r}$ functions $\vartheta_j:\real^d\to [0,1]$ such that 
$\sum_j \vartheta_j(x)\le 1$ for all $x\in \real^d$ is called a $C^{r}$ sub-partition of unity.
The fragmentation and reconstitution lemmas  of \cite{Baladi2007} and \cite{Bbook} extend  straightforwardly to our 
anisotropic spaces.
The first lemma is a variant of \cite[Lemma 7.1]{Baladi2007}:

\begin{lemma}[Fragmentation]\label{partexpdag0}
Let  $1<p<\infty$ and let $s$, $q$, $t$ as in \eqref{regul-sqt}
   be real numbers, and let $K\subset \real^d$
be compact. For any $s',t', q'\in \integer$, there exists
 $C<\infty$ such that, for
any $C^{r}$ sub-partition of unity $\{\vartheta_j\}_{j=1,\ldots J}$ 
of $K$ with intersection multiplicity
$\nu$, there exists $\tilde C_\vartheta<\infty$
such that (in the applications, we take  $s'<s$, $t'<t$, and $q'\le q$)
  \begin{equation}
\biggl \| \sum_{j=1}^J  \vartheta_j v \biggr \|_{W_{p,\Theta}^{s,t,q}}  
\le C \nu^{(p-1)/p}
  \bigl ( \sum_ j \| \vartheta_j v\|^p_{W_{p,\Theta}^{s,t,q}} \bigr )^{1/p}+ 
  \tilde C_\vartheta \sum_{j=1}^J \|\vartheta_j v\|_{W_{p,\Theta}^{ s',t',q'}} \, .
  \end{equation}
\end{lemma}

The last lemma, a variant of \cite[Prop. 7.2]{Baladi2007} (see also \cite[Lemma 4.29]{Bbook}), 
is useful to group partitions of unity
associated with a fixed cone system:

\begin{lemma}[Reconstitution]\label{partgroupdag}
Let $1<p<\infty$, let   $s$, $q$, $t$ as in \eqref{regul-sqt} be real, and let $K\subset \real^d$
be compact. If $\Theta'< \Theta$ then for any $s', q', t'\in \integer$, there exists
$C<\infty$ such that, for
any $C^{r}$ sub-partition of unity $\{\vartheta_j\}_{j=1,\ldots J}$ 
of $K$ with intersection multiplicity $\nu$, there exists  $C'_\vartheta<\infty$
such that (in the applications, we take  $s'<s$, $t'<t$, and $q'\le q$)
  \begin{equation}
\biggl ( \sum_{j=1}^J  \|\vartheta_j v\|^p_{W_{p,\Theta'}^{s,t,q}} \biggr ) ^{1/p}\le C 
\nu^{1/p}
  \| v\|_{W_{p, \Theta}^{s,t,q}} + C'_\vartheta  \|v\|_{W_{p, \Theta}^{s',t',q'}}\, .
  \end{equation}
\end{lemma}

\section{Interpolation, Mollification, and Approximations of the Identity}
\label{mollapp}

Let $[\BB_1,\BB_2]_{u}$, for $u \in (0,1)$, denote the complex (Calder\'on) interpolation of an interpolation pair of Banach spaces
$\BB_1$, $\BB_2$.  The 
Banach spaces  in this paper are complex interpolation spaces:

\begin{lemma}[Interpolation]\label{interpolation}
Setting
$w(x_1,x_2)=(1-w)x_1+w x_2$ for $w \in (0,1)$, we have for any
$p\in (1,\infty)$, all $t_j-(r-1)<s_j<0<q_j\le t_j$, $j=1, 2$, and all $w \in (0,1)$ that
$$
[W^{s_1,t_1,q_1}_p(M), W^{s_2,t_2,q_2}_p(M)]_w=W^{s,t,q}_p(M)\, ,\,\,
s=w(s_1,s_2)\, ,\,\, t=w(t_1,t_2)\, ,\,\,
q=w(q_1,q_2)\, .
$$
\end{lemma}

\begin{proof}
The norms on $\real$ given by
$\|x\|_{n,u}= 2^{nw} |x|$, $n\in \naturall$,
form a complex interpolation scale with respect to $w\in \real$. The lemma thus follows from
using \cite[Thms 1.18.1, 1.18.4]{Trr}
to show that the local norms $\| \cdot\|_{W^{s,t,q}_{\Theta_\omega,p}}$ have the 
desired interpolation property, and
then applying \cite[Thm~1.18.4]{Trr} to the function
$\alpha \mapsto 
\|(\vartheta_\omega \LL_{\alpha,V}\varphi )\circ \kappa_\omega^{-1}\|_{W^{s,t,q}_{\Theta_\omega,p}}$ 
on $[0,\alpha_0]$.
\end{proof}

Recall the finite atlas $\AAA$,  indexed by $\omega\in\Omega$, and the pair $\{\Theta_\omega\}$, $\{\Theta'_\omega\}$  of adapted cone systems
from Lemma~\ref{l:excone} and Remark~\ref{forLeibniz}.
Let $\{U_\omega\}$ be an open cover of $M$ with $\overline U_\omega\subset V_\omega$.
Let $\nu$ be as in Lemma~\ref{l:partint:reg} and set
$\nu_\epsilon(x)=\epsilon^{-d}\nu(x/\epsilon)$.  Fix $C^\infty$ functions $\bar \vartheta_{\omega}:M \to [0,1]$, with $\bar \vartheta_\omega$
supported in  $U_\omega$, 
such that 
$\sum_{\omega} \bar \vartheta_{\omega} (x)=1$ for all 
$x\in M$.
Finally, let  $\epsilon>0$ be such that the $\epsilon$-neigbourhood of
$\kappa_\omega(U_\omega)$ is contained in $\kappa_\omega(V_\omega)$ for each $\omega$.
As in \cite[(5.4)]{Baladi2012},  define a mollifier operator $\MMM_\epsilon$,  by  setting,
for any distribution $\varphi$ of order at most $r$ on $M$,  
\begin{align}
\nonumber & (\MMM_\epsilon (\varphi))_\omega(u)
= 
\int_{\real^{d}} \nu_\epsilon ( u-v) 
\psi ( \kappa_{\omega}^{-1}(v))
\, \dd v = [\nu_\epsilon *  (\psi \circ \kappa_\omega^{-1})](u)\, ,\quad \omega\in \Omega\, ,\,
u \in \kappa_\omega(U_\omega)\, , 
\\
\label{defmoll}
&\MMM_\epsilon (\varphi)=\sum_{\omega\in \Omega} \bar \vartheta_\omega \cdot (( (\MMM_\epsilon (\varphi))_\omega \circ  \kappa_\omega) \, .
\end{align}

Since $\{\Theta_\omega\}$ and  $\{\Theta'_\omega\}$ are adapted to $\AAA$ and $g_\alpha$,  the fact that $\Theta'_{\omega}<\Theta_{\omega}$ in the next lemma is not a problem:

\begin{lemma}[Approximation of the Identity]\label{moll}
For any $p\in (1,\infty)$, all  $s',t',q'\in \real$ and all $\eta>0$ such that 
$-(r-1)+t'+\eta<s'<-\eta <0<q'<t'$, 
there exists $C<\infty$ such that, letting $W^{s',t',q'}_{p}(\Theta')$
be the space constructed with $\Theta'$,
$$\|\MMM_\epsilon\varphi-\varphi\|
_{W^{s',t',q'}_{p}(\Theta')}
\le C \epsilon^ \eta
\|\varphi\|_{W^{s'+\eta,t'+\eta,q'+\eta}_{p}}\,,
\,\, 
\forall \varphi\,,\,\,\forall \epsilon >0\, .$$
\end{lemma}

\begin{proof}
Minkowski-type integral bounds
hold for the local norms $W^{s,t,q}_{p, \Theta'_\omega}$: There exists $C_M<\infty$ such that
 for any    $\psi : \real^{d}\to \real_+$ and
any  family $\varphi_{y} \in W^{s,t,q}_{p, \Theta'_\omega}$, uniformly
bounded in $y$,
\begin{align}\label{minko}
\| \int_{\real^{d}} \psi(y) \varphi_{y}(\cdot )\, dy\|_{W^{s,t,q}_{p, \Theta'_\omega}}
&\le C_M \|\psi\|_{L^1(\real^{d})}
\sup_{y} \|\varphi_{y}\|_{W^{s,t,q}_{p, \Theta'_\omega}}\, .
\end{align}
 (See \cite[Remark 5.1]{Baladi2012}. We already established interpolation for our spaces.)
So we proceed as in the proof of
 \cite[Lemma 5.4]{Baladi2012}.
The changes of charts
$\kappa_{\omega'}\circ \kappa_{\omega}^{-1}$  (one chart is for the mollifier and the other for the norm)  are cone-hyperbolic from $ \Theta_\omega$
to $\Theta'_{\omega'}$ by construction.
\end{proof}

\begin{remark}\label{C.3}
If we attempted to show  Dolgopyat bounds
using mollifiers  through isotropic spaces as in \cite[Lemma 5.4, (7.5)--(7.6)]{Baladi2012},
we would face a factor 
\begin{equation}\label{smalls'}\|\RR^n_{a+\iunit b+\htop}\|_{{W_p^{s',s',s'}}}\le \frac C
{(a+s'\log \Theta)^n}
\end{equation}
instead of $Ca ^{-n}$  in \eqref{ani}. After applying \eqref{triviall} with  $\beta=\kappa(s-s')-1>0$, we would  end up with
an upper bound $\gamma' <
\frac{\kappa(s-s') -1
}{\log \bigl (1+\frac { \lambda_{\max}-\delta}a)-\log \bigl (1-\frac{|s'|\log \Theta}{a})}$. 
In our main application, Proposition~\ref{p:Anosov}, we need to take $s'$ close to $-1$ to guarantee
$\lambda_{\min}<\lambda_{\max}$. The upper bound would then conflict with \eqref{conf}.
This is why  (proving \eqref{conj} would give another  solution to this problem) we used mollification through anisotropic spaces
as\footnote{For another use of approximations of identity
with anisotropic norms see \cite[(2.77) proof of Cor. 2.3, sentence after (3.27) in the proof
of Lemma 3.5]{GJ}.} in \cite[Lemma D.2]{Giulietti_2013}, taking advantage of the exact growth from Lemma~\ref{l:invmeas}. Our norms are different from those of \cite{Giulietti_2013}: Their drawback is that we need to go through
charts twice  to prove Lemma~\ref{moll}
(we have Remark~\ref{forLeibniz} to save us). Their strength is that we can
use the interpolation\footnote{Interpolation was not available in \cite{Baladi2012} due to the presence of
a supremum in the norm there.} Lemma~\ref{interpolation} and Minkowski inequalities
to prove Lemma~\ref{moll}.
\end{remark}
\end{appendix}

\end{document}